\documentclass[11pt,a4paper,draft,twoside]{article}

\usepackage{amsmath, exscale, epsfig,amssymb , 
}

\small
\normalsize
\usepackage[latin1]{inputenc}
\usepackage[T1]{fontenc}
\usepackage{amsfonts}
\usepackage{amsbsy}
\usepackage{fancyheadings}
\usepackage{amscd}
\usepackage{theorem}

\usepackage{epsf}
\usepackage{psfrag}
\usepackage{graphicx}

\usepackage{multicol}
\allowdisplaybreaks

\def\build#1_#2^#3{\mathrel{\mathop{\kern 0pt#1}\limits_{#2}^{#3}}}
\def\noi{{\noindent}}
\def\be{\begin{equation}}
\def\ee{\end{equation}}
\def\ba{\begin{eqnarray*}}
\def\ea{\end{eqnarray*}}

\def\cqfd{ \hfill $\blacksquare$ }
\def\cq{ \hfill $\square$ }

\def\cA{{\cal A}}
\def\cH{{\cal H}}
\def\cC{{\cal C}}

\def\cI{{\cal I}}
\def\cF{{\cal F}}
\def\cG{{\cal G}}

\def\cE{{\cal E}}

\def\ee{\epsilon}

\def\noi{\noindent}

\newcommand{\lgeo}{[\![}
\newcommand{\rgeo}{]\!]}

\def\build#1_#2^#3{\mathrel{\mathop{\kern 0pt#1}\limits_{#2}^{#3}}}

\newtheorem{theorem}{Theorem}[section]
\newtheorem{lemma}[theorem]{Lemma}
\newtheorem{proposition}[theorem]{Proposition}
\newtheorem{corollary}[theorem]{Corollary}
\newtheorem{definition}{Definition}[section]
{\theorembodyfont{\rmfamily}\newtheorem{remark}{Remark}[section]}
{\theorembodyfont{\rmfamily}\newtheorem{comment}{Comment}[section]}
{\theorembodyfont{\rmfamily}\newtheorem{example}{Example}[section]}
\newcommand{\R}{\mathbb{R}}

\newcommand{\bP}{\mathbb{P}}

\newcommand{\N}{\mathbb{N}}
\newcommand{\bN}{\mathbb{N}}
\newcommand{\U}{\mathbb{U}}

\newcommand{\un}{\boldsymbol{1}}

\begin{document}

\title{ {\bf   THE CODING OF COMPACT REAL TREES BY REAL VALUED FUNCTIONS } }
\author{ by \\
Thomas {\sc Duquesne}
 \thanks{Université Paris-Sud, Mathématiques, 91405 Orsay Cedex,
   France; email: thomas.duquesne@math.u-psud.fr} }
\vspace{4mm}
\date{\today} 

\maketitle

\begin{abstract} This paper is a detailled study of the coding of real trees 
by real valued functions that is motivated by probabilistic 
problems related to continuum random trees. Indeed it is known 
since the works of Aldous \cite{Al2} 
and Le Gall \cite{LG93} that a continuous non-negative function $h$ on $[0,1]$ such that $h(0)=0$ can be
seen as the contour process of a compact real tree. This particular coding of
a compact real tree provides additional structures, namely 
a root that is
the vertex corresponding to $0\in [0,1]$, a linear order 
inherited from the usual
order on $[0,1]$ and a measure induced by the Lebesgue measure on
$[0,1]$; of course, the root, the linear order and the measure obtained by such a coding
have to satisfy some compatibility conditions. In 
this paper, we prove that any compact real tree equipped with a root, a linear
order and a measure that are compatible can be encoded by a non-negative 
function $h$
defined on a finite interval $[0, M]$, that is assumed to be left-continuous with right-limit, 
without positive jump and such that $h(0+)=h(0)=0$. 
Moreover, this function is unique if we assume 
that the exploration of the tree induced by such a coding backtracks as less as
possible. We also prove that a measure-change on the tree corresponds to a 
re-parametrization of the coding
function. In addition, we describe 
several path-properties of the coding function in
terms of the metric properties of the real tree.

\vspace{4mm}

\noi {\it MSC 2000 subject classifications}:  Primary: 54F50, 60B99, 26A99;
secondary: 05C05, 05C12, 06A05, 26A46, 54B15, 54C30, 54E35, 54E45, 54E70,
54F05. 

\vspace{4mm}

\noi {\it Key words and phrases}: real tree, continuum random tree, height process, contour process,
  coding of a tree, linear order.

\end{abstract}

\section{Introduction}

    Real trees form a class of loop-free length spaces, which turns out to be the class of limiting objects 
of many combinatorial and discrete trees. More precisely,  
we say that a metric space $(T,d)$ is 
{\it a real tree} if it satisfies the following conditions:  
\begin{itemize} 
  \item For all $\sigma,\sigma' \in T$, there is an isometry $f_{\sigma,\sigma'}:[0,d(\sigma,\sigma')]\rightarrow T$ 
such that 
    $f_{\sigma,\sigma'}(0)=\sigma $ and $f_{\sigma,\sigma'}(d(\sigma,\sigma'))
=\sigma'$. We introduce the following notation 
$$\lgeo \sigma,\sigma'\rgeo
:=f_{\sigma,\sigma'}([0, d(\sigma,\sigma')])\; . $$

  \item If $q$ is a continuous injective map from $[0,1]$ into $T$, we have
    $$ q([0,1])= \lgeo q(0), q(1)\rgeo \; .$$
\end{itemize}

   Let us introduce some notation: we denote by 
$\rgeo \sigma,\sigma'\rgeo $, $\lgeo \sigma,\sigma' \lgeo$ and 
$\rgeo \sigma,\sigma'\lgeo$ the images of resp. 
$(0,d(\sigma,\sigma')] $, $[0,d(\sigma,\sigma')) $ and 
$(0,d(\sigma,\sigma')) $ by $f_{\sigma,\sigma'}$. 
For any $\sigma \in T$ we denote by ${\rm n}(\sigma ,T)$ 
the {\it degree} of $\sigma$, namely the 
(possibly infinite) 
number of connected components of $T \backslash \{ \sigma \}$. For convenience
of notation, we often denote ${\rm n}(\sigma ,T)$ simply by ${\rm n}(\sigma)$ 
when there is no risk of confusion. We denote by 
$${\rm Lf}(T)=\{ \sigma \in T\backslash \{ \rho \}: \; {\rm n}(\sigma ,T)=
1\} \quad  {\rm and } 
\quad {\rm Br}(T)=\{ \sigma\in T\backslash \{ \rho \}: \; 
{\rm n}(\sigma ,T)\geq 3\}$$
respectively the set of {\it leaves} of $T$ and 
the set of {\it branching points} of $T$. By convention, the
root $\rho$ is neither a leaf nor a branching point. We also denote by 
${\rm Sk} (T)$ the {\it internal skeleton} of $T$: ${\rm Sk} (T)= T\backslash
{\rm Lf}(T)$. We can easily prove that for any 
sequence $\sigma_n$, $n\geq 1$, that is dense in $T$, we have 
\begin{equation}
\label{intskel}
{\rm Sk} (T)=\bigcup_{n\geq 1} \lgeo \rho, \sigma_n \lgeo \; . 
\end{equation}
\noindent
Since $T$ is compact, we easily show that 
${\rm Br}(T)$ is at most countable (see Lemma 3.1 in \cite{DuWi1}).
We shall also need to introduce {\it the length measure} of a real tree $(T,
d)$ denoted by $\ell_T$. The length measure is defined on the trace on ${\rm Sk}(T)$
of the Borel sigma-field of $T$ and it is characterized by 
$$ \ell_T (\rgeo \sigma , \sigma' \lgeo )= d(\sigma, \sigma') \; . $$
The length measure can also be seen as the one-dimensional Hausdorff measure
on $T$.

  Real trees have a characterization called the {\it four points condition} 
that asserts that if $(X,d)$ is complete path-connected metric space then it is a real tree iff 
\begin{equation}
\label{fourpoint}
d(\sigma_1 , \sigma_2) +d(\sigma_3 , \sigma_4) \leq 
(d(\sigma_1 , \sigma_3) +d(\sigma_2 , \sigma_4))\vee (d(\sigma_1 , \sigma_4) +d(\sigma_2 , \sigma_3)), 
\end{equation}
for all $\sigma_1 , \sigma_2, \sigma_3 , \sigma_4 \in T$. The four points
condition has been first investigated independently by 
K. A. Zareckii \cite{Zar}, J.M.S. Sim{\~o}es Pereira \cite{SimPer} and
P. Buneman \cite{Bun}. See A. Dress, V. Moulton and W. Terhalle 
\cite{Dress84,DMT96,DT96} and also \cite{Chis, MayOv} for general
results concerning real trees. We also refer to the works of D. Aldous \cite{Al1, Al2} and 
of J-F. Le Gall \cite{LG2} for a first study of the Brownian 
Continuum Random Tree (CRT for short) that is 
coded by the normalized Brownian excursion. We refer to 
S. N. Evans \cite{Ev00} for the first
explicit use of real tree to construct random trees (see also 
\cite{EvPitWin,EvWin}); 
the reader interested in applications to phylogenetic models 
may consult the books of J. Felsenstein \cite{Fels} and of C. Semple and
M. Steel \cite{SemSte}; in a different direction, for applications to 
Super-Brownian motion, see \cite{LG93,LG1}; see also \cite{Krebs, Croy}
for a study of the Brownian motion on the Brownian CRT. 
We refer to the work of J-F. Le Gall and Y. Le Jan \cite{LGLJ1} 
for the definition of L\'evy trees that are random
real trees generalizing Aldous's Brownian CRT; see also \cite{LGLJ2,DuLG} for
application to general superprocesses and \cite{Du2,DuLG2,DuLG3,DuWi1,Weil} for
fractal and probabilistic properties of L\'evy trees. We refer to the  
work of D. Aldous, J. Pitman and G. Miermont \cite{AlMiPi1,AlMiPi2,AlMiPi3} 
for a detailed account on 
inhomogeneous continuum random trees (inhomogeneous continuum random trees
generalize the CRT and they are the possible scaling limits of 
interesting discrete combinatorial trees in
connection with random mappings). See the papers of B. Hass and G. Miermont \cite{Mier03, Mier05, HaaMi} for 
fragmentation processes linked with real trees. 
Let us also mention   
that in \cite{LyoHam} T. Lyons and B. Hambly  use real trees and tree-like 
paths for rough path 
integration theory.

\vspace{4mm}

   It has been shown by S.N. Evans, J. Pitman and A. Winter in \cite{EvPitWin} that the set of isometry 
classes of compact real trees endowed with the 
Gromov-Hausdorff distance is a complete separable metric space. However
there seems to be no natural way to choose a 
representent in a given isometry 
class. This 
contrasts with the discrete case. Indeed, if we consider a finite ordered 
rooted tree that is a finite planar graph without cycle with a 
distinguished vertex, then it is possible to label its vertices with 
words written with positive integers (see \cite{Ne}). More precisely, set 
$\U = \{ \varnothing \} \cup \bigcup_{n\geq 1} (\N^* )^n $, where 
$\N^*$ is the set 
$ \{ 1, 2, \ldots \} $ of positive integers and where $\varnothing$ 
stands for the empty word; an ordered rooted tree can be viewed as a 
subset $\tau$ of $\U $ satisfying the following conditions: 
\begin{description}
\item{(i)} $\varnothing \in \tau$ and $\varnothing $ is called the 
{\it root} of $\tau$.

\item{(ii)} If $v=(v_1, \ldots , v_n)\in \tau$ then, $(v_1, \ldots , v_{k})
\in \tau $ for any $1\leq k \leq n$.  

\item{(iii)} For every $v=(v_1, \ldots , v_n)\in \tau $, 
there exists $k_v (\tau) \geq 0 $ such that 
$(v_1, \ldots , v_n , j)\in \tau$ for every $1\leq j\leq k_v (\tau)$. 
\end{description}
If $v=(v_1, \ldots , v_n)\in \tau$, then $|v|=n$ is its {\it height} 
in $\tau$, that is its distance from the root (so we set $|\varnothing|=0$). 
Observe that $\U$ is linearly (or totally) ordered by the {\it  lexicographical order} 
denoted by $\leq $. If $\tau$ is finite, then we can list its vertices in an 
increasing sequence with respect to the lexicographical order, namely 
$\varnothing =v(0) < v(1) < \ldots < v(\# \tau -1)$. We define the 
{\it height process} of $\tau$ by 
$$ H_n(\tau)= |v(n)| \; ,\quad 0\leq n < \# \tau \; .$$
Clearly $H(\tau)$ characterizes the tree $\tau$ and in particular for any 
 $0\leq m \leq n < \# \tau$, the youngest common ancestor of $v(m)$ and $v(n)$
is situated at height $\min \{ H_k(\tau); m\leq k \leq n  \}$. 
Thus, the distance 
between  $v(m)$ and $v(n)$ is given in terms of $H(\tau)$ by 
$$ {\rm dist} (v(m) ,v(n))= H_m(\tau)+ H_n(\tau)-2\min_{m\leq k \leq n } 
H_k(\tau). $$

   One of the aim of this paper is to provide a similar coding for compact 
real trees and also an uniqueness result for such a representation. It turns 
out that the relevant class of coding functions for compact real trees are 
the left-continuous with right-limit functions: 
such functions are called {\it caglad functions} in the standard 
probabilistic terminology (caglad standing for ``continu \`a gauche et avec 
limite \`a droite'' in french). We shall explain further why 
the set of caglad functions is the right class of coding 
functions to consider (see Comment \ref{nonconcon}).

   Let us be more specific: for any $M\geq 0$, 
let us denote by $\cH_M$ the set 
of non-negative caglad functions $h$ on $[0,M]$ such that 
$h(0)=h(0+)=0$ and $h(t)-h(t+)\geq 0$, $t\in [0, M)$. The set 
$\cH=\cup_{M\geq 0}\cH_M$ is called the {\it set of height functions} and 
if $h\in \cH_M$, $\zeta (h)=M$ is called the {\it lifetime} of $h$. 
Let $h\in \cH_M$. For every $s,t\in [0,M]$, we set
$$m_h(s,t)=\inf_{r\in[s\wedge t,s\vee t]}h(r)$$
and
$$d_h(s,t)=h(s)+h(t)-2m_h(s,t).$$
We can easily show that for any $s_1, s_2, s_3, s_4 \in 
[0,M]$ we get 
\begin{equation} 
\label{fourtimes}
d_h(s_1, s_2) +d_h(s_3, s_4) \leq 
(d_h(s_1, s_3) +d_h(s_2, s_4))\vee (d_h(s_3, s_2) +d_h(s_1, s_4)).
\end{equation}
In particular, it implies the triangle inequality
by taking $s_3=s_4$. We introduce the equivalence relation $\sim_h $ defined
by $s\sim_h t$ iff $d_h(s,t)=0$ 
(or equivalently iff $h(s)=h(t)=m_h(s,t)$). Let
$T_h$ be the quotient space
$$T_h=[0,M]/ \sim_h.$$
The function $d_h$ induces a distance on $T_h$ that is also denoted by $d_h$. 
Thus, $(T_h, d_h)$ is a metric 
space satisfying the four points conditions. 
Denote by
$p_h:[0,M]\longrightarrow
T_h$ the canonical projection. $p_h$ is continuous when $h$ is
continuous and then $(T_h, d_h)$ is compact and 
path-connected. Now observe that 
this construction can be done with any non-negative function on $[0,M]$ but the resulting 
metric space may not be path-connected (take for instance 
an increasing function with a unique positive jump). 
We shall prove in Lemma \ref{premiereetape} that if $h\in \cH_M$,  then $(T_h,
d_h)$ is a compact real tree. 

\vspace{3mm}

   Observe that the construction of $(T_h, d_h)$ provides interesting additional structures:
\begin{itemize} 
  
\item Firstly, the construction provides a special vertex 
$\rho_h= p_h (0)$ called the {\it root} of the tree. $T_h$ can be viewed as family tree and the root as the ancestor of the 
family; it induces a partial order $\preccurlyeq$ given by 
$$ \sigma\preccurlyeq\sigma'\qquad\hbox{iff}\qquad \sigma 
\in \lgeo \rho_h , \sigma '\rgeo .$$
This order is called the {\it genealogical order} associated with the 
rooted tree $(T_h, d_h , \rho_h)$.

\item Secondly, the construction provides 
a relation $\leq_h$ 
on the tree that is induced by the usual order on $[0, M]$. More precisely,  
$$ \sigma \leq_h \sigma' \qquad\hbox{iff}\qquad \inf p_h^{-1} ( \{ \sigma \}) 
\leq \inf p_h^{-1} (\{ \sigma ' \}) .$$
The relation $\leq_h$ is actually a linear (or total) order (antisymmetry is the only non
obvious point to prove: to that end use Lemma \ref{leftlimi}). 
This order is the analogue of the lexicographical order on discrete rooted 
ordered trees. 

\item  Thirdly, the construction provides a measure $\mu_h$ that is the 
measure on $T_h$ induced by the Lebesgue measure $\lambda$ on $[0, M]$. More precisely, for any 
Borel set $A$ in $T_h$: 
$$ \mu_h (A)= \lambda \left(  p_h^{-1} (A)\right) . $$

\end{itemize}

We call a rooted, 
linearly ordered and measured compact real tree a {\it structured} compact real 
tree. 
In this paper we investigate the problem to know which of the structured 
compact real trees can be obtained by such a 
construction. More precisely, we say that two structured compact real 
trees $(T, d, \rho, \leq , \mu )$ and
 $(T', d', \rho', \leq' , \mu' )$ are {\it equivalent} iff there exists 
an isometry $f$ from $(T, d)$ onto $(T', d')$ that preserves roots
(i.e. $f(\rho)=\rho'$), that preserves orders (i.e. 
$f(\sigma_1)\leq' f(\sigma_2)$ as soon as 
$\sigma_1 \leq \sigma_2$) and that preserves measures (i.e. $\mu'=\mu \circ
f^{-1}$). Let us introduce notation $\sigma \wedge \sigma'$ for the branching point of
$\sigma$ and $\sigma'$ that is defined by 
$$ \lgeo \rho , \sigma \wedge \sigma'\rgeo =\lgeo \rho , \sigma\rgeo \cap
\lgeo \rho , \sigma' \rgeo \; . $$
Let assume 
that $(T, d, \rho, \leq , \mu )$ is equivalent to a structured tree obtained 
by the coding via a height function $h$; of course the order and the measure
have to satisfy some compatibility conditions. More precisely, 
we claim that necessarily, $\leq $ and $\mu $ have to 
satisfy the following conditions: 
\begin{itemize} 
  
\item{({\bf Or1})} $\;$   For any 
$\sigma_1 , \sigma_2$ in $T$, if $\sigma_1 \in 
\lgeo \rho , \sigma_2 \rgeo$, then $\sigma_1 \leq \sigma_2$. 

\item{({\bf Or2})} $\;$ If $\sigma_1 \leq \sigma_2 \leq \sigma_3$, then 
$\gamma \in \lgeo \rho , \sigma_2 \rgeo$ where $\gamma $ stands for the branching 
point of $\sigma_1$ on the subtree spanned by $\rho$,  $\sigma_2 $ and 
$\sigma_3$, namely:
$$ \lgeo \rho , \gamma \rgeo= \lgeo \rho , \sigma_1 \rgeo \cap \left( 
\lgeo \rho , \sigma_2 \rgeo \cup \lgeo \rho , \sigma_3 \rgeo \right) \; . $$
(see Figure \ref{obstructionorder} and Remark \ref{metritri}).

\item{({\bf Mes})} $\;$   For any distinct 
$\sigma_1$ and $\sigma_2$ in $T$ such that $\sigma_1 < \sigma_2$, we have 
$$ \mu (\{ \sigma \in T\; : \; \sigma_1 < \sigma <\sigma_2 \} ) >0 . $$

\end{itemize} 
This claim shall be proved in Lemma \ref{necessary}. A linear order satisfying
(Or1) and (Or2) is said to be {\it compatible} (with the metric and the choice
of a root) and a measure satisfying
(Mes) is also said to be {\it compatible} (with the metric, the root and 
$\leq$).

\begin{remark}
\label{metritri}
If $\leq $ satisfies (Or2), then $\sigma_1 \leq \sigma_2 \leq \sigma_3$ implies 
$$ d(\sigma_1 , \sigma_3)\geq d(\sigma_2 \wedge \sigma_3 ,\sigma_3 ) \; $$
and $\gamma=\sigma_1 \wedge \sigma_2$. \cq 
\begin{figure}[ht]
\psfrag{ro}{$\rho$}
\psfrag{gaga}{$\gamma$}
\psfrag{s1}{$\sigma_1$}
\psfrag{ss1}{$\sigma_1 '$}
\psfrag{s2}{$\sigma_2$}
\psfrag{s3}{$\sigma_3$}
\psfrag{s23}{$\sigma_2\wedge \sigma_3$}
\centerline{\epsfbox{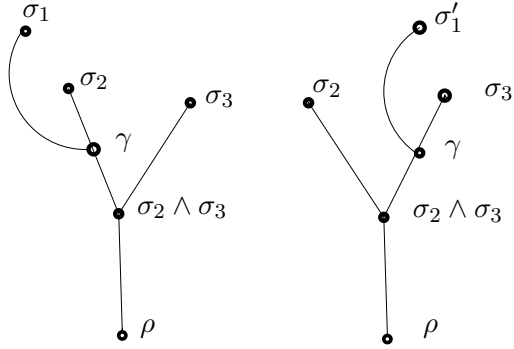}}
\caption{{\small In the first case we 
can find an order $\leq $ that satisfies (Or2) and 
$\sigma_1<\sigma_2<\sigma_3$, while in the second case, it is not 
possible to find an order $\leq $ that satisfies (Or2) and 
$\sigma_1' <\sigma_2<\sigma_3$.}}
\label{obstructionorder}
\end{figure}

\end{remark}
\vspace{4mm}

\begin{remark}
\label{evident}
Observe that if $\mu$ satisfies (Mes), then its 
topological support ${\rm supp} (\mu)$ is $T$. But the converse is not true:
see Figure \ref{poiludegarni}. \cq
\begin{figure}[ht]
\psfrag{roro}{$\rho$}
\psfrag{sigg}{$\sigma$}
\psfrag{ss}{$\sigma '$}
\psfrag{sss}{$\sigma \wedge \sigma'$}
\centerline{\epsfbox{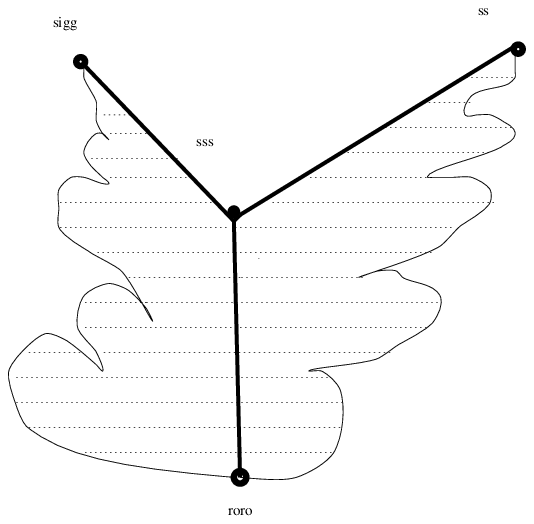}}
\caption{{\small 
Assume that $ \sigma ,\sigma'  $ are two distinct leaves of the compact rooted real tree $(T, d, \rho)$ and
assume that a bush is grafted on each vertex of a dense subset of the subtree 
$T'=\lgeo \rho  , \sigma \rgeo \cup \lgeo \rho  , \sigma' \rgeo $. Thus,  
$T\backslash T'$ (that is the dashed part of the tree) is dense in $T$. Then, 
there exists a finite Borel measure $\mu$ whose topological support is $T$ 
and such that $\mu(T')=0$. Now assume that the bushes grafted on 
$\rgeo \sigma \wedge \sigma' , \sigma \rgeo$ are all on the left and that 
the bushes grafted on $\lgeo \sigma \wedge \sigma' , \sigma' \rgeo$ are all on
the right so that $\rgeo \sigma \wedge \sigma' , \sigma' \lgeo=\{\xi \in T :
\sigma <\xi <\sigma' \}$. Now, observe that $\mu$ does not satisfy (Mes). }}
\label{poiludegarni}
\end{figure}

\end{remark}

\vspace{4mm}

\begin{remark}
\label{remorder}
We shall describe in Proposition
\ref{converseorder} all the compatible orders that can be defined on a given 
rooted compact real tree. We also explain in Section \ref{shuffle} that there is a natural way to 
pick uniformly at random these compatible orders. More precisely, on any fixed rooted compact real tree 
$(T,d, \rho)$ we shall construct a random compatible order 
denoted by $\leq_{{\rm Sh}}$ such that for any 
$\sigma_1 , \ldots , \sigma_n$ distinct elements of $T$, the random 
ordering induced on $\{ \sigma_1 , \ldots , \sigma_n \}$ by $\leq_{{\rm Sh}}$ is
uniformly distributed among all the distinct orderings of $\{ \sigma_1 ,
\ldots , \sigma_n \}$ induced by linear orders satisfying 
(Or1) and (Or2) (see Proposition \ref{unifshuffle} for details). This random 
order $\leq_{{\rm Sh}}$ is called the {\it uniform random shuffling of} $(T,d, \rho)$. We shall also prove in 
Proposition \ref{Shufflemes} that a.s. any finite 
Borel measure $\mu$ on $T$ whose topological support is $T$ satisfies 
(Mes) with respect to $(T, d, \rho ,  \leq_{{\rm Sh}})$.  \cq
\end{remark}

\begin{remark}
\label{nonuni}
Observe that different functions in $\cH$ may correspond to the 
same structured tree: consider for instance $h_1 \in \cH_1$  that is the
non-decreasing, continuous and piecewise linear height function such that $h_1(1/3)=h_1(2/3)=1/2$ and 
$h_1(1)=1$, and define $h_2 $ by 
$$ h_2(t)=\un_{[0, 1/2]} (t) \, h_1(2t) + \un_{(1/2 , 1]} (t) \, h_1 (2(1-t)) . $$
See Figure \ref{premfigu}. 
First note that $h_2\in \cH_1$ and then observe that $h_1$ and $h_2$ code 
the same structured compact real tree $(T,d, \rho,\leq , \mu)$ where 
$T$ is $[0,1]$, $d$ is the usual distance on $[0,1]$, $\rho$ is $0$, $\leq $
is the usual order on the line and $\mu =2/3 \lambda + 1/3 \delta_{1/2}$ 
(here $\lambda$ stands for the Lebesgue measure on $[0, 1]$ and 
$\delta_{1/2}$ is the Dirac mass at $1/2$).  \cq
\begin{figure}[ht]
\psfrag{un}{$1$}
\psfrag{tiers}{$\frac{1}{3}$}
\psfrag{2tiers}{$\frac{2}{3}$}
\psfrag{six}{$\frac{1}{6}$}
\psfrag{dem}{$\frac{1}{2}$}
\psfrag{cinq}{$\frac{5}{6}$}
\psfrag{zero}{$0$}
\psfrag{h1}{$h_1$}
\psfrag{h2}{$h_2$}
\centerline{\epsfbox{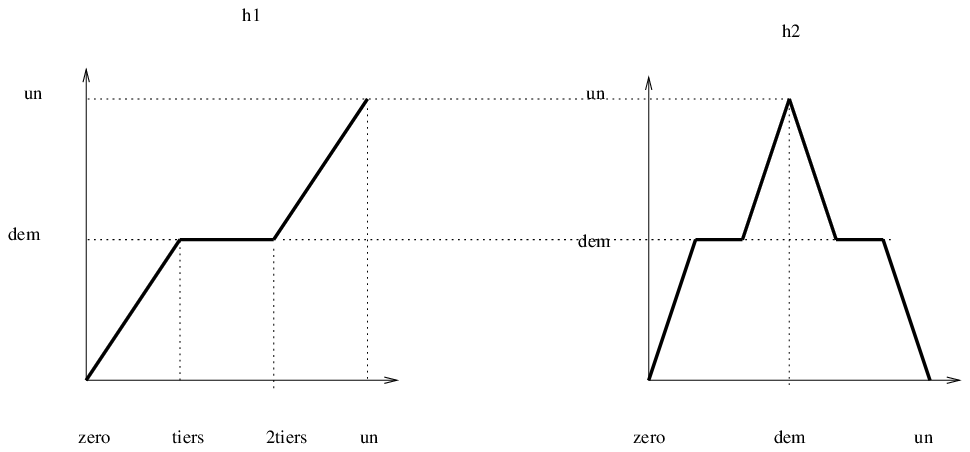}}
\caption{{\small $h_1$ and 
$h_2$ are two functions coding the same structured compact real tree. 
Observe that $F(h_1)=[0,1]$ and that $F(h_2)=[0,1/2]$. Thus, 
$h_1$ satisfies (Min) while $h_2$ does not.}}
\label{premfigu}
\end{figure}

\end{remark}

\vspace{4mm}

   As in the discrete case where the height process only visits the vertices 
once in the lexicographical order, we get uniqueness of the coding by 
requiring that the height fonction backtracks as less as possible. More 
precisely, let $h\in \cH$; for any $\sigma \in T_h$, set 
$$\ell (\sigma) =\inf p_h^{-1} (\{ \sigma  \} ) \quad {\rm and} 
\quad  r( \sigma)= \inf \{ t> \ell (\sigma ) \;  : \; p_h(t)\neq \sigma  \} . $$
We shall prove in Lemma \ref{leftlimi} that $p_h(\ell (\sigma ))=\sigma$ 
and the left-continuity of $h$ implies that $p_h(r(\sigma ))=\sigma$. So, we call the two sets 
$$ F(h)= \bigcup_{\sigma \in T_h} [\ell (\sigma ), r( \sigma) ] \quad {\rm and } 
\quad S(h) =[0, \zeta (h)] \backslash F(h)  $$
resp. the set of {\it times of first visit} and the set of {\it times of 
latter visit}. We shall prove in Lemma \ref{lattermesu} 
that $F(h)$ and $S(h)$ are Borel sets of the real line. 
We introduce the following property: 
\begin{itemize} 
  \item{({\bf Min})} $\;$   The height function $h\in \cH$ is said 
{\it minimal} iff $\lambda (S(h))=0$. 

\end{itemize}

\noindent
One of the two main results of the paper is the following.
\begin{theorem}
\label{main}
Let $(T,d,\rho, \leq, \mu)$ be a structured compact real tree such that $\leq
$ satisfies (Or1) and (Or2) and such that $\mu$ satisfies (Mes). There exists 
a unique $h\in \cH$ satisfying (Min) such that $(T,d,\rho, \leq, \mu)$ and 
$(T_h,d_h,\rho_h, \leq_h, \mu_h)$ are equivalent. 
\end{theorem}

\noindent
This theorem is proved in Section \ref{construc}. 

\begin{corollary}
\label{concontinu}
Let $(T,d)$ be a rooted compact real tree. There exists a continuous function 
$c \in \cH$ such that $(T_c,d_c)$ and $(T, d)$ are isometric. 
\end{corollary}

\noindent
{\bf Proof of Corollary \ref{concontinu}:} We first fix a root $\rho \in
T$. We can always find a probability measure $\mu$ whose topological support is $T$: consider for 
instance a sequence $\sigma_n $, $ n\geq 1$, that is dense in $(T,d)$ and 
define 
$$ \mu (d\sigma ) = \sum_{n\geq 1} 2^{-n} \delta_{\sigma_n} \, (d\sigma ) \; . $$
 As already mentioned in Remark \ref{remorder}, we can always find a linear order $\leq$
on $T$ such that the structured tree 
$(T,d,\rho, \leq, \mu)$ satisfies (Or1), (Or2) and (Mes). Denote by 
$h$ the coding height function associated with $(T,d,\rho, \leq, \mu)$ by Theorem \ref{main}. 
It implies in particular that 
$(T,d,\rho, \leq )$ and 
$(T_h,d_h,\rho_h, \leq_h )$ are equivalent. We conclude thanks to Lemma \ref{continuous} proved in Section \ref{propriri} 
that asserts that with any $h \in \cH$, we can always associate a (non-unique) continuous $c\in \cH$ 
such that $(T_h,d_h ,\rho_h , \leq_h)$ and $(T_c,d_c, \rho_c , \leq_c)$ are equivalent. \cqfd

\begin{remark}
\label{continu}
The result of the corollary has been proved independently by J-F. Le Gall 
\cite{LGman} by an approximation procedure. \cq 
\end{remark}

\begin{comment}
\label{nonconcon}
Let $(T,d,\rho, \leq, \mu)$ be a structured compact real tree such that $\leq
$ satisfies (Or1) and (Or2) and such that $\mu$ satisfies (Mes). Assume that
all the mass of $\mu$ is on the leaves of $T$, namely 
\begin{equation}
\label{feufeu}
\mu ({\rm Sk} (T)) =0 \; .
\end{equation}
Let $h\in \cH$ be such that $(T_{h},d_{h},\rho_{h}, \leq_{h}, \mu_{h})$ is
equivalent to $(T,d,\rho, \leq, \mu)$. Since 
$$p_h^{-1}({\rm Lf} (T)) \subset F(h) \; ,$$
$h$ automatically satisfies (Min). Consequently, there is a unique height function
that codes a structured compact real tree $(T,d,\rho, \leq,
\mu)$ that satisfies (Or1), (Or2), (Mes) and (\ref{feufeu}). (Let us mention that 
it is the case of the Continuum Random
Tree that is coded by the normalized Brownian excursion (see \cite{Al1, LG2})).

   We have shown in the proof of Corollary \ref{concontinu} that it is always
possible to find a continuous height function $c\in \cH$ that codes $(T,d,\rho ,\leq )$ 
where $\leq $ satisfies (Or1), (Or2). But if we fix $\mu$ on $T$ satisfying
(Mes), there may be no continuous function coding the structured tree 
$(T,d,\rho, \leq, \mu)$. For instance, assume that $\mu$ 
satisfies (\ref{feufeu}) and denote by $h$ the unique height function that 
codes $(T,d,\rho, \leq, \mu)$; choose $t\in (0, \zeta (h))$ such that 
$$ p_h (t) \in {\rm Lf} (T_h) \quad {\rm and} \quad h(t) >0 \; . $$
Define $h' \in \cH$ by 
$$ h'(s) =\un_{[0,t]} (s)\; h(s) + \frac{1}{2} \un_{(t, \zeta (h)]} (s) \; h(s) \; .$$
It is easy to check that $\mu_{h'} ({\rm Sk} (T_{h'}))=0$. Thus, $h'$ is the
unique height function that codes $(T_{h'},d_{h'},\rho_{h'}, \leq_{h'},
\mu_{h'})$ and obviously $h'$ is not continuous at $t$. \cq

\end{comment}

In view of the previous corollary and of Lemma \ref{continuous} in Section \ref{propriri}, 
let us note that there are many height functions coding the same compact real tree. 
More precisely, let $h'\in \cH$ and let  
$\varphi$ be an increasing continuous mapping from a finite interval $[0, M]$ onto $[0, \zeta (h')]$. 
Then, observe that the function defined by $h=h'\circ \varphi$ is in 
$ \cH_M$ and also observe that the two ordered rooted compact real trees 
$(T_h,d_h ,\rho_h , \leq_h)$ and $(T_{h'},d_{h'} ,\rho_{h'} , \leq_{h'})$ are equivalent. 
The time-change $\varphi$ only affects the measures $\mu_h$ and $\mu_{h'}$. Of course a measure-change does not always 
correspond to a re-parametrization of the coding functions: consider for
instance $h_1$ and $h_2$ as in Remark \ref{nonuni} and observe that $h_2$ cannot 
be obtained from $h_1$ by a time-change, while 
$h_1$ and $h_2$ code the same structured tree $(T, d, \rho, \leq ,
\mu)$. However, if we require that the 
height functions $h$ and $h'$ both satisfy (Min), then we can 
have a precise result explained in the following theorem that is proved in Section \ref{construc}. 

\begin{theorem}
\label{timechange}
Let $(T, d, \rho)$ be a rooted compact real tree and let $\leq $ be a linear order satisfying (Or1) and 
(Or2). Let $\mu$ and $\mu'$ be two finite Borel measures on $T$ that both satisfy (Mes). Denote by $h$ and $h'$ the height 
functions associated with resp. $(T, d, \rho , \leq, \mu)$ and $(T, d, \rho , \leq, \mu' )$ by 
Theorem \ref{main} ($h$ and $h'$ then satisfy (Min)). Then, there exists a non-decreasing and left-continuous mapping 
$\varphi: [0 , \mu (T)] \rightarrow [0 ,\infty ) $ such that 
$$ \varphi (0)=0 \; , \quad {\rm and } 
\quad h=h'\circ \varphi  \; .$$
Moreover, the following assertions are true:
\begin{itemize}
\item{(i)}   The time-change $\varphi$ is unique iff for any $ \sigma \in T$, 
$\mu (\{ \sigma \}) \mu' (\{ \sigma \}) =0$, namely iff $\mu$ and $\mu'$ do not
share any atom. 
 
\item{(ii)}  If $\mu' $ has no atom, then $\varphi$ is continuous. 

\item{(iii)} If $\mu$ has no atom, then $\varphi$ is increasing.

\item{(iv)} For any $t\in [0, \zeta (h))$ such that $\varphi (t) < \varphi (t+)$,  
$$ h'(u)=h(t) \; , \quad u \in [\varphi (t),  \varphi (t+)] \; .$$

\end{itemize}
\end{theorem}

\vspace{4mm}

\begin{remark}
\label{sophi}
Observe that if $\mu (\{ \rho \})=0$ and $\mu' (\{ \rho \}) >0$, 
then $\varphi (0+)>0$.  \cq 
\end{remark}

\begin{comment}
\label{comaa}
If $h, h'\in \cH$ are as in Theorem \ref{timechange}, then 
$h$ is continuous iff $h'$ is continuous. As already mentioned, 
Lemma \ref{continuous} asserts that any ordered compact tree $(T, d, \rho ,\leq )$ can be coded by a continuous height function 
$c \in \cH$. But the latter observation implies that in certain cases, it is impossible to find such a 
continuous function satisfying (Min). Consider for instance 
a ``Y-shaped'' tree rooted at the foot
of the ``Y'' (namely a tree with two leaves and a root distinct from the
branching point of the leaves) equipped with any of the two compatible
orders. This tree cannot be coded by a continuous function 
satisfying (Min) (see Figure \ref{Yarbre}). See Section \ref{propriri} for detailled
results about continuity properties of height functions.

   More generally, properties of 
$h$ as a path (i.e. properties that do not depend on any parametrization
of $h$) only concern $(T, d,\rho , \leq)$. In that vein, we prove in
Proposition \ref{vartot} that the total length $\ell_T(T)$ of a compact real tree $(T,d)$ is finite iff 
there exists a height function $h$ with bounded variation that 
codes $(T,d)$. \cq

\begin{figure}[ht]
\psfrag{rhho}{$\rho$}
\psfrag{sig}{$\sigma$}
\psfrag{sigprim}{$\sigma'$}
\psfrag{ssigpri}{$\sigma \wedge \sigma'$}
\psfrag{zer}{$0 $}
\psfrag{rsi}{$r(\sigma)$}
\psfrag{tmin}{$t_0 $}
\psfrag{lsig}{$\ell (\sigma')$}
\psfrag{h}{$h$}
\centerline{\epsfbox{yyarbre.eps}}
\caption{{\small Set 
$t_0=\inf \{ t \in (r(\sigma), \ell (\sigma' ) ) : \;  
h(t)= d_h (\rho , \sigma \wedge \sigma' ) \}$. Observe that if $h$ is 
continuous, then $t_0 > r(\sigma)$ but note that necessarily,  
$(r(\sigma),t_0)\subset S(h)$. Thus, $h$ does not satisfies (Min).}}

\label{Yarbre}

\end{figure}

\end{comment}

\vspace{4mm}

\begin{comment}
\label{combb}
Observe that if $h, h' \in \cH$ are as in Theorem \ref{timechange} and
if $\mu$ and $\mu'$ have no atom, then there exists a unique
increasing and continuous time-change $\varphi$ mapping $[0 , \zeta (h)]$ onto $ [0 , \zeta (h')]$ and 
such that $h=h'\circ \varphi $. \cq
\end{comment}

\vspace{4mm}

\begin{comment}
\label{comcc}
In Aldous's terminology, a {\it continuum tree} is a rooted compact real 
tree $(T, d, \rho)$ equipped with a probability measure $\mu$ that satisfies the following conditions: 
\begin{itemize} 
\item{{\bf (CT1)}} $\; $ The topological support of $\mu$ is $T$. 

\item{{\bf (CT2)}} $\; $ The measure $\mu$ is non-atomic.

\item{{\bf (CT3)}} $\; $ The measure $\mu$ is supported by the set of leaves of $T$, namely
  $\mu ({\rm Sk} (T))=0 $. 

\end{itemize}
First, let us mention that the definition of a continuum tree given here 
is slightly different from the definition 
given by Conditions (a), (b), (c) p. 253 and Assumption 
(17) p. 265 in \cite{Al2}: however the difference is inessential for our
purpose. Next, observe that
(CT1), (CT2), (CT3) imply a certain topological constraint on $(T,d,
\rho)$, namely that the leaves are dense in $T$ :
\begin{equation}
\label{denseleaves} 
\overline{{\rm Lf} (T)} =T \; .
\end{equation}
Conversely, if $T$ satisfies (\ref{denseleaves}), then we can find a
probability measure $\mu$ such that $(T, d, \rho , \mu)$ is a continuum tree
(See Proposition \ref{existconti}).

  Let us put the uniform random shuffling
$\leq_{{\rm Sh}}$ on $(T,d,\rho)$ and denote by $h_{{\rm Sh}}$ the height function
associated with $(T,d,\rho,\leq_{{\rm Sh}} , \mu)$ by Theorem \ref{main}. As a
direct consequence of Proposition \ref{shuffleconti} and Remark \ref{lamesu} 
in Section \ref{propriri}, a.s. the height function $h_{{\rm Sh}}$ 
is continuous and a.s. it is the
unique element of $\cH$ that codes $(T,d,\rho,\leq_{{\rm Sh}} , \mu)$. Moreover,
by Theorem \ref{timechange} and Comment \ref{combb}, a.s. for any $\mu'$ on $T$ satisfying 
(CT1), (CT2) and (CT3), the height function $h_{{\rm Sh}}'$ coding 
$(T,d,\rho,\leq_{{\rm Sh}} , \mu')$ is obtained from $h_{{\rm Sh}}$ by a 
unique increasing
continuous time-change $\varphi$. \cq 
\end{comment}

\vspace{4mm}

The paper is organized as follows: in Section \ref{constr}, we describe all
the possible linear orders satisfying (Or1) and (Or2); in Section
\ref{shuffle}, we define the
uniform random shuffling of a compact rooted real tree; in Section \ref{topol}
we prove topological properties of compatible linear orders. 
Section \ref{construc} is devoted to the proofs of Theorem
\ref{main} and Theorem \ref{timechange}. In the last section, we discuss special properties of
height functions and we
make the connection with an earlier result of Aldous (namely Theorem 15 in \cite{Al2}) that 
provides a randomized construction of height functions in the special case of
continuum trees. We
conclude Section \ref{propriri} with a probabilistic example illustrating the
effect of order-change on height functions.   

\vspace{4mm}


\section{Compatible linear orders.}
\label{complinord}

\subsection{Construction.}
\label{constr}

Let $(T, d , \rho )$ be a rooted compact real tree. In this section we 
construct all compatible linear orders on $T$. But first, let us 
prove the following proposition.
\begin{lemma}
\label{premiereetape}
For any $h\in \cH $, $(T_h,d_h)$ is a compact real tree.
\end{lemma}
{\bf Proof:}  Let us prove that if $h$ is in $\cH_M$, then $(T_h,d_h)$ is 
path-connected. Recall that $p_h$ stands for the canonical projection from 
$[0,M]$ onto $T_h$; note that $p_h$ is not necessary 
continuous; set $\rho_h= p_h(0)$; let $t_0\in [0, M]$ and set 
$\sigma =p_h(t_0)$; for any $s\in [0, h(t_0)]$ define $\sigma (s)= 
p_h(i(s))$   
where $i(s)$ is given by 
$$ i(s)=\sup \{ t \in [0, t_0] \; : \; h(t) \leq s \} \; .$$
Since $h$ is caglad without positive jump, $h(i(s))=h(i(s)+)=s$. Moreover, 
it is also easy to check that 
\begin{equation}
\label{easytocheck}
 \sigma (0)=\rho \; , \; \sigma ( h(t_0))=\sigma \quad {\rm and }
\quad d_h (\sigma (s),\sigma (s') )=|s-s'| \; , \; s,s'\in [0, h(t_0)] . 
\end{equation}
Thus,  $(T_h, d_h)$ is path-connected.

Let us prove now that $(T_h, d_h)$ is compact. 
Let $\sigma_n$, $n\geq 1$, be a $T_h$-valued sequence; let 
$t_n \in [0, \zeta (h)]$ , $n\geq 1$, be such that 
$p_h (t_n)= \sigma_n$; we can always find a monotone subsequence 
$t_{k_n}$, $n\geq 1$, that converges to $t\in [0, \zeta (h)]$ (say). If 
$t_{k_n}$, $n\geq 1$, is non-decreasing, then the left-continuity 
of $h$ implies that  
$d_h (t_{k_n} , t)$ goes to $0$ and consequently,  
$$ \lim_{n\rightarrow \infty} d_h (\sigma_{k_n} , p_h(t))=0 \, .$$ 
Assume that $t_{k_n}$, $n\geq 1$, is non-increasing. By definition $h(t+) \leq h(t)$ and set 
$$ t'=\sup \{ s \in [0, t] \; : \; h(s) \leq h(t+) \} \; .$$
Then observe that $d_h (t_{k_n} , t')$ goes to $0$, which implies that 
$$ \lim_{n\rightarrow \infty} d_h (\sigma_{k_n} , p_h(t'))=0 \, .$$ 
Thus, it proves that $(T_h, d_h)$ is compact. Use the four points condition to
complete the proof of the lemma. \cqfd

\vspace{4mm}

Next let us prove that for any $h\in \cH $, the relation 
$\leq_h$ is an order that satisfies (Or1), (Or2) and that 
the measure $\mu_h$ satisfy (Mes). Recall that 
$$\ell (\sigma) =\inf p_h^{-1} (\{ \sigma  \} ) \; .$$
We need the following lemma.
\begin{lemma}
\label{leftlimi} For any $\sigma \in T_h$, we have $ p_h (\ell (\sigma ) )=\sigma $. 
\end{lemma}
{\bf Proof:} If $h$ is continuous at $\ell (\sigma)$, then the result is
obvious. Assume that $h(\ell (\sigma)) >h(\ell (\sigma)+)$ and suppose that 
$ p_h (\ell (\sigma ) ) \neq \sigma $. Then, by definition of $\ell (\sigma)$,
there is a decreasing sequence
$t_n$, $n\geq 1$, converging to $\ell (\sigma )$ such that $ p_h (t_n ) =
\sigma $. Then, $ h(t_n)= d(\rho, \sigma)$ for any $n\geq 1$, and   
$$ \lim_{n\rightarrow \infty} h(t_n)= h(\ell (\sigma ) +) = d(\rho, \sigma) <h(\ell (\sigma ) ) \; .$$
Set $s_0=\sup \{ s \in [0, \ell (\sigma ) ) \; :\; h(s)\leq h(\ell (\sigma)+) \}$. Observe that 
$m_h (s_0 , t_n)=h(\ell (\sigma ) +)$, $n\geq 1$. It implies that $p_h
(s_0)= p_h (t_n ) =\sigma $. But clearly $s_0 < \ell (\sigma ) $ which contradicts the
definition of $\ell (\sigma ) $. Thus, $ p_h (\ell (\sigma ) ) =\sigma $,
which completes the proof of the lemma. \cqfd

\vspace{4mm}

This lemma implies that the relation $\leq_h$ is antisymmetric, which is the only
non-obvious point to justify in order to prove that $\leq_h$ is a linear
order. 

\begin{proposition}
\label{necessary}
Let $h\in \cH$. Then, $\leq_h$ satisfies (Or1) and (Or2) and $\mu_h$ satisfies 
(Mes). 
\end{proposition}
\noindent
{\bf Proof:} Let 
$\sigma_1 , \sigma_2 \in T_h$ be such that $\sigma_1 \in 
\lgeo \rho_h , \sigma_2 \lgeo $. Set 
$$t_0=\sup \{ 0 \leq s \leq \ell (\sigma_2) \; : \; h(s) \leq d_h( \rho_h , 
\sigma_1 )\} .$$
Since $h$ is caglad without negative jump, we get $h(t_0)=d_h (\rho_h,
\sigma_1)=h(t_0 +)$. Thus by the previous lemma 
$$ m_h(t_0, \ell (\sigma_2))= h(t_0)=d_h (\rho_h,
\sigma_1) < d_h (\rho_h,
\sigma_2)=h(\ell (\sigma_2 )) $$
and $p_h(t_0)=\sigma_1$. Then 
$\ell (\sigma_1)\leq t_0 <\ell (\sigma_2)$. Thus, $\sigma_1 \leq_h \sigma_2$, which proves 
that $\leq_h$ satisfies (Or1).

   Let us prove that $\leq_h$ satisfies (Or2): let $\sigma_1
\leq_h \sigma_2 \leq_h \sigma_3$, which implies that $\ell (\sigma_1)
\leq \ell (\sigma_2) \leq \ell (\sigma_3)$ by definition. By the previous
lemma,  $ p_h (\ell (\sigma_i ) )=\sigma_i $, $i\in \{ 1,2,3\}$. Consequently,  
$$ d_h( \rho_h , \sigma_1 \wedge \sigma_3) =
\min_{\ell (\sigma_1)\leq t \leq \ell (\sigma_3)} h(t) \leq 
\min_{\ell (\sigma_1)\leq t \leq \ell (\sigma_2)}h(t)=
d_h( \rho_h , \sigma_1 \wedge \sigma_2) \; ,  $$
which implies (Or2).

   Let us prove that $\mu_h$ satisfies (Mes): let $\sigma_1 <_h \sigma_2$. By
Lemma \ref{leftlimi}, we get $ m_h (\ell (\sigma_1), \ell (\sigma_2))= d_h (\rho_h , 
\sigma_1\wedge \sigma_2) $. The left-continuity of $h$ implies that there
exists $t_0 \in [ \ell (\sigma_1), \ell (\sigma_2))$ such that 
$$ \forall t\in [t_0 , \ell (\sigma_2)] \; , \quad h (t) > m_h (\ell (\sigma_1),
\ell (\sigma_2)) \; . $$
Choose $t\in (t_0,  \ell (\sigma_2)]$ and set $\sigma =p_h(t)$. Observe that
for any $s\sim_h t$ and any $s'\in [0, \ell (\sigma_1)]$
$$ m_h (s', \ell (\sigma_2)) \leq  m_h (\ell (\sigma_1), \ell (\sigma_2))= d_h (\rho_h , 
\sigma_1\wedge \sigma_2) <m_h (t, \ell (\sigma_2)) = m_h (s, \ell (\sigma_2)).$$
It implies that $\inf p_h^{-1} (\{ \sigma \}) \geq t_0 \geq \ell (\sigma_1)$ and
thus $\sigma_1 <_h \sigma$. Consequently, 
$$ (t_0 ,  \ell (\sigma_2)) \subset  p_h^{-1} (\{ \sigma \; : \; 
\sigma_1 <_h  \sigma <_h \sigma_2 \} ) , $$
which implies (Mes) since $t_0 <\ell (\sigma_2)$.  \cqfd

\vspace{4mm}

   Let $(T, d, \rho) $ be a rooted compact real tree. To avoid trivialities,
we assume that $T$ is not a point. We now construct a 
compatible linear order on $T$. To that end we need to introduce some 
notation: for any $\sigma \in T$, we denote by $\cC_{\sigma}$ the set of the connected 
components of $T\backslash \{ \sigma \}$ that do not contain the root $\rho$.
Observe that  $\cC_{\sigma}$ is empty iff $\sigma$ is a leaf. 
We also introduce the following set $\cI$ 
$$ \cI =\left\{ \emptyset ; \bN^* ; \{1, \ldots , n \} \; , \; n\in 
\bN^{*} \right\} . $$
We think of $\cI$ as a familly of indexing sets. 
More precisely, with any $\sigma \in T$ we associate $I_{\sigma} \in \cI$ and a bijection 
$C_{\sigma}$ from $I_{\sigma} $ onto $\cC_{\sigma}$ such that  
$$ \cC_{\sigma} = \{C_{\sigma} (k) \; , \; k\in I_{\sigma} \} .$$
Recall that ${\rm Br}(T)$ stands for the set of branching points of $T$ and
that $\rho \notin {\rm Br}(T)$, by convention. 
For any $\sigma \in {\rm Br}(T)\cup \{ \rho \}$ choose a linear order on $I_{\sigma}$ that is denoted 
by $\lhd_{\sigma}$. Set  
${\bf O}=\{\lhd_{\sigma} \; , \;  \sigma \in {\rm Br}(T)\cup \{ \rho \}\}$. We define a binary relation ${\bf \lhd_{o}}$ on $T$ 
in the following way: let $\sigma $ and 
$\sigma'$ be two distinct elements of $T$.
\begin{itemize} 
  
\item{({\bf Def1)}} $\;$  If $\sigma \wedge\sigma' \in \{ \sigma ,\sigma' \}$ 
then we set $\sigma {\bf \lhd_o} \sigma'$ if 
$\sigma \in \lgeo \rho , \sigma' \rgeo $ and we set $\sigma' {\bf \lhd_o} \sigma$ 
if $\sigma' \in \lgeo \rho , \sigma \rgeo $.

\item{({\bf Def2)}} $\;$  If $\sigma \wedge\sigma' \notin 
\{ \sigma ,\sigma' \}$, then $\sigma \wedge\sigma' \in {\rm Br}(T)\cup \{ \rho \}$ and 
there exist two distinct integers $k$ and $k'$ in 
$I_{\sigma \wedge\sigma'}$ such that $\sigma \in C_{\sigma \wedge\sigma'}
(k)$ and $\sigma' \in C_{\sigma \wedge\sigma'}
(k')$; then we set $\sigma {\bf \lhd_o} \sigma'$ if 
$ k\lhd_{\sigma \wedge\sigma'} k'$ and  $\sigma' {\bf \lhd_o} \sigma$ if 
$ k'\lhd_{\sigma \wedge\sigma'} k$. 

\end{itemize}
\noindent

\begin{proposition} 
\label{constrordre}
${\bf \lhd_o}$ is a linear order satisfying (Or1) and (Or2).
\end{proposition}
\noindent 
{\bf Proof:} By definition, for any $\sigma , \sigma'$ in $T$, we either 
have $\sigma {\bf \lhd_o} \sigma'$ or $\sigma' {\bf \lhd_o} \sigma$ so that 
the relation ${\bf \lhd_o}$ is linear. Let us prove that ${\bf \lhd_o}$ is antisymmetric. 
Suppose that 
$\sigma {\bf \lhd_o} \sigma'$ and $\sigma' {\bf \lhd_o} \sigma$; if 
$\sigma \wedge\sigma' \in \{ \sigma ,\sigma' \}$ then (Def1) easily implies 
that $\sigma =\sigma'$; suppose that $\sigma \wedge\sigma' \notin 
\{ \sigma ,\sigma' \}$; then with the notation of (Def2), we should have 
$ k\lhd_{\sigma \wedge\sigma'} k'$ and $ k'\lhd_{\sigma \wedge\sigma'} k$, 
which implies that $k=k'$; thus $\sigma $ and $\sigma'$ would be in the same 
connected component of $T\backslash \{ \sigma \wedge\sigma'\}$, 
which is impossible by definition of $ \sigma \wedge\sigma'$.

  Let us prove that ${\bf \lhd_o}$ is transitive. Let $\sigma_1 , 
\sigma_2, \sigma_3 \in T$ be such that $\sigma_1 {\bf \lhd_o}\sigma_2$ and 
 $\sigma_2 {\bf \lhd_o}\sigma_3$. To avoid trivialities, we assume that 
$\sigma_1 \wedge \sigma_3 \notin \{ \sigma_1 , \sigma_3 \}$ and that 
$\sigma_1$ and $\sigma_2$ are distinct. 
Let $\gamma $ be such that 
$$ \lgeo \rho ,\gamma \rgeo  = \lgeo \rho ,\sigma_1\rgeo \cap 
\left( \lgeo \rho ,\sigma_2\rgeo  \cup \lgeo \rho ,\sigma_3\rgeo \right) . $$  
First assume that $\sigma_1 \wedge \sigma_2 \in 
\lgeo \sigma_1 \wedge \sigma_3 ,\sigma_1\rgeo $. Then 
$\sigma_2 \wedge \sigma_3 =\sigma_1 \wedge \sigma_3$ and $\sigma_1$ and 
$\sigma_2$ are in the same connected component of 
$\cC_{\sigma_1 \wedge \sigma_3}$. There exist two distinct integers 
$k$ and $k'$ in 
$I_{\sigma_1 \wedge\sigma_3}$ such that 
$$\sigma_1 , \sigma_2 \in C_{\sigma_1 \wedge\sigma_3} (k)
\quad {\rm and } \quad \sigma_3 \in C_{\sigma_1 \wedge\sigma_3} (k') \; .$$
By definition, $\sigma_2 {\bf \lhd_o}\sigma_3$ implies 
$k \lhd_{\sigma_1 \wedge\sigma_3} k'$. Since 
$\sigma_2 \wedge \sigma_3 =\sigma_1 \wedge \sigma_3$, we also get 
$\sigma_1 {\bf \lhd_o}\sigma_3$. 
In addition, observe that $\gamma = \sigma_1 \wedge \sigma_2 $ so that 
(Or2) is verified.

   If we assume next that $\sigma_2 \wedge \sigma_3 \in 
\lgeo \sigma_1 \wedge \sigma_3 ,\sigma_3\rgeo $, then we can show by similar 
arguments that $\sigma_1 {\bf \lhd_o}\sigma_3$ and that 
$\gamma \in \lgeo \rho , \sigma_2 \rgeo $.

   It remains to consider the case $\sigma_2 \wedge \sigma_3 \in 
\lgeo \rho ,\sigma_1 \wedge \sigma_3 \lgeo$: if 
$\sigma_2 =\sigma_2 \wedge \sigma_3$, then $\sigma_2 \in \lgeo \rho , 
 \sigma_1 \rgeo$. By (Def1), it 
implies $\sigma_2 {\bf \lhd_o}\sigma_1$; we have shown that it implies 
$\sigma_2 =\sigma_1$ which contradicts the assumption that 
$\sigma_1$ and $\sigma_2$ are distinct. 
Thus, $\sigma_2 \neq \sigma_2 \wedge \sigma_3$. Consequently, 
$\sigma_1$ and 
$\sigma_3$ are in the same connected component of 
$\cC_{\sigma_2 \wedge \sigma_3}$; so there exist two distinct integers 
$k$ and $k'$ in 
$I_{\sigma_2 \wedge\sigma_3}$ such that 
$$\sigma_1 , \sigma_3 \in C_{\sigma_2 \wedge\sigma_3} (k)
\quad {\rm and } \quad \sigma_2 \in C_{\sigma_2 \wedge\sigma_3}(k') \; .$$
But $\sigma_2 {\bf \lhd_o}\sigma_3$ implies 
$k \lhd_{\sigma_2 \wedge\sigma_3} k'$ and 
$\sigma_1 {\bf \lhd_o}\sigma_2$ implies 
$k' \lhd_{\sigma_2 \wedge\sigma_3} k$, which rises a contradiction. 
Thus, we cannot have  $\sigma_2 \wedge \sigma_3 \in 
\lgeo \rho ,\sigma_1 \wedge \sigma_3 \lgeo$.

   We have proved that ${\bf \lhd_o}$ is a linear order satisfying (Or2). 
Observe now that (Or1) is a direct consequence of (Def1), which completes 
the proof of the proposition. \cqfd

\vspace{4mm}

  Consider now a compatible $\leq $ linear order on $(T, d, \rho)$. Let 
$\sigma_0 \in T$ be such that $T\backslash \{ \sigma_0\} $ has at least 
two connected components $C$ and $C'$ that do not contain the root $\rho$. 
\begin{lemma}
\label{compoorder} Either $\sigma \leq \sigma'$ 
for any $\sigma\in C$ and any $\sigma'\in C'$, which is denoted by 
$C\leq C'$; either  $\sigma' \leq \sigma$ 
for any $\sigma\in C$ and any $\sigma'\in C'$, which is denoted by 
$C'\leq C$.
\end{lemma}
\noindent
{\bf Proof:} Suppose that we can find $\sigma_1 , \sigma_3 \in C$ and 
 $\sigma_2 \in C'$ such that $\sigma_1< \sigma_2 <\sigma_3$. Let $\gamma$ 
be as in (Or2). We have $\sigma_1 \wedge \sigma_3 \in C$; it implies 
$\gamma \in C$. Thus $\gamma \notin \lgeo \rho , \sigma_2 \rgeo $, which
 contradicts (Or2). Exchange the role of 
$C$ and $C'$ in the previous argument to complete the proof of the lemma. \cqfd

\vspace{4mm}

The lemma implies that for any $\sigma \in {\rm Br}(T)\cup \{ \rho \}$, we can find a linear order 
$\lhd_{\sigma}$ on $I_{\sigma}$ such that for any $k, l\in I_{\sigma}$
that satisfy $k\lhd_{\sigma} l$, we have $C_{\sigma} (k) \leq C_{\sigma} (l)$. 
Consequently, we have proved the following proposition. 
\begin{proposition}
\label{converseorder}
Any compatible linear order on $(T, d, \rho)$ is of the form 
${\bf \lhd_o}$, where 
$${\bf O}=\{ \lhd_{\sigma} \; , \; \sigma \in  {\rm Br}(T)\cup \{ \rho \}\} $$
stands for a certain choice of linear orders 
on the $I_{\sigma}$'s , $ \sigma \in {\rm Br}(T)\cup \{ \rho \}$. Moreover, this representation of
$\leq$ is unique.
\end{proposition}

As a consequence of this proposition, two compatible linear orders are obtained one from another by
re-ordering each set of indices $I_{\sigma}$ ,  $ \sigma \in {\rm Br}(T)\cup \{ \rho \}$. 

\subsection{Uniform random shuffling of real trees.}
\label{shuffle}
In this section we explain how to ``pick'' a compatible linear order 
uniformly at random. Let $(\Omega , \cF, \bP)$ be a probability space 
on which all the random variables that we consider are defined. Let $S$ be a set. 
We formally define a random order $\leq $ on $S$
as a random mapping 
$F: \Omega \times (S \times S) \rightarrow \{0, 1\}$ such that 
$F(\omega ; \sigma, \sigma')= \un_{\{ \sigma \leq  \sigma' \}} (\omega )$ and such that 
$\omega\rightarrow F(\omega; \cdot , \cdot)$ is $(\cF,\cG)$-measurable, where
we have set $$ \cG =\cE^{\otimes (S\times S)} \quad {\rm with} \quad \cE=\{ \emptyset,
\{0\}, \{1\},  \{0,1\} \} \; . $$

\begin{example}
\label{linnoo}
Consider for instance $S=\bN$ and denote by ${\cal O}$ the random uniform order 
defined by the following property: for any $k_1, \ldots , k_n$ distinct elements of $\bN$, the random ordering on 
$\{ k_1, \ldots , k_n\}$ induced by ${\cal O}$ is uniformly distributed among 
the $n!$ possible ones. ${\cal O}$ is unique in distribution and it can be 
constructed as follows: let $U_n $, $n\in \bN$ be i.i.d. random variables that
are uniformly distributed on $[0, 1]$; we set $n\, {\cal O}\, m$ iff $U_n \leq U_m$, 
$m, n \in \bN$, where $\leq $ stands here for the usual order on $[0,1]$ 
(see Lemma 10 in \cite{Al2}).  \cq 
\end{example}

  We shall give a similar construction for the uniform random compatible 
order on a fixed rooted compact real tree $(T, d, \rho)$ called the 
{\it shuffling of} $T$: Recall that ${\rm Br}(T)$ is at most countable; 
let $\{ U_{\sigma , k} \, ; \, \sigma\in {\rm Br}(T)\cup \{ \rho \}\, , \, k\in I_{\sigma }\}$ be 
a (countable) family 
of i.i.d. random variables that are uniformly distributed on $[0, 1]$. Define 
${\bf O}=\{ \lhd_{\sigma} \, , \, \sigma\in {\rm Br}(T)\cup \{ \rho \} \}$ by 
$$ k \lhd_{\sigma}l \quad {\rm iff} \quad  U_{\sigma , k} 
\leq  U_{\sigma , l} \; . $$
We define the random uniform shuffling of $T$ by $\leq_{{\rm Sh}} \, = \, 
\lhd_{{\bf O}}$. 
\begin{proposition}  
\label{unifshuffle}For any $\sigma_1 , \ldots , \sigma_n$  
distinct elements of $T$, the random 
ordering of the set $\{ \sigma_1 , \ldots , \sigma_n \}$ 
induced by $\leq_{{\rm Sh}}$ is uniformly distributed among 
the orderings of this set induced by linear orders satisfying (Or1) and (Or2).
\end{proposition}
\noindent
{\bf Proof:} This is a consequence of the construction of 
$\lhd_{{\bf O}}$ given in the previous section and of the result asserted in
Example \ref{linnoo}. The 
details are left to 
the reader. \cqfd

\vspace{4mm}

Observe that this proposition implies that $\leq_{{\rm Sh}}$ is unique in 
distribution. We now prove the following proposition:

\begin{proposition}
\label{Shufflemes}
Almost surely, any finite Borel 
measure $\mu $ whose topological support is $T$
satisfies (Mes) with respect to $\leq_{{\rm Sh}}$.
\end{proposition}
\noindent
{\bf Proof:} Let $\mu $ be a finite Borel measure whose topological support is $T$
and  let $\sigma_1$ and $\sigma_2$ be such that $\sigma_2 \notin 
\lgeo  \rho , \sigma_1 \rgeo $. Thus by definition of $\leq_{{\rm Sh}}$, 
$\bP(\sigma_1 <_{{\rm Sh}} \sigma_2  ) \geq 1/2 $. Let us fix 
$\omega_0 \in \{\omega \in \Omega \; : \;  \sigma_1 <_{{\rm Sh}} \sigma_2 \}$. There are two cases to consider: 
suppose first that 
$\mu (\rgeo \sigma_1 \wedge \sigma_2 , \sigma_2 \lgeo) >0$; then since 
we have fixed $\omega_0 \in \{ \omega \in \Omega \; : \;  \sigma_1 <_{{\rm Sh}} \sigma_2 \}$, we get 
$$\rgeo \sigma_1 \wedge \sigma_2 , \sigma_2 \lgeo \subset 
\{\sigma \in T \; : \;\sigma_1 <_{{\rm Sh}} \sigma <_{{\rm Sh}} \sigma_2  \}, $$  
which implies that 
$\mu \left( \{\sigma \in T \; : \;\sigma_1 \leq_{{\rm Sh}} \sigma 
\leq_{{\rm Sh}} \sigma_2  \}\right) $ is non-zero.

   Let us now suppose that 
$\mu (\rgeo \sigma_1 \wedge \sigma_2 , \sigma_2 \lgeo) =0$. Since the 
topological support of $\mu$ is $T$, there exists a sequence $s_n $ , 
$n\geq 1$, of branching points in 
$\rgeo \sigma_1 \wedge \sigma_2 , \sigma_2 \lgeo$ 
that is dense in this set; then for any $n\geq 1$, denote by $k_n \in I_{s_n}$ the 
index such that $\sigma_2 \in C_{s_n} (k_n)$. Fix $n\geq 1$, and take $\sigma  
\in C_{s_n} (k)$ with $ k \in I_{s_n} \backslash \{k_n\}$. Suppose that 
$\sigma \leq _{{\rm Sh}} \sigma_2$. Then, Lemma 
\ref{compoorder} implies that  $C_{s_n} (k) \leq _{{\rm Sh}} C_{s_n} (k_n)$ 
and since we have fixed 
$\omega_0 \in \{ \omega \in \Omega \; : \; \sigma_1 <_{{\rm Sh}} \sigma_2 \}$, (Or2) implies that $\sigma_1 
< _{{\rm Sh}} C_{s_n} (k)$. Consequently 
$$ \mu \left( \{\sigma \in T \; : \;\sigma_1 <_{{\rm Sh}} \sigma 
<_{{\rm Sh}} \sigma_2  \} \right) 
\geq \mu \left( C_{s_n} (k) \right) >0 \; . $$
Thus, if we fix 
$\omega_0 \in \{ \omega \in \Omega \; : \; \sigma_1 <_{{\rm Sh}} \sigma_2 \}$
and if $\mu \left( \{\sigma \in T \; : \;\sigma_1 <_{{\rm Sh}} \sigma 
<_{{\rm Sh}} \sigma_2  \} \right)=0$, then 
\begin{itemize}

\item{(a)} $\mu (\rgeo \sigma_1 \wedge \sigma_2 , \sigma_2 \lgeo) =0$.

\item{(b)} The denumerable set 
$$ \{ s_n \;, \;  n\geq 1 \}={\rm Br} (T)\cap \rgeo \sigma_1 \wedge \sigma_2 , \sigma_2
  \lgeo$$ 
is dense in  $\rgeo \sigma_1 \wedge \sigma_2 , \sigma_2 \lgeo $.

\item{(c)} For any $n\geq 1$, and any $k \in I_{s_n} \backslash \{k_n \}$, we 
have $ C_{s_n} (k_n)\leq _{{\rm Sh}} C_{s_n} (k) $. 
Thus, by definition of $\leq _{{\rm Sh}}$, it implies that  
$$ U_{s_n , k_n}(\omega_0 ) < U_{s_n , k} (\omega_0 )\; , \quad 
k \in I_{s_n} \backslash \{k_n \} \; , \; n\geq 1 \; .$$

\end{itemize}

\noindent Let $\sigma , \sigma' \in T$ such that 
$\sigma \in \lgeo \rho , \sigma' \lgeo$; for any $s\in {\rm Br} (T) 
\cap \rgeo \sigma, \sigma' \lgeo$, denote by $k(s) \in I_s$ the index such 
that $\sigma' \in C_s (k(s))$. Define
$$ A_{\sigma , \sigma'}= \bigcap_{s\in {\rm Br} (T) 
\cap \rgeo \sigma, \sigma' \lgeo} \{ \omega \in \Omega \; : \; 
 U_{s , k(s)}(\omega ) <  U_{s, k} (\omega ) \; , \; 
k \in I_{s} \backslash \{ k(s) \} \} , $$
with the convention $ A_{\sigma , \sigma'}= \emptyset $ if ${\rm Br} (T) 
\cap \rgeo \sigma, \sigma' \lgeo =\emptyset $. Clearly, 
if $ \# \;  {\rm Br} (T) \cap \rgeo \sigma, \sigma' \lgeo =\infty $, 
then $\bP(A_{\sigma , \sigma'})=0$. Thus, if we set 
$$ B=\bigcup \{ A_{\sigma , \sigma'} \; :  \quad \sigma , \sigma' 
\in  {\rm Br} (T)\cup \{ \rho \} \; : \; \sigma \in \lgeo \rho , \sigma' \lgeo 
\; {\rm and }\;  
\# \; {\rm Br} (T) \cap \rgeo \sigma, \sigma' \lgeo  =\infty \; \} \; , $$
then $\bP (B)=0$ since ${\rm Br} (T)$ is at most countable. Now observe that
(a), (b) and (c) imply that 
$$ \{ \omega \in \Omega \; : \; \sigma_1 <_{{\rm Sh}} \sigma_2 \} \cap 
\{ \omega \in \Omega \; : \; 
\mu ( \{\sigma \in T \; : \; \sigma_1 \leq_{{\rm Sh}} \sigma 
\leq_{{\rm Sh}} \sigma_2  \} ) =0  \}  \; \subset \; B , $$
which implies the proposition since $B$ does not depend on $\mu$, $\sigma_1$ 
or $\sigma_2$.  \cqfd

\subsection{Topological properties of compatible linear orders.}
\label{topol}
 
   In this section we prove properties of 
 compatible linear orders that shall be needed in the next section. 
Let  $\leq $ be a compatible linear order 
on the rooted compact tree $(T, d, \rho)$. By Proposition \ref{converseorder},
$\leq $ is of the form $\lhd_{{\bf O}}$, 
for a certain choice ${\bf O}=\{ \lhd_{\sigma} , \sigma
\in {\rm Br}(T)\cup \{ \rho \} \}$ of linear orders on the indexing sets $I_{\sigma}$ , $\sigma \in
{\rm Br}(T)\cup \{ \rho \}$.

    Let us first introduce some notation. Fix $\sigma \in T$. Recall
that $\cC_{\sigma}=\{ C_{\sigma} (k) , k\in I_{\sigma}\}$ stands for the set
of the connected components of $T\backslash \{ \sigma\}$ that do not contain
the root. Consider the connected components of 
$T\backslash \lgeo \rho ,\sigma \rgeo $ that are grafted on 
$\lgeo \rho ,\sigma \lgeo$: by Lemma \ref{compoorder}, either all the points
of a such component are smaller than $\sigma$, either all points of the connected component are
greater than $\sigma$. So we denote by 
$\cC_{\sigma}^{-}=\{ C^-_j , j\in J_{\sigma}^{-}\}$ the set of connected components grafted on 
$\lgeo \rho ,\sigma \lgeo$ that are smaller than $\sigma$ and we denote 
by $\cC_{\sigma}^{+}=\{ C^+_j , j\in J_{\sigma}^{+}\}$ the set of connected components grafted on 
$\lgeo \rho ,\sigma \lgeo$ that are greater than $\sigma$. For any $j\in J_{\sigma}^{\pm}$ we
denote by $\gamma_j^{\pm}$ the point of $\lgeo \rho ,\sigma \lgeo$ on which the
component $C_j^{\pm}$ is grafted. Observe that $\{\gamma_j^{\pm} \}\cup
C_j^{\pm}$ is the closure of $C_j^{\pm }$. Note that different components may be
grafted on the same point. We thus get 
$$ T\backslash \lgeo \rho ,\sigma \rgeo =\bigcup_{C\in \cC_{\sigma}^{-}\cup \cC_{\sigma}
\cup \cC_{\sigma}^{+}} C \; \; .$$

\begin{lemma}
\label{echelon}The following assertions are true.

\begin{itemize} 
\item For any $j_1 \in J_{\sigma}^{-}$, any $k\in I_{\sigma}$ and any $j_2 \in
  J_{\sigma}^{+}$, we get  
$$ C_{j_1}^- \leq C_{\sigma} (k) \leq C_{j_2}^+ \; . $$

\item For any $j_1,j_2  \in J_{\sigma}^{-}$ such that $d(\rho , \gamma_{j_1}^{-})
  <d(\rho , \gamma_{j_2}^{-})$, we get 
$$  C_{j_1}^- \leq  C_{j_2}^- \; .$$
 
\item For any $j_1,j_2  \in J_{\sigma}^{+}$ such that $d(\rho , \gamma_{j_1}^{+})
  > d(\rho , \gamma_{j_2}^{+})$, we get 
$$  C_{j_1}^+ \leq  C_{j_2}^+ \; .$$
\end{itemize}
\end{lemma}
\noindent
{\bf Proof:} This is a direct consequence of (Def1) and (Def2). \cqfd

\vspace{4mm}

   We consider the family of subsets of $T$ denoted by $ L_{\sigma}$,
$\sigma\in T$, and defined by 
$$ L_{\sigma}= \{\sigma' \in T\; : \;  \sigma' \leq  \sigma \} . $$
These subsets are called the {\it left sets of} $T$. We first prove the following
proposition. 
\begin{proposition}
\label{compa} 
For any $\sigma
\in T$, $ L_{\sigma}$ is a compact set.
\end{proposition}
\noindent
{\bf Proof:} Observe that 
$$ T\backslash L_{\sigma}=  \left( \bigcup_{C\in \cC_{\sigma}} C \right) 
\cup  \left(
    \bigcup_{ C^+ \in \cC^+_{\sigma}} C^+  \right) \; . $$
Thus $T\backslash L_{\sigma}$ is an open set, which proves the
proposition. \cqfd 

\vspace{4mm}

\begin{proposition}
\label{monot}
Every ($\leq \, $)-monotone sequence in $T$ is convergent. 
\end{proposition}
\noindent
{\bf Proof:} Let $\sigma_n \in T$ , $n\geq 1$, be a ($\leq \, $)-monotone
sequence. Suppose that it has two distinct limit points $\sigma$ and
$\sigma'$. Assume that $\sigma<\sigma'$ and choose 
$\sigma_0 \in \rgeo \sigma \wedge \sigma' , \sigma'\lgeo $. Denote by 
$C$ and $C'$ the two distinct connected components of 
$T\backslash \{ \sigma_0\}$ that contain respectively $\sigma$ and
$\sigma'$. Note that $C$ also contains the root $\rho$. Since the sequence is
monotone, we can find $n_1, n_2, n_3 \geq 1$, such that 
$$ \sigma_{n_1}<\sigma_{n_2}<\sigma_{n_3} \quad , \quad \sigma_{n_1},
\sigma_{n_3} \in C' \quad {\rm and} \quad 
\sigma_{n_2} \in C .$$
It first implies that $\sigma_{n_1}\wedge
\sigma_{n_3} \in C'$. Then observe that $\lgeo \rho , \sigma_{n_2} \rgeo \subset
C$, which implies that $\sigma_{n_2}\wedge \sigma_{n_3} \in C$ and which contradicts (Or2). Thus, the sequence has only one
limit point. \cqfd 

\vspace{4mm}

We first consider a non-decreasing sequence $\sigma_n \in T$ , $n\geq 1$, that
converges to $\sigma$. We distinguish three cases:
\begin{itemize} 
\item{{\bf Case (I)}}:  $\;$   $\sigma_n <\sigma$, for all $n\geq 1$. 

\item{{\bf Case (II)}}:  $\;$ $\sigma_n =\sigma$, for all sufficiently large $n\geq 1$.
 
\item{{\bf Case (III)}}:  $\;$ $\sigma <\sigma_n  $, for all sufficiently large $n\geq 1$.
\end{itemize}
We also set $D=\{ \sigma' \in T \; : \; \sigma_n <\sigma' \; , \; n\geq 1
\}$. We prove the following lemma. 
\begin{lemma}
\label{increa} 
The following assertions are true.
\begin{itemize}
\item In Case (I), we get 
$$ \bigcup_{n\geq 1} L_{\sigma_n} = L_{\sigma} \backslash \{ \sigma \} \quad
{\rm and } \quad \bigcap_{\sigma' \in D}  L_{\sigma'} =  L_{\sigma} .$$

\item In Case (II), we get 
$$  \bigcup_{n\geq 1} L_{\sigma_n} = L_{\sigma} =\bigcap_{\sigma' \in D}  L_{\sigma'} .$$

\item In Case (III), we get $\sigma \notin {\rm Lf}(T)$ and  
$$  \bigcup_{n\geq 1} L_{\sigma_n} =  L_{\sigma} \cup \left( \bigcup_{k\in K}
C_{\sigma}(k) \right)= \bigcap_{\sigma' \in D}  L_{\sigma'} ,$$
where we have set 
$$K= \{ k\in I_{\sigma} \; : \; \exists n\geq 1 \, , \,
\exists l\in I_{\sigma} \; : \; \sigma_n \in C_{\sigma}(l) \; {\rm and } \;
k\lhd_{\sigma} l  \} $$
if $\sigma \in {\rm Br}(T)\cup \{ \rho \}$ and $K=I_{\sigma }=\{ 1\}$ if ${\rm n} (\sigma)=2$.  
\end{itemize} 
\end{lemma} 
\noindent
{\bf Proof:} Let us first consider Case (I): suppose that there exists 
$\sigma'$ such that $\sigma_n <\sigma' <\sigma$ for all $n\geq 1$. Define 
$\gamma_n$ by 
$$ \lgeo \rho , \gamma_n \rgeo =\lgeo \rho , \sigma_n \rgeo \cap 
\left( \lgeo \rho , \sigma' \rgeo \cup  \lgeo \rho , \sigma\rgeo \right) .$$
Then we get $\gamma_n \in \lgeo \rho , \sigma' \rgeo $ by (Or2), which implies
that $d(\gamma_n , \sigma) \geq d(\sigma \wedge \sigma' , \sigma )$. Thus, for
any $n\geq 1$,  
$$ d(\sigma_n , \sigma)= d(\sigma_n , \gamma_n )+ d(\gamma_n , \sigma) \geq 
d(\sigma \wedge \sigma' , \sigma ) >0 , $$
which rises a contradiction. This prove the first point of the lemma. 

\vspace{3mm}

   The first equality in Case (II) is obvious. Suppose there exists $\sigma'' 
\in \bigcap_{\sigma' \in D}  L_{\sigma'} \backslash  L_{\sigma}$. Then,
$\sigma'' $ is a minimal element of $D$. But we can always 
find $\xi \in T$ such that $\sigma  < \xi  < \sigma''  $. It implies that $\xi$
is also a minimal element of $D$ distinct from $\sigma''$, which is 
absurd since $\leq $ is linear. It proves the second equality in Case (II).

\vspace{3mm}

Let us consider Case (III): We first suppose that there
are $n_0 \geq 1$, and $j_0 \in  J_{\sigma}^{+}$ such that $\sigma_{n_0} \in 
C_{j_0}^{+}$. Recall that $\gamma_{j_0}^+$ stands for the point of 
$\lgeo \rho , \sigma \lgeo$ on which the connected component $C_{j_0}^{+}$ is
grafted. Lemma \ref{echelon} implies that for any $n\geq n_0$, 
$\sigma_n$ is in a connected component of $\cC_{\sigma}^{+}$ that is grafted on a 
point of 
$\lgeo \rho , \gamma_{j_0}^+ \rgeo $. Thus, for all $ n \geq n_0 $
$$ d( \sigma_n , \sigma ) > d( \gamma_{j_0} ^+, \sigma) >0 \; ,$$
which rises a contradiction. It shows that the sequence $\sigma_n$, $n\geq 1$, has no
term in any of the connected
components of $\cC_{\sigma}^+$. It implies that $\sigma \notin {\rm Lf}(T)$ and that    
$$ \bigcup_{n\geq 1} L_{\sigma_n} \subset  \; L_{\sigma} \cup \left( \bigcup_{k\in
  I_{\sigma}} C_{\sigma}(k) \right) \; .$$
By similar arguments we also get 
$$ \bigcup_{n\geq 1} L_{\sigma_n} \subset   \; L_{\sigma}  \cup \left( \bigcup_{k\in K}
C_{\sigma}(k) \right) \; . $$
We now prove the reversed inclusion: suppose that there are $k_0 \in
I_{\sigma}$, $n_0 \geq 1$ and $\sigma'\in C_{\sigma}(k_0)$ such that
$$ \sigma_{n_0} \in C_{\sigma}(k_0) \quad {\rm and } \quad \sigma_n <\sigma'
\; , \quad n\geq 1 \; .$$
Then, $\sigma_{n} \in C_{\sigma}(k_0)$ for any $n\geq n_0$
and by definition of $K$ we get $K=\{ k\in I_{\sigma} \; : \; k\lhd_{\sigma}
k_0\}$. Since $\sigma_{n_0} \leq \sigma_{n} <\sigma' $, (Or2) implies that 
$\sigma_{n_0} \wedge \sigma' \in \lgeo \sigma , \sigma_{n} \wedge
\sigma'\rgeo $
and consequently 
$$ d(\sigma_{n}, \sigma ) \geq d(\sigma_{n} \wedge \sigma', \sigma  ) \geq d( \sigma_{n_0} \wedge \sigma' ,\sigma  ) >0 \; , $$
which rises a contradiction. So, it proves that if $\sigma_{n_0} \in
C_{\sigma}(k_0)$, then $ C_{\sigma}(k_0) \subset \bigcup_{n\geq 1}
L_{\sigma_n}$. This implies 
$$ L_{\sigma}  \cup \left( \bigcup_{k\in K}
C_{\sigma}(k) \right)
\subset \bigcup_{n\geq 1} L_{\sigma_n} \; , $$
and the first equality of Case (III) follows. Let us prove the second equality of Case (III): clearly, we have 
\begin{equation}
\label{triv}
\bigcup_{n\geq 1} L_{\sigma_n} \subset \bigcap_{\sigma' \in D}  L_{\sigma'} .
\end{equation}
Suppose now that there is a point $\sigma''$ 
in $\bigcap_{\sigma' \in D}  L_{\sigma'}$ such that $\sigma_n
<\sigma''$ for all $n\geq 1$. It implies that $\sigma'' \in D$. 
Since $\sigma''$ is in $\bigcap_{\sigma' \in D}  L_{\sigma'}$, we have 
$\sigma''\leq \sigma'$, $\sigma' \in D$. $\sigma''$ is the minimal point
of $D$ (there is at most one since the order is linear). Thus, 
for any $\xi \in T$, if $\xi <\sigma''$, then 
$\xi \notin D$ and, by definition, there exists
$n_0\geq 1$, such that $\xi \leq \sigma_{n_0}$. This implies 
\begin{equation}
\label{contrapo}
\bigcup_{n\geq 1} L_{\sigma_n} = L_{\sigma''} \backslash \{ \sigma'' \} \; .
\end{equation}
But the first equality of Case (III) implies
that $\bigcup_{n\geq 1} L_{\sigma_n} $ has to be a compact set, which
contradicts (\ref{contrapo}). So it proves that there is no point $\sigma''$ 
in $\bigcap_{\sigma' \in D}  L_{\sigma'}$ such that $\sigma_n
<\sigma''$ for all $n\geq 1$. This, combined with (\ref{triv}), implies 
$$ \bigcup_{n\geq 1} L_{\sigma_n} = \bigcap_{\sigma' \in D}  L_{\sigma'} \; , $$
which completes the proof of the Lemma. \cqfd

\vspace{4mm}

Consider now a non-increasing sequence $\sigma_n \in T$ , $n\geq 1$ that converges to $\sigma$. Set 
$$G=\{ \sigma' \in T \; : \;  \sigma' < \sigma_n \; , \; n\geq 1
\} \; . $$ 
We prove the following lemma. 
\begin{lemma}
\label{decrea}
We have 
$$ \bigcap_{n\geq 1} L_{\sigma_n}=  L_{\sigma} \cup \left( \bigcup_{k\in I_{\sigma}
  \backslash K}
C_{\sigma}(k) \right) \; , $$
where we have set 
$$K=\{ k\in I_{\sigma} \; : \; \exists n\geq 1 \, , \,
\exists l\in I_{\sigma} \; : \; \sigma_n \in C_{\sigma}(l) \; {\rm and } \;
l\lhd_{\sigma} k  \} $$ 
if $\sigma \in {\rm Br}(T)\cup \{ \rho \}$, 
$K = I_{\sigma} =\emptyset $ if $\sigma \in {\rm Lf}(T )$ and
$K=I_{\sigma}=\{ 1\}$ if ${\rm n}(\sigma)=2$ (observe that $I_{\sigma}
\backslash K$ may be empty). 
Moreover, if $G$ is non-empty and if $\sigma_n=\sigma$ for all sufficiently large $n$, then we get  
$$ \bigcup_{\sigma' \in G}  L_{\sigma'} =  L_{\sigma}  \backslash \{ \sigma \}
\subset  L_{\sigma} = \bigcap_{n\geq 1} L_{\sigma_n} . $$
Otherwise, we get 
$$  \bigcup_{\sigma' \in G}  L_{\sigma'} = \bigcap_{n\geq 1} L_{\sigma_n} . $$
\end{lemma}
\noindent
{\bf Proof:} The arguments are similar to those used to prove Lemma
\ref{increa}. The details are left to the reader. \cqfd 

\vspace{4mm}

We shall need the following lemma in Section \ref{construc}. 
\begin{lemma}
\label{pisyst}
The collection of sets $\{ \;   \emptyset  \, ; \,  T  \, ; \,   L_{\sigma} \; , \; \sigma \in T \; \}$ is a
$\pi$-system that generates the Borel sigma-field of $T$. 
\end{lemma}
\noindent
{\bf Proof:} Clearly, $\{ L_{\sigma} \; , \; \sigma \in T \}$ is closed under
finite intersection. This implies that 
$$\{ \;  \emptyset  \, ; \,  
T  \, ; \,   L_{\sigma} \; , \; \sigma \in T \; \}$$ 
is a $\pi$-system. Denote by $\cA$ the sigma-field
generated by this $\pi$-system. Let $\sigma \in T$. Deduce from Lemma
\ref{increa} that any connected component of $T \backslash \{ \sigma \}$ is in
$\cA$ (the details are left to the reader). Let $r>0$ and denote by
$\overline{B}(\sigma , r)$ the closed ball with radius $r$ and with center $\sigma$. Then, 
$$ T \backslash \overline{B}(\sigma , r)=\{ \sigma' \in T \; : \; d(\sigma , \sigma'
) >r \} . $$
Let $C$ be a connected component of $T \backslash \overline{B} (\sigma , r)$. There exists
$\sigma_0$ such that $d(\sigma, \sigma_0) =r$ and such that $C\cup \{
\sigma_0\}$ is the closure of $C$. It implies that $C$ is a connected
component of $T \backslash \{ \sigma_0\}$. Thus, any connected component $C$ of $T \backslash \overline{B}(\sigma , r)$
is in $\cA$. So is $T \backslash \overline{B}(\sigma , r)$, for any $\sigma
\in T$ and any $r>0$, which easily completes the proof of
the lemma. \cqfd

\section{Construction of the height function.}
\label{construc}

In this section we prove Theorem \ref{main} and Theorem \ref{timechange}. Let us consider 
a rooted compact real tree $(T, d, \rho) $ equipped with a compatible linear
order $\leq$ and a compatible measure $\mu$. 
To avoid trivialities, we assume that $T$ is not a point. Observe that (Mes) implies 

\begin{itemize} 
\item{{\bf (Inc)}}: $\; $ 
For any $\sigma <\sigma' $ we get  
$$ \mu ( L_{\sigma}) < \mu ( L_{\sigma'})  \; .$$
\end{itemize}

\vspace{4mm}
Let us set $ M =\mu (T)$. Observe that 
$$ M =\sup  \{ \mu(L_{\sigma}) \; , \; \sigma \in T\}. $$
Since $T\neq \{ \rho \}$, $M$ is positive. For any $t\in [0, M]$, we define 
$$ G_t =\{ \sigma \in T \; : \; \mu(L_{\sigma}) \leq t\} .$$
We use the following notation: 
$$ m_t= \sup_{ \sigma \in G_t } \mu (L_{\sigma})\quad {\rm and} \quad  M_t=
\inf_{ \sigma \in T \backslash G_t } \mu (L_{\sigma}) , $$
with the convention that $m_t=0$ if $G_t=\emptyset$ and that $M_t=M$ if
$G_t=T$. Observe that $G_t=T$ iff $t=M$. 
Clearly, $ m_t \leq t\leq M_t$. We also introduce 
$$ D_t= \bigcap_{ \sigma \in T \backslash G_t} L_{\sigma} \; , $$
if $t<M$ and $D_t=T$ if $t=M$. Here is the key lemma used in the proof of Theorem \ref{main}. 
\begin{lemma}
\label{keylem} Fix $t\in  [0, M]$. The following assertions are true.
\begin{itemize} 
\item{(i)} $\mu (G_t)=m_t$ and $\mu (D_t)=M_t$. 

\item{(ii)} There exist $\sigma_-$ and $\sigma_+$ such that 
$\sigma_+ \in \lgeo \rho ,  \sigma_- \rgeo$ and such that

   $\quad - \; $   for any
non-decreasing sequence $\sigma^{-}_n $ , $n\geq 1$, that satisfies 
$$\lim_{n\rightarrow \infty} \mu (L_{\sigma^{-}_n })=m_t $$
one has  
$$\lim_{n\rightarrow \infty} d(\sigma^{-}_n ,\sigma_-
)=0 \quad {\rm and} \quad  
G_t \backslash \{ \sigma_{-} \} \subset \bigcup_{n\geq 1} L_{\sigma^{-}_n}
\subset G_t  \; ;$$

 $\quad - \; $  for any decreasing sequence $\sigma^{+}_n $ , $n\geq 1$, that satisfies 
$$\lim_{n\rightarrow \infty} \mu (L_{\sigma^{+}_n })=M_t$$
one has  
$$\lim_{n\rightarrow \infty} d(\sigma^{+}_n ,\sigma_+
)=0 \quad {\rm and} \quad \bigcap_{n\geq 1} L_{\sigma^{+}_n}=D_t \; . $$


\item{(iii)} We also get 
$$ D_t \backslash \{ \sigma_{-} \} \subset G_t \subset D_t . $$
\end{itemize}
\end{lemma}

\begin{remark}
\label{attentionM}
There is no non-decreasing sequence $\sigma^{-}_n $ , $n\geq 1$ satisfying the
condition of $(ii)$ iff $G_t=\emptyset$. In this case, we set $\sigma_- =
\rho$, by convention. Note also that is no decreasing sequence $\sigma^{+}_n $ , $n\geq 1$ satisfying the
condition of $(ii)$ iff $M_t=M$. In this case, we set $\sigma_+ =\sigma_- $, by convention.
\end{remark}

\noindent
{\bf Proof}: First note that $G_t \subset D_t$. Suppose that there exists 
$\widetilde{\sigma} \in  D_t \backslash G_t$: then, by definition of $D_t$, for any $\sigma ' 
\in T \backslash G_t$, we get $\widetilde{\sigma} \leq \sigma '$. 
Consequently, $\widetilde{\sigma}$ is the minimal element of 
$ T \backslash G_t$. Observe next that (Inc) implies that 
$ \sigma <\widetilde{\sigma}$
for any $\sigma \in  G_t$. Thus we have proved that

-  either $T \backslash G_t$ has no minimal point and then $G_t =D_t$,

-  either $T \backslash G_t$ has a minimal point denoted by $\widetilde{\sigma} $ and then 
\begin{equation}
\label{prems}
G_t= L_{\widetilde{\sigma} }\backslash \{ \widetilde{\sigma} \} \quad 
{\rm and} \quad D_t = L_{\widetilde{\sigma} }.
\end{equation}
\noindent
We distinguish several cases in the proof:  

\noindent
$\bullet$ {\it Case 1}: Suppose that there exists $\sigma_* \in G_t$ such that 
$m_t=\mu ( L_{\sigma_*})$. Then, by (Inc) we get $G_t=L_{\sigma_*}$. We first claim that 
\begin{equation}
\label{deuz}
 D_t = G_t=L_{\sigma_*} \; . 
\end{equation}
Clearly $L_{\sigma_*}\subset D_t$. Suppose that there exists 
$\sigma' \in D_t \backslash L_{\sigma_*}$. Then, $\sigma'$ is the 
minimal element of $T\backslash G_t$. Then by (\ref{prems}), $G_t$ is not
compact, which rises a contradiction. Then it implies (\ref{deuz}) and $(i)$ and
$(iii)$ follow.

Let $\sigma^{-}_n $ , $n\geq 1$ be as in the lemma. By Lemma \ref{monot}, 
it has a limit that we denote by $\sigma_{-}$. We claim that 
\begin{equation}
\label{troiz}
 \sigma_{-}=\sigma_* \; .   
\end{equation}
Indeed, by Lemma \ref{increa}, $\sigma_{-}$ is in the closure of 
$\bigcup_{n\geq 1} L_{\sigma^{-}_n } $ and since 
$\sigma^{-}_n \in L_{\sigma_*}$ for any $n\geq 1$, we get 
$\sigma_{-}\in L_{\sigma_*}$. Suppose that we are in Case (I), Case (II) or Case (III) with 
$\sigma_- \notin  {\rm Br} (T)\cup \{ \rho \} $ of Lemma \ref{increa}. Then, the closure of 
$\bigcup_{n\geq 1} L_{\sigma^{-}_n } $ is exactly $L_{\sigma_-}$ and (Inc) implies that 
$$ m_t= \mu \left( \bigcup_{n\geq 1} L_{\sigma^{-}_n }  \right) 
\leq \mu ( L_{\sigma_-}) \leq \mu ( L_{\sigma_*} )=m_t .  $$
Thus, $m_t= \mu ( L_{\sigma_-})=\mu ( L_{\sigma_*})$ and $\sigma_{-}=\sigma_* $ by 
(Inc) again.

   Assume now that we are in Case (III) of Lemma \ref{increa} 
with $\sigma_- \in  {\rm Br} (T)\cup \{ \rho \} $. We keep the same
notation. We easily get
$$ \bigcup_{n\geq 1} L_{\sigma^{-}_n } =
\overline{\bigcup_{n\geq 1} L_{\sigma^{-}_n }  }= L_{\sigma_-} \cup \left( 
\bigcup_{i\in K} C_{\sigma_-}(i) \right) \subset L_{\sigma_*}. $$
Suppose that $\sigma_{-} \neq   \sigma_*$. Then $\sigma_{-} <   \sigma_*$. It
implies that $\sigma_* $ is in $T \backslash \lgeo \rho , \sigma_- \rgeo $. 
Let $C$ be the connected component of 
$T \backslash \lgeo \rho , \sigma_- \rgeo $ such that $\sigma_* \in
C$. Clearly, we get  
$$ L_{\sigma_-} \cup \left( 
\bigcup_{i\in K} C_{\sigma_-}(i) \right) \, \leq  \, C \; . $$
Observe that it is always possible to find $\sigma'\in C$ such that $\sigma' < \sigma_*$. So we get  
$$ \bigcup_{n\geq 1} L_{\sigma^{-}_n } =  L_{\sigma_-} \cup \left(  
\bigcup_{i\in K} C_{\sigma_-}(i) \right) \subset  L_{\sigma'} 
\subset_{\neq } L_{\sigma_*} \; .$$ 
But (Inc) implies 
$$ m_t= \mu \left( \bigcup_{n\geq 1} L_{\sigma^{-}_n } \right) \leq 
\mu ( L_{\sigma'} ) < \mu ( L_{\sigma_*} )=m_t ,$$ 
which is impossible. Consequently, (\ref{troiz}) holds and we also have proved that
\begin{equation}
\label{quatrz}
 L_{\sigma_-} \backslash \{ \sigma_- \} \subset 
\bigcup_{n\geq 1} L_{\sigma^{-}_n } \subset  L_{\sigma_-} =G_t .
\end{equation}
Observe that $\sigma_* =\sigma_-$ does not depend 
on any sequence $\sigma^{-}_n $ , $n\geq 1$ satisfying the assumptions of 
the lemma. Consequently, (\ref{quatrz}) remains true for any such sequence. We next claim that
\begin{equation}
\label{quatrzbis}
 \forall \sigma' \in T \backslash G_t \; , \quad \mu (L_{\sigma'})\, > \, M_t \; .
\end{equation}
Indeed, suppose that there is 
$\sigma' \in T \backslash G_t$ such that $\mu (L_{\sigma'})=M_t$. It implies 
$$M_t= \mu (L_{\sigma'}) >t \geq m_t= \mu (L_{\sigma_-}) \; . $$
Thus, $\sigma_- <\sigma '$ by (Inc). But we can always find 
$\sigma ''$ such that $\sigma_- <\sigma '' <\sigma ' $. Since 
$G_t=L_{\sigma_-}$, it implies that $\sigma'' \in T \backslash G_t$ and then 
by (Inc)
$$ M_t \leq \mu (L_{\sigma'' }) < \mu (L_{\sigma' })=M_t , $$
which is impossible. Therefore (\ref{quatrzbis}) holds.

\vspace{3mm}
  
   Let $\sigma^{+}_n $ , $n\geq 1$, satisfying the assumptions of 
the lemma. By Lemma \ref{decrea}, it has a limit that we denote by 
$\sigma_+$. Let $\sigma^{,+}_n $ , $n\geq 1$, satisfying the same 
assumptions. (\ref{quatrzbis}) implies that for any $n\geq 1$, 
$$\mu ( L_{\sigma^{+}_n }) \;{\rm and} \;  \mu (L_{\sigma^{,+}_n }) >M_t \quad {\rm
  and} \quad \lim_{n\rightarrow \infty }\mu ( L_{\sigma^{+}_n }) =  
\lim_{n\rightarrow \infty }\mu ( L_{\sigma^{,+}_n }) = M_t .$$ 
So, by (Inc) we can construct a sequence 
 $\sigma^{,,+}_n $ , $n\geq 1$, that also satisfies the assumptions of 
the lemma and such that it contains an infinite number of terms of the 
two sequences $\sigma^{+}_n $ , $n\geq 1$, and $\sigma^{,+}_n $ , $n\geq 1$. 
Lemma \ref{decrea} implies that 
$\sigma^{,,+}_n $ , $n\geq 1$, is convergent. Therefore the limit of
$\sigma^{,+}_n $ , $n\geq 1$ has to be also $\sigma_+$. 
Thus, $\sigma_+$ does not depend on a choice of a sequence satisfying the assumptions of the lemma.

   Let us fix such a sequence $ \sigma^{+}_n $ , $n\geq 1$. We claim that 
\begin{equation}
\label{cinqz}
D_t = \bigcap_{n\geq 1} L_{\sigma^{+}_n } .
\end{equation}
Indeed, observe first that 
$$ D_t \subset \bigcap_{n\geq 1} L_{\sigma^{+}_n } .$$
Suppose that there is $\sigma' \in 
\bigcap_{n\geq 1} L_{\sigma^{+}_n } \backslash D_t$. Then, we have 
$\mu ( L_{\sigma^{+}_n }) \geq \mu ( L_{\sigma'})$ for any 
$n \geq 1$. It implies that 
$M_t  \geq \mu ( L_{\sigma'})$, which actually implies $M_t  > \mu ( L_{\sigma'})$
by (\ref{quatrzbis}). Thus, $\sigma' $ is in $G_t$. But, it implies
$$  L_{\sigma'} \subset \bigcap_{T \backslash G_t}  L_{\sigma}  =D_t  , $$
which rises a contradiction. Consequently (\ref{cinqz}) holds true.

   To complete the proof of the lemma in the first case, it remains to prove 
that $\sigma_+ \in \lgeo \rho , \sigma_- \rgeo$: we clearly have 
$$  L_{\sigma_+} \subset 
\bigcap_{n\geq 1} L_{\sigma^{+}_n } =D_t .$$
Thus, $\mu (L_{\sigma_+}) \leq M_t$, which actually implies $\mu (L_{\sigma_+}) <
M_t$ by (\ref{quatrzbis}) and we get $\sigma_+ \in 
G_t= L_{\sigma_-}$.  Observe now that $\sigma_+ \in 
\overline{T\backslash L_{\sigma_-}}$ as the limit of the $\sigma^{+}_n $'s. Thus, 
$$ \sigma_+ \in L_{\sigma_-} \cap 
\overline{ T \backslash L_{\sigma_-}} \subset  
\lgeo \rho , \sigma_- \rgeo , $$
which completes the proof of the lemma in Case 1. 

\vspace{4mm} 

$\bullet$ {\it Case 2}: We now suppose 
\begin{equation}
\label{bisz}
 \forall \sigma' \in G_t \; , \quad \mu (L_{\sigma'}) <m_t \; .
\end{equation}
\noindent
By arguments similar to those used to prove uniqueness for $\sigma_{+}$ in Case 1, we prove that there exists 
$\sigma_{-} \in T$ such that 
any non-decreasing sequence $\sigma^{-}_n $ , $n\geq 1$, that satisfies the assumptions of 
the lemma converges to $\sigma_{-}$. 

    Consider such a sequence $\sigma^{-}_n $ , $n\geq 1$, and note that 
$$ \bigcup_{n\geq 1} L_{\sigma^{-}_n } \subset G_t \; . $$
Suppose there exists $\sigma' \in  G_t \backslash 
\bigcup_{n\geq 1} L_{\sigma^{-}_n } \; $ . Then for any $n \geq 1$, $\mu (L_{\sigma'}) \geq 
\mu (L_{\sigma_n^{-}} )$ and we get $\mu (L_{\sigma'}) \geq m_t$, which contradicts (\ref{bisz}). Consequently, 
\begin{equation}
\label{sixz}
\bigcup_{n\geq 1} L_{\sigma^{-}_n } = G_t \; .
\end{equation}
We distinguish two subcases whether it exists $\sigma_* \in T \backslash 
G_t $ such that $\mu (L_{\sigma_*})=M_t$ or not. 

\vspace{4mm}

${\bf -}$ {\it Case 2.1}: Suppose there exists such a 
$\sigma_* \in T \backslash G_t $. Then $\sigma_*$ is the minimal element of $T
\backslash G_t$ and by (\ref{prems}) we get 
\begin{equation}
\label{septz}
D_t= L_{\sigma_*} \quad {\rm and} \quad G_t= L_{\sigma_*} \backslash  \{ 
\sigma_*\} \; 
\end{equation}
and we get 
\begin{equation}
\label{huitz}
 \widetilde{\sigma}=\sigma_*=\sigma_- \; .
\end{equation}
Assume that there exists $\sigma^{+}_n $ , $n\geq 1$, a decreasing sequence 
satisfying assumptions of 
the lemma. It has a limit denoted by $\sigma_{+}$. By previously used arguments, we can prove that any 
decreasing sequence satisfying assumptions of the lemma converges to $\sigma_{+}$. 
Recall that we have $D_t=L_{\sigma_* }=L_{\sigma_- }$. Suppose that there is 
$\sigma' \in \bigcap_{n\geq 1} L_{\sigma_n^{+}} \backslash D_t $. Then we get by (Inc) 
$$ M_t =\mu (L_{\sigma_- }) < \mu (L_{\sigma' }) \leq  \mu \left( \bigcap_{n\geq 1}
L_{\sigma_n^{+}} \right) =M_t , $$
which is absurd. It implies $D_t=L_{\sigma_* }=L_{\sigma_- }=\bigcap_{n\geq
  0} L_{\sigma_n^{+}} $. By (\ref{septz}) and by the form of $\bigcap_{n\geq 1} L_{\sigma_n^{+}}$ given by Lemma
\ref{decrea}, either $\sigma_-= \sigma_+$, either $\sigma_-$ is in a connected component of 
$T\backslash \{ \sigma_+ \} $ that does not contain the root. Thus, it shows 
$$  \sigma_{+} \in \lgeo \rho , \sigma_- \rgeo \; ,$$
which completes the proof of the lemma in Case 2.1.  

\vspace{4mm}

${\bf -}$ {\it Case 2.2}: We suppose that (\ref{quatrzbis}) holds. 
By arguments similar to those used previously, we prove that there exists 
$\sigma_+ \in T$ such that 
any decreasing sequence 
$\sigma^{+}_n $ , $n\geq 1$, that satisfies the assumptions of 
the lemma converges to $\sigma_{+}$. 
Consider such a sequence $\sigma^{+}_n $ , $n\geq 1$. We claim that 
\begin{equation}
\label{dixz}
D_t = \bigcap_{n\geq 1} 
L_{\sigma_n^{+}} =G_t \; .
\end{equation}
First note that $ D_t \subset \bigcap_{n\geq 1} L_{\sigma_n^{+}} $. If  
$\sigma \in  
\bigcap_{n\geq 1} L_{\sigma_n^{+}} $, then  
$\mu ( L_{\sigma}) \leq \mu ( L_{\sigma_n^{+}})$ for any $n\geq 1$. So 
$\mu ( L_{\sigma}) \leq M_t$ and thus, $\sigma \in G_t$ by (\ref{quatrzbis}).
But clearly $G_t \subset D_t$, which completes the proof of(\ref{dixz}) .

  It remains to prove 
 \begin{equation}
\label{onzez}
\sigma_+ \in \lgeo \rho , \sigma_-\rgeo \; .
\end{equation}
First, observe that by (\ref{dixz}), $\sigma_+ \in G_t$. But it is the limit 
of the $\sigma_n^{+}$'s that are in $T \backslash G_t$. Thus, 
$$ \sigma_+ \in \overline{T \backslash G_t } \cap G_t \;  .$$
Recall (\ref{sixz}). Then by  Lemma
\ref{increa}. Deduce that  
$$\overline{T \backslash G_t } \cap G_t  \subset 
\lgeo \rho , \sigma_-\rgeo  \; ,$$
which implies (\ref{onzez}) and the proof of the lemma is completed. \cqfd

\vspace{4mm}

%

\begin{definition}
\label{exploration}(Exploration mapping associated with $\mu $) 
For any $t\in [0, M]$, we set $\phi (t)=\sigma_-$, where $\sigma_-$ 
is defined by Lemma \ref{keylem} and Remark \ref{attentionM}. 
\end{definition}

\begin{lemma}
\label{regul}
The exploration mapping $\phi$ is left-continuous with right-limits. Moreover, 
$\phi(0+)=\rho$ and for any $t\in [0, M)$, $\phi (t+) \in \lgeo \rho , \phi (t)\rgeo $. 
\end{lemma}
\noindent
{\bf Proof}: Fix $t\in (0, M]$. 
We first prove that $\phi$ is left-continuous at $t$.  
We define $t_0$ by 

$$ t_0 = \sup \; \{ \mu (L_{\sigma }) \; , \; \sigma \in T \; : \;  
\mu (L_{\sigma }) <t \} \; . $$
Recall that 
$$ m_t = \sup \; \{ \mu (L_{\sigma }) \; , \; \sigma \in T \; : \;  
\mu (L_{\sigma }) \leq t \} \; . $$
We distinguish three cases: 

- {\it Case 1}: $t_0 <m_t$. Then, clearly $t=m_t =\mu (L_{\phi (t)})$ and $G_t=L_{\phi(t)}$. Next observe 
that for any $s\in [t_0, t)$, $G_s=L_{\phi (t)} \backslash \{ \phi (t)\}$ and 
thus $\phi (s)=\phi(t)$, which implies that $\phi$ is left-continuous at $t$.

\vspace{3mm}

- {\it Case 2}: $t_0 =m_t <t$. If $s\in [m_t, t]$, then $m_s=m_t$ and $\phi(s)= \phi(t)$, which also implies that 
$\phi $ is left-continuous at $t$.  

\vspace{3mm}

- {\it Case 3}: $t_0=m_t=t$.  
Thus we can find an {\it increasing} sequence 
$\sigma^{-}_n $ , $n\geq 1$, converging to $\sigma_-=\phi (t)$ and such that 
$\lim_{n \rightarrow \infty} \mu (L_{\sigma_n^{-}} )= t $. Let $t_k$ , 
$k\geq 1$, be any increasing sequence of $[0, M]$ converging to $t$. We first claim 
that 
\begin{equation}
\label{croiss}
\lim_{k \rightarrow \infty} m_{t_k}=t \; . 
\end{equation}
Clearly $m_{t_k} \leq t_k <t$. Set $t_n^-= \mu (L_{\sigma_n^{-}} )$. By definition, 
$t_n^-=m_{t_n^-} $. Since $t_n^- <t$, then for any $n\geq 1$, we can find $k_n\geq 1$, such 
that $t_n^- \leq t_{k_n}$, 
which implies that $m_{t_n^-}=t_n^- \leq m_{t_{k_n}}$. Thus, (\ref{croiss})
follows since $\lim_{n \rightarrow \infty} m_{t_n^-}=m_t=t$.

   Without loss of generality we can assume that the sequence 
$m_{t_k}$ , $k\geq 1$, is increasing. 
Use Lemma \ref{keylem} at each $t_k$ to find a sequence $\gamma_k$ , $k\geq 1$, 
such that 
\begin{equation}
\label{sousapprox}
 d(\gamma_k , \phi (t_k) )\leq 2^{-k}\quad {\rm and}  \quad 
0\leq m_{t_k} -\mu (L_{\gamma_k}) 
 < 2^{-k}  \wedge \left( m_{t_k} -\mu (L_{\gamma_{k-1}})   \right)   
\end{equation}
(observe that strict inequality is possible because $m_{t_{k-1}}
<m_{t_k}$). Then, we get 
$$\mu (L_{\gamma_{k-1}} ) <\mu (L_{\gamma_{k}} )  \leq m_t=t \; , $$
which implies that $\gamma_k$ , $ k\geq 1$, is an increasing sequence of 
$G_t$. Moreover by (\ref{croiss}), 
we get $\lim_{k \rightarrow \infty} \mu (L_{\gamma_{k}} ) =t=m_t$. Lemma
\ref{keylem} implies that  
$\lim_{k \rightarrow \infty}d(\gamma_{k} , \phi(t))=0$ and, by definition of 
the sequence $\gamma_k$ , $ k\geq 1$, it implies that 
$\lim_{k \rightarrow \infty}d(\phi(t_k), \phi(t))=0$. This proves that 
$\phi$ is left-continuous at $t$ in Case 3. 

\vspace{4mm}

Existence of right-limits of $\phi $ at $t\in [0, M)$ is treated similarly: if 
$t< M_t$, then for any $s \in [t, M_t)$, we clearly have $m_s=m_t$, $M_s=M_t$ and 
$\phi ( s)=\phi (t)$. Thus, $\phi $ has a right-limit at $t$, which is $\phi (t)$.

   If we now assume that $M_t=t$, then we can find a {\it decreasing} sequence 
$\sigma^{+}_n $ , $n\geq 1$, converging to $\sigma_+$ and such that 
$\lim_{n \rightarrow \infty} \mu (L_{\sigma_n^{+}} )= M_t= t $. 
Let $t_k$ , $k\geq 1$, be any decreasing sequence of $(0, M]$ converging to $t$. 
Set $t_n^+ = \mu (L_{\sigma_n^{+}} )$ , $n \geq 1$. 
We can find two increasing 
subsequences of indices $n(1,k)$ , $n(2,k)$, $k\geq 1$ such that 
$$ t_{n(1,k+1)}^+ = \mu (L_{\sigma_{n(1,k+1)}^{+}} ) < t_{n(2,k)} < t_{n(1,k)}^+ = \mu (L_{\sigma_{n(1,k)}^{+}} )  \; . $$
Consequently
\begin{equation}
\label{entrelac}
 t_{n(1,k+1)}^+ \leq m_{t_{n(2,k)}}\leq M_{t_{n(2,k)}} \leq t_{n(1,k)}^+  \; . 
\end{equation}
It implies that $m_{t_k}$, $k\geq 1$ converges to $t$. Without loss of
generality we can assume that $m_{t_k}$, $k\geq 1$ is a decreasing 
sequence. Use Lemma \ref{keylem} at each $t_k$ to find a sequence $\gamma_k$ , $k\geq 1$ 
such that 
\begin{equation}
\label{sousapprox2}
 d(\gamma_k , \phi (t_k) )\leq 2^{-k}\quad {\rm and}  \quad 
0\leq m_{t_k} -\mu (L_{\gamma_k}) < 2^{-k} \wedge \left(\mu (L_{\gamma_{k-1}}) -m_{t_k}  \right)  \; .  
\end{equation}
Then, by (\ref{entrelac}) and (\ref{sousapprox2})
$$   t < m_{t_{k+1}} <\mu (L_{\gamma_{k}} ) < \mu (L_{\gamma_{k-1}} ) \leq m_{t_{k-1}} \; , $$
which implies that $\gamma_k$ , $ k\geq 1$ is a decreasing sequence of 
$T\backslash G_t$ such that $\lim_{k \rightarrow \infty} \mu (L_{\gamma_{k}} ) = M_t =t$.
Lemma \ref{keylem} implies 
$\lim_{k \rightarrow \infty}d(\gamma_{k} , \sigma_+)=0$. By definition of 
the sequence $\gamma_k$ , $ k\geq 1$, it implies that 
$\lim_{k \rightarrow \infty}d(\phi(t_k), \sigma_+)=0$. 
This proves that $\phi$ has a right-limit at $t$ and also that $\phi (t+)=\sigma_+$, where $\sigma_+$ is the point 
associated with $t$ as defined in Lemma \ref{keylem}.

\vspace{4mm}

  It remains to prove that $\phi$ is right-continuous at $0$. If $\mu (\{ \rho \}) >0$, then $M_0 >0=m_0$ and we are in Case 1 or in 
Case 2. Assume that  $\mu (\{ \rho \}) =0$. Fix a sequence 
$t_k $, $k \geq 1$ that decreases to $0$. Let $\sigma_n$, 
$n\geq 1$, be a decreasing sequence of $T$ such that $\{ \rho\}= 
\bigcap_{n\geq 1} L_{\sigma_n}$. (Inc) implies that 
$$ \lim_{n \rightarrow \infty} \mu (  L_{\sigma_n})=\mu (\{ \rho\})=0 \; .$$
Let $n \geq 1$. For all sufficiently large $k$ we get 
$t_k <  \mu (  L_{\sigma_n})$, which implies that 
$\phi(t_k) \in D_{t_k} \subset  L_{\sigma_n}$. Then any 
limit point $\gamma$ of the sequence $\phi (t_k)$ , $k\geq 1$, is 
in $L_{\sigma_n}$, for any $n\geq 1$. This implies that $\rho$ 
is the only limit point of the sequence $\phi (t_k)$ , $k\geq 1$ and the proof of the 
lemma is now completed. \cqfd

\vspace{4mm}

Let us set 
$$ h(t) = d(\rho , \phi (t)) \; , \quad t\in [0, M] .$$
Clearly $h$ is left-continuous with right-limit; we also have 
$h (0)=h(0+)=0$. Recall that $\phi (t)=\sigma_-$ and that if $\phi (t)\neq \phi (t+)$, then 
$\phi (t+)=\sigma_+$ with the notation of 
Lemma \ref{keylem}. Since $\sigma_+ \in \lgeo \rho , 
\sigma_-\rgeo$, we get $h(t+) \leq h(t)$ (note that if $\phi (t)= \phi (t+)$,
then $\phi (t+)$ is not necessarily equal to $\sigma_+ $). Thus, $h$ is in 
$\cH_M$. 

\begin{proposition}
\label{g-isome} There exits an isometry $\jmath_h$ from $(T_h,d_h)$ 
onto $(T,d)$ such that $\jmath (\rho_h)=\rho $ and such that 
$$ \xi_1 \leq _h \xi_2 \Longrightarrow \jmath_h (\xi_1) \leq 
 \jmath_h (\xi_2) \; .$$
\end{proposition}
\noindent
{\bf Proof}: We first claim that for any $t_1 <t_2<t_3$ in $[0, M]$, 
\begin{equation}
\label{supe}
d(\rho , \phi (t_2)) \geq d(\rho , \phi (t_1) \wedge  \phi (t_3)) \; .
\end{equation}
First observe that if $t_2\in [m_{t_1}, M_{t_1})$, then clearly 
$\phi (t_2)=\phi (t_1)$. By left-continuity, we also get $\phi (M_{t_1})=\phi (t_1)$. 
Similarly, if $t_2\in [m_{t_3}, M_{t_3}]$, then 
$\phi (t_2)=\phi (t_3)$ . Consequently, (\ref{supe}) holds for any 
$t_2\in [m_{t_1}, M_{t_1}] \cup [m_{t_3}, M_{t_3}]$.

   Let us assume that $M_{t_1} < t_2 <m_{t_3}$, which implies 
\begin{equation}
\label{interw}
m_{t_1} <m_{t_2}<m_{t_3} \; .
\end{equation}
By Lemma \ref{keylem}, we can find three non-decreasing sequences 
$\sigma_n^{-} (i)$ , $n\geq 1$, $i\in \{ 1,2,3 \}$ 
such that 
$$ \sigma_n^{-} (i) \in G_{t_i} \; , \quad 
\lim_{n\rightarrow \infty } d( \sigma_n^{-} (i) , \phi (t_i) )=0 \; , 
\quad \lim_{n\rightarrow \infty } \mu (L_{\sigma_n^{-} (i) } )= m_{t_i} \; .$$
Inequality (\ref{interw}) implies that for all sufficiently large 
$n$, $\mu (L_{\sigma_n^{-} (1) }) < \mu (L_{\sigma_n^{-} (2) }) 
< \mu (L_{\sigma_n^{-} (3) })$ and by (Inc) 
\begin{equation}
\label{interw2}
\sigma_n^{-} (1) < \sigma_n^{-} (2) <\sigma_n^{-} (3) \; .
\end{equation}
Set $\sigma_0= \phi(t_1) \wedge \phi(t_3) $ and let $\gamma $ be such that 
$$ \lgeo \rho , \gamma \rgeo = \lgeo \rho , \phi (t_2) \rgeo 
\cap \left( \lgeo \rho , \phi (t_1) \rgeo \cup  \lgeo \rho , \phi (t_3) \rgeo \right) .$$
Suppose that $\gamma \in \lgeo 
\rho , \sigma_0 \lgeo $. Then $\phi(t_2) \wedge \phi(t_3) 
=\phi(t_2) \wedge \phi(t_1) =\gamma $. Then $\phi (t_2) \wedge \sigma_0 =
\gamma $. 
Now observe that for any $\sigma \in T$ 
$$ 2d(\sigma \wedge \sigma_0 , \sigma_0  ) = d(\sigma_0 , \sigma )+ 
d(\rho , \sigma_0 ) -d(\sigma , \rho ) \; . $$
Thus, the application $\sigma \rightarrow d(\sigma 
\wedge \sigma_0  , \sigma_0 )$  is continuous.
Consequently, for all sufficiently large $n$ 
\begin{equation}
\label{excluss}
 d(\sigma_n^{-}(2) \wedge \sigma_0 , \sigma_0) > 
\frac{2}{3} d(\gamma , \sigma_0) \quad {\rm and} \quad 
 d(\sigma_n^{-}(i) \wedge \sigma_0 , \sigma_0) < 
\frac{1}{3} d(\gamma , \sigma_0)  \; , \; i\in \{ 1,3 \} \; .
\end{equation}
Let $\gamma_n $ be such that 
$$ \lgeo \rho , \gamma_n \rgeo = \lgeo \rho , \sigma_n^{-}(1) \rgeo 
\cap \left( \lgeo \rho , \sigma_n^{-}(2) \rgeo \cup  
\lgeo \rho ,\sigma_n^{-}(3)  \rgeo \right) .$$
(\ref{excluss}) implies that for all sufficiently large $n$, the point 
$\gamma_n $ is not in $\lgeo \rho , \sigma_n^{-}(2) \rgeo  $, which 
contradicts (\ref{interw2}) by (Or2). Then, $\gamma \notin \lgeo 
\rho , \sigma_0 \lgeo $ and $\gamma$ is necessarily  in $\lgeo \phi (t_1) , \phi (t_3) \rgeo $. Then, we get 
$$ d(\rho , \phi (t_2)) = d(\rho , \sigma_0) +d(\gamma , \sigma_0 ) 
+  d(\gamma , \phi (t_2)) , $$
which implies (\ref{supe}).

We keep notation $\sigma_0 = \phi (t_1) \wedge  \phi (t_3)$  and we now prove that 
\begin{equation}
\label{infatteint}
\inf_{t\in [t_1 , t_3]} d( \rho , \phi (t)) = d(\rho , \sigma_0) \; .
\end{equation}
To avoid triviality we suppose that $\sigma_0 \notin \{ 
\phi(t_1) ,  \phi(t_3) \}$. By Lemma \ref{keylem} and the form of $D_{t_1}$ and $G_{t_3}$ given by 
Lemmas \ref{increa} and 
\ref{decrea} we get $D_{t_1} \subset G_{t_3}$, which implies that  
$\phi (t_1)< \phi (t_3)$.  Let $\sigma_n $, $n\geq 1$, be a 
sequence in $\rgeo \sigma_0 , \phi(t_3) \rgeo $ that decreases to 
$\sigma_0$. By 
Lemma \ref{keylem} and the form of $D_{t_1}$ and $G_{t_3}$ given by 
Lemma \ref{increa} and Lemma \ref{decrea}, we get for all 
$n\geq 1$, 
$$ D_{t_1} \subset_{\neq} L_{\sigma_n}  \subset_{\neq} G_{t_3} \;  .$$
Set $s_n=\mu (L_{\sigma_n})$ , $n\geq 1$. The previous observation implies that for any 
$n\geq 1$,  
$$ t_1 \leq M_{t_1} \leq s_n \leq m_{t_3} \leq t_{3} \; .$$
Set $t=\mu (\bigcap_{n\geq 1} L_{\sigma_n})=\lim_{n\rightarrow \infty } s_n $. 
Then $t\in [t_1, t_3]$. Next, by definition of $\phi$, we have  
$\phi (s_n)=\sigma_n$. Since the sequence $s_n$ , 
$n\geq 1$, decreases to $t$, we get $\phi (t+)=\sigma_0$. This, combined 
with (\ref{supe}) implies (\ref{infatteint}).

   Now observe that (\ref{infatteint}) easily implies that for any $s ,t  \in [0, M]$,
\begin{equation}
\label{isoso}
d(\phi (s) , \phi (t) )= h(s) +h(t) -2 
\inf_{u\in [s\wedge t ,s\vee t ]} h(u) =d_h(s,t)\; .
\end{equation}
Recall that $p_h$ stands for the canonical projection from $[0, M]$ to $T_h$. 
It implies that if $p_h(s)=p_h(t)$, then $\phi (s)=\phi (t)$. Thus, it makes 
sense to define $\jmath_h: T_h \rightarrow T$ by 
$\jmath_h (\xi)=\phi (s)$ for any $s\in p_h^{-1}(\{ \xi \})$. Then 
(\ref{isoso}) implies that $\jmath_h$ is an isometry from $(T_h, d_h)$ 
onto $(T,d)$. Moreover we get $\jmath_h ( \rho_h)=\rho$.

It remains to prove that $\jmath_h $ is increasing: let $ \sigma \in T \backslash \{ \rho \}$. It is always possible 
to find an increasing sequence $\sigma_n^-$ , $n\geq 1$, that converges to $\sigma$. 
Lemma \ref{increa} Case (I) implies that 
$$ \bigcup_{n\geq 1} L_{\sigma_n^-} = L_{\sigma}  \backslash \{ \sigma \} $$
and thus  
$$ \lim_{n\rightarrow \infty} \mu \left(  L_{\sigma_n^-} \right) = \mu \left(   L_{\sigma}  \backslash \{ \sigma \} \right)  \; . $$
Lemma \ref{keylem} implies that $\phi (\mu (L_{\sigma}  \backslash \{ \sigma \} ) )= \sigma$. Now observe 
that if $\phi (t)=\sigma$, then Lemma \ref{keylem} easily implies that 
$$  L_{\sigma}  \backslash \{ \sigma \} \, \subset \, G_t  \; . $$
Thus,
$$ \mu \left(   L_{\sigma}  \backslash \{ \sigma \} \right) \leq  \mu ( G_t) =m_t \leq t \; . $$
It proves that for any $\sigma \in T$
\begin{equation}
\label{debutprems}
 \inf \phi^{-1} (\{ \sigma \}) =\mu (L_{\sigma} \backslash \{ \sigma \}) \; .
\end{equation}
Consequently, if $\sigma_1 <\sigma_2$, then 
$$\inf \phi^{-1} (\{ \sigma_1 \}) =\mu (L_{\sigma_1} \backslash \{ \sigma_1 \}) \leq 
\mu (L_{\sigma_1}) \leq \mu (L_{\sigma_2 } \backslash \{ \sigma_2 \})=\inf \phi^{-1} (\{ \sigma_2 \}) \; .$$
Apply this inequality to $\sigma_1 = \jmath_h (\xi_1)$ and $\sigma_2 =\jmath_h (\xi_2)$ and observe that 
$$ \inf \phi^{-1} (\{ \sigma_i \}) = \inf p_h^{-1} (\{ \xi_i \}) \; , \quad i\in \{ 1,2 \} $$
to complete the proof of the proposition.  \cqfd 

\vspace{4mm}

   We next prove the following proposition. 
\begin{proposition}
\label{mesurtheo} We have $\mu = \mu_h \circ \jmath_h ^{-1}$. 
Furthermore the function $h$ satisfies (Min). 
\end{proposition}
\noindent
{\bf Proof}: We first introduce some notation. Fix $\sigma \in T$ and recall
notation $\cC_{\sigma}$ and $\cC_{\sigma}^{+}$ from Section \ref{topol}. By Lemmas \ref{compoorder} and 
\ref{echelon}, all the connected components in $\cC_{\sigma} $ and $\cC_{\sigma}^{+}$  can be ordered by 
$\leq$. Define the following collection $\Sigma_{\sigma}$ 
of families of connected components: 
$$ \Sigma_{\sigma}= \left\{ \; S \subset \cC_{\sigma} \cup \cC_{\sigma}^{+} \; : \quad
\forall  C\in S \; , \;  \forall  C' \in \cC_{\sigma} \cup \cC_{\sigma}^{+} \;
: \; (C' \leq C) \Longrightarrow   (C'\in S ) \; \right\}  \; .$$
We next define the two following sets of real numbers:
$$ A_{\sigma}= \left\{ \mu (L_{\sigma }) + 
\sum_{C\in S} \mu (C) \; , \; S \in \Sigma_{\sigma} \right\} \quad {\rm and} \quad 
B_{\sigma}= \left\{ \mu (L_{\sigma }) + \sum_{C\in S \cap \cC_{\sigma} } 
\mu (C) \; , \; S \in \Sigma_{\sigma} \right\}  \; . $$
We first prove the following lemma
\begin{lemma}
\label{qdmeme} For any $\sigma \in T$, one has 
\begin{equation}
\label{preee}
F(h)\cap \phi^{-1} (\{ \sigma \})=[\mu (L_{\sigma } \backslash \{ \sigma \}), \mu (L_{\sigma })] \; ,
\end{equation}

\begin{equation}
\label{sesecond}
S(h)\cap \phi^{-1} (\{ \sigma \}) \subset B_{\sigma} \backslash \{  \mu (L_{\sigma })
\}\; , 
\end{equation}
and 
\begin{equation}
\label{secondlinea}
\phi^{-1} (\lgeo \rho ,\sigma \rgeo) \cap (\mu (L_{\sigma }), M]  \subset  (A_{\sigma}
\backslash \{  \mu (L_{\sigma }) \}) \; \cap S(h) \; .
\end{equation}
\end{lemma}
\noindent
{\bf Proof}: For any $t\in [\mu (L_{\sigma } \backslash \{ \sigma \}), \mu (L_{\sigma })\, )$, we get 
$m_t=\mu (L_{\sigma } \backslash \{ \sigma \})$ and $M_t=\mu (L_{\sigma })$; the definition of 
$\phi$ and (\ref{debutprems}) imply that 
\begin{equation}
\label{incluetap}
[\mu (L_{\sigma } \backslash \{ \sigma \}), \mu (L_{\sigma })] \subset F(h)\cap \phi^{-1} (\{ \sigma \}) \; .
\end{equation}
Let $t\in \phi^{-1} (\{ \sigma \}) \cap ( \, \mu (L_{\sigma }) \, , \, M]$. Thus 
\begin{equation}
\label{ineqlatter}
\mu (L_{\sigma }) \leq m_t \leq t \leq M_t \; .
\end{equation}
Let $\sigma_n^{-}$, $n\geq 1$ be a non-decreasing sequence such that 
$$ \lim_{n\rightarrow \infty} \mu (L_{\sigma_n^{-}}) = m_t \; .$$
Then Lemma \ref{keylem} implies that $ \lim_{n\rightarrow \infty} d(\sigma_n^{-} , \sigma ) =0$ and 
\begin{equation}
\label{inclulatter}
D_t\backslash \{ \sigma \}=
G_t\backslash \{ \sigma \} \subset \bigcup_{n\geq 1}L_{\sigma_n^{-}} \subset G_t \subset D_t \; .
\end{equation}
Suppose that the sequence $\sigma_n^{-}$, $n\geq 1$ corresponds to Case (I) or
Case (II) in Lemma \ref{increa}. Then 
(\ref{inclulatter}) implies 
$$G_t\backslash \{ \sigma \}= D_t\backslash \{ \sigma \}
=L_{\sigma} \backslash \{ \sigma \} \; .$$
It implies that $M_t= \mu (D_t) \leq \mu (L_{\sigma} )$, 
which contradicts (\ref{ineqlatter}). Then, 
$\sigma_n^{-}$, $n\geq 1$ corresponds to Case (III) in Lemma \ref{increa}. Consequently 
\begin{equation}
\label{equallatter}
D_t=G_t = L_{\sigma} \cup \left( \bigcup_{k\in K} C_{\sigma } (k)\right) \; ,
\end{equation}
with the same definition of $K\subset I_{\sigma}$ as in Lemma \ref{increa}. Thus it implies 
$$ m_t = \mu (G_t)= \mu (D_t)=M_t =t \; .$$
Set $S=\{ C_{\sigma } (k)\; , \; k\in K \}$. The definition of $K$ in  Lemma \ref{increa} implies that 
$S\in \Sigma_{\sigma}$. Thus, 
$$ t= \mu (G_t)= \mu (L_{\sigma}) + \sum_{k\in K} \mu (C_{\sigma } (k)) \; \in \; B_{\sigma} 
\backslash \{ \mu (L_{\sigma }) \}  \; .$$
Let $k_0 \in K$ and $\sigma' \in C_{\sigma} (k_0)$. Observe that 
$$ L_{\sigma} \subset L_{\sigma'} \subset G_t \; .$$
Set $t'=\mu (L_{\sigma'})$. Then $\mu (L_{\sigma}) <t'\leq t$ and $\phi (t')=\sigma' \neq \sigma$. It 
implies that $t \in S(h)$. It completes the proof of (\ref{preee}) and (\ref{sesecond}). 

\vspace{3mm}

    Let us prove (\ref{secondlinea}). First note that 
$$\phi^{-1}( \lgeo \rho , \sigma \lgeo ) \cap (\mu (L_{\sigma }) \, , \, M] \subset S(h) $$
and (\ref{sesecond}) implies that 
$\phi^{-1}( \{ \sigma \} )\cap (\mu (L_{\sigma }) \, , \, M] \subset S(h) \; .$
Thus 
$$ \phi^{-1}( \lgeo \rho , \sigma \rgeo ) \cap (\mu (L_{\sigma }) \, , \, M] \subset S(h) \; .$$
Let $t\in (\mu (L_{\sigma }) \, , \, M] $ be such that 
$\jmath_h (t)=\sigma' \in \lgeo \rho , \sigma \lgeo $. Then 
$t\in \phi^{-1}( \{ \sigma' \} ) \cap  S(h) $. By (\ref{sesecond}), it implies that 
$t \in B_{\sigma'} \backslash \{  \mu (L_{\sigma' }) \}$. Then there exists $S'\in \Sigma_{\sigma'}$ such that 
\begin{equation}
\label{lattersum}
t=\mu (L_{\sigma'}) + \sum_{C' \in S'} \mu (C') \; .
\end{equation}
Since $  \sigma' \in \lgeo \rho , \sigma \lgeo $, there is $C_0 \in S'$ such that $\sigma \in C_0$. Recall from Section 
\ref{topol} notation 
$ C^+_j$ and $\gamma^+_j$, $j\in   J^+_{\sigma}$ and set 
$$A = \{ C^+_j \; , \; j\in   J^+_{\sigma} \; : \; d(\sigma , \gamma^+_j) <d(\sigma , \sigma')  \} $$
and 
$$ B=  \{ C^+_j \; , \; j\in   J^+_{\sigma} \; : \; d(\sigma , \gamma^+_j) =d(\sigma , \sigma') \quad 
{\rm and } \quad C^+_j \in S' \} \; .$$ 
Observe that $\cC_{\sigma} \cup A\cup B \in \Sigma_{\sigma}$ and that 
$$ L_{\sigma' } \cup \left(   \bigcup_{C' \in S'} C' \right)= 
 L_{\sigma } \cup \left(   \bigcup_{k\in I_{\sigma }} C_{\sigma} (k) \right) \cup \left( \bigcup_{C'' \in A\cup B } 
C'' \right) \; .$$
It implies that $t\in A_{\sigma} \backslash \{  \mu (L_{\sigma }) \}$ and it completes the proof of the lemma. \cqfd


\vspace{4mm}
Observe that by (\ref{debutprems}), we get $\phi (t) \leq \sigma $ for any 
$t\in [0,\mu (L_{\sigma } \backslash \{ \sigma \}) ]$. If $\sigma' \in C $ with $C\in \cC_{\sigma}^-$ , then 
we get 
\begin{equation}
\label{incluinclu}
 L_{\sigma' } \cup \left( \bigcup_{C' \in \cC_{\sigma'}} C' \right) \subset L_{\sigma } \backslash \{ \sigma \} \; .
\end{equation}
Thus, by (\ref{sesecond}) applied to $\sigma '$ we get 
$$ \sup \phi^{-1} (\{ \sigma' \}) \leq \mu (L_{\sigma' } ) + \sum_{C' \in \cC_{\sigma'}}\mu ( C')  $$
and (\ref{incluinclu}) implies 
$$ \sup \phi^{-1} (\{ \sigma' \}) \leq  \mu (L_{\sigma } \backslash \{ \sigma \}) \; .  $$
Consequently, 
$$  \phi^{-1} \left(L_{\sigma } \backslash \lgeo \rho ,\sigma \rgeo \right)= \bigcup_{C \in \cC_{\sigma}^-} \phi^{-1} 
\left(  C \right) \subset [ 0 , \mu (L_{\sigma } \backslash \{ \sigma \})] \; . $$
This, combined with (\ref{secondlinea}), implies that  
\begin{equation}
\label{inclureci}
[0, \mu (L_{\sigma })] \subset \phi^{-1} (L_{\sigma}) \subset [0, \mu (L_{\sigma })] \cup A_{\sigma} \; .
\end{equation}
We need the following lemma.  
\begin{lemma}
\label{lebenulle} For any $\sigma \in T$, $A_{\sigma}$ is a Lebesgue null set. 
\end{lemma}
\noindent
{\bf Proof}: First observe that 
$$ \mu (L_{\sigma}) =\inf A_{\sigma } \leq \sup A_{\sigma }= \mu (L_{\sigma }) + 
\sum_{C\in \cC_{\sigma} \cup \cC_{\sigma}^{+}} \mu (C)   =\mu (T)= M\; .$$
Set for any $C\in \cC_{\sigma} \cup \cC_{\sigma}^{+}$, $S(C)= \{ C' \in \cC_{\sigma} \cup \cC_{\sigma}^{+} \; : \; C' <C\}$. 
Clearly $S(C)$ and $S(C) \cup \{ C  \}$ are in $\Sigma_{\sigma}$. Thus, the real
numbers $a(C)$ and $b(C)$ given by 
$$a(C)=   \mu (L_{\sigma }) + \sum_{C'\in S(C)} \mu (C') \quad {\rm and} \quad   b(C)= \mu (C) + a(C) $$
are in $A_{\sigma}$. Observe now that $ \overline{A}_{\sigma} \cap (a(C) , b(C))=\emptyset $, 
where $\overline{A}_{\sigma}$ stands for the closure of the set $A_{\sigma}$. Thus 
\begin{equation}
\label{ineqqq}
\bigcup_{C\in \cC_{\sigma} \cup \cC_{\sigma}^{+} } (\, a(C) , b(C) \, ) \,
\subset \, [ \,  \mu (L_{\sigma }) ,\mu (T) \, ] \backslash 
\overline{A}_{\sigma} \; . 
\end{equation}
Note that 
$$ \lambda \left ( \bigcup_{C\in \cC_{\sigma} \cup \cC_{\sigma}^{+} } (\, a(C)
  , b(C)\, )\right) = 
\sum_{C\in \cC_{\sigma} \cup \cC_{\sigma}^{+} } \mu (C)= \mu (T) -\mu (L_{\sigma}) \; .$$
Thus, (\ref{ineqqq}) implies that $\lambda (\overline{A}_{\sigma})=0$, which completes the proof of the lemma. \cqfd

\vspace{4mm}

The previous lemma and (\ref{inclureci}) imply that for any $\sigma \in T$, 
$$ \lambda \left( \phi^{-1}(L_{\sigma}) \right)= \mu(L_{\sigma}) \; .$$ 
Thus, $\mu=\lambda \circ \phi^{-1}$ by Lemma \ref{pisyst}. Consequently 
$\mu = \mu_h \circ  \jmath_h ^{-1}$, by definition of 
$\jmath_h $.

\vspace{3mm}

\begin{lemma}
\label{lattermesu} For any $h\in \cH$, the set of times of first visit 
$F(h)$ and the set of times of latter visit $S(h)$ are Borel sets of the line.
\end{lemma} 
\noindent
{\bf Proof}: As already noticed, we have $p_h^{-1}({\rm Lf}(T_h)) \subset
F(h)$. Thus,  
\begin{equation}
\label{latterske}
 S(h) \; \subset \; p_h^{-1}({\rm Sk}(T_h)) \; .
\end{equation}
Let $t\in [0, \zeta(h)] $. We set 
$$ E_t= \left\{  s\in [t, \zeta (h)] \; : \; h(s)=\inf_{t\leq u\leq s} h(u)
  \quad {\rm and} \quad h(s) <h(t) \; \right\} \; .$$
Clearly, $E_t$ is a (possibly empty) Borel set of the line. Suppose that $E_t\neq \emptyset$. 
Let $s \in E_t$. Set 
$\ell =\sup \{ u\leq t \; : \; h(u) \leq h(s)\}$. Then 
$$ h(\ell )= h(s)=\inf_{\ell \leq u\leq s} h(u) \; ,$$
that is $p_h(\ell)= p_h(s)$. Now observe that since $h(t) >h(s)$, 
$p_h(t) \neq  p_h(s)$. Thus, $s\in S(h)$ since $\ell < t\leq s$.  
So we get 
$$ E_t \subset S(h) \; .$$

Let $t_n$, $n\geq 1$, be a sequence that is dense in $[0, \zeta (h)]$ and let 
$s\in S(h) $. Then, there exists $s'<s$ such that $p_h(s')=p_h(s)$. Since 
$p_h (s)$ is not a leaf, there exists $t_n \in (s',s)$ such that $h(t_n)>h(s)$. Then 
$$h(s)=\inf_{t_n\leq u\leq s} h(u) \, <  \,  h(t_n) \; , $$ 
which implies that $s\in E_{t_n}$. We thus have proved that 
\begin{equation}
\label{latteridinc}
 S(h)=\bigcup_{n\geq 1} E_{t_n} \; , 
\end{equation}
which implies the lemma. \cqfd 

%

\vspace{4mm}

   It remains to prove that $h$ satisfies (Min). Observe that for any $t\in [0, \zeta (h)]$, one has 
$$ E_t \subset \phi^{-1} (\lgeo \rho , p_h(t) \rgeo ) \cap ( \mu (L_{p_h (t)}), M]  \subset  A_{p_h(t)} \; .$$
Then (\ref{latteridinc}) implies 
$$ S(h) \subset \bigcup_{n\geq 1} A_{p_h(t_n)} \; , $$
which implies (Min) by Lemma \ref{lebenulle}.  This completes the
proof of the proposition. \cqfd 


\vspace{4mm}

The following proposition completes the proof of Theorem \ref{main}. 
\begin{proposition}
\label{uniq}
Let $h_1$, $h_2$ be two functions in $\cH$ that satisfy (Min) and such that the two structured trees 
$$(T_{h_1}, d_{h_1} , \rho_{h_1} , \leq_{h_1} ,\mu_{h_1}) \quad 
{\rm and } \quad (T_{h_2},
d_{h_2} , \rho_{h_2} , \leq_{h_2} ,\mu_{h_2})$$ 
are equivalent. Then, $h_1=h_2$. 
\end{proposition}
\noindent
{\bf Proof}: To simplify notation we assume that 
$$(T_{h_1}, d_{h_1} , \rho_{h_1} , \leq_{h_1} ,\mu_{h_1}) = (T_{h_2}, d_{h_2} , \rho_{h_2} , 
\leq_{h_2} ,\mu_{h_2})=(T,d, \rho , \leq ,\mu) .$$
First  observe that $\zeta (h_1) =\zeta (h_2)= \mu (T)=M$. 
Set for any $\sigma \in T$ and for $i\in \{ 1,2\}$ 
$$ \ell_i (\sigma)= \inf p_{h_i}^{-1} (\{ \sigma \}) \quad {\rm and } \quad 
r_i (\sigma)= \inf \{ t > \ell_i (\sigma) \; : \; p_{h_i} (t) \neq \sigma \}  $$
and recall that by Lemma \ref{leftlimi}, $ p_{h_i} (\ell_i (\sigma))=
\sigma$. By definition, if $\sigma' <\sigma$, then $\ell_i (\sigma')< \ell_i
(\sigma)$, for $i\in \{ 1,2\}$. Thus 
\begin{equation}
\label{morceaun}
p_{h_i} \left( [0, \ell_i (\sigma)) \right) =L_{\sigma} \backslash \{ \sigma \}
\; , \quad \sigma \in T \; , \; i\in \{ 1,2\}.
\end{equation}
Observe that 
\begin{equation}
\label{morceadeux}
p_{h_i} \left( [\ell_i (\sigma) , r_i (\sigma)] \right) =\{ \sigma \} \; , \quad \sigma \in T \; , \; i\in \{ 1,2\}.
\end{equation}
Let $t\in (r_i (\sigma) ,M]$ be such that $p_{h_i} (t) \in L_{\sigma} $;
observe that $t$ is necessarily a time of latter visit. Thus,  
\begin{equation}
\label{morceatrois}
p_{h_i}^{-1} \left(L_{\sigma} \right) \cap (r_i (\sigma) ,M] \subset S(h_i)
\; , \quad \sigma \in T \; , \; i\in \{ 1,2\}.
\end{equation}
Then (\ref{morceaun}), (\ref{morceadeux}) and (\ref{morceatrois}) imply for $i\in \{ 1,2\}$ 
$$ [0, \ell_i (\sigma)) \subset p_{h_i}^{-1} \left( L_{\sigma} \backslash \{
  \sigma \} \right) \subset [0, \ell_i (\sigma)) \cup  S(h_i) \quad {\rm and} \quad 
[0, r_i (\sigma)] \subset p_{h_i}^{-1}  \left( L_{\sigma}  \right)  \subset
[0, r_i (\sigma)] \cup  S(h_i) .$$
Consequently, 
$$ \mu ( L_{\sigma} \backslash \{
  \sigma \} )= \lambda \left( p_{h_i}^{-1} \left( L_{\sigma} \backslash \{
  \sigma \} \right) \right)=\ell_i (\sigma) \quad {\rm and} \quad  
\mu ( L_{\sigma} )= \lambda \left( p_{h_i}^{-1} \left( L_{\sigma}  \right)
\right)=r_i (\sigma) ,$$
since $h_1$ and $h_2$ satisfy (Min). (\ref{morceadeux}) then implies that $h_1$ and $h_2$ 
coincide on the set 
$$ F(h_1)=F(h_2)=\bigcup_{\sigma \in T} [\mu ( L_{\sigma} \backslash \{
  \sigma \} ),\mu ( L_{\sigma} )] $$
that is a set of full Lebesgue measure in $[0, M]$ and $h_1=h_2$ follows 
since $h_1$ and $h_2$ are left-continuous. \cqfd

\vspace{4mm}

We now prove Theorem \ref{timechange}. 

\vspace{4mm}

\noindent 
{\bf Proof of Theorem \ref{timechange}}: Recall that $\jmath_h$ (resp. $\jmath_{h'}$) stands for 
the isometry that maps the
structured tree $(T_{h}, d_{h} , \rho_{h} , \leq_{h}, \mu_{h})$ 
(resp. $(T_{h'}, d_{h'} , \rho_{h'} , \leq_{h'} , \mu_{h'})$) onto $(T, d, \rho
, \leq, \mu)$ (resp. onto $(T, d , \rho , \leq , \mu')$). Obviously, $\mu(T)=\zeta(h)$ and $\mu '(T)=\zeta(h')$. Let us
denote by $\phi$ (resp. $\phi'$) the exploration mapping from $[0, \mu(T)]$ 
(resp. $[0, \mu'(T)]$) onto $T$ associated with $\mu$ (resp. $\mu'$) as in 
Definition \ref{exploration}. Recall that Theorem \ref{main} implies that 
\begin{equation}
\label{specibij}
\jmath_h \circ p_h =\phi \quad {\rm and} \quad \jmath_{h'}\circ p_{h'} =\phi'
\; .
\end{equation}
Let us first prove the existence of the time-change. Recall that (\ref{preee}) implies 
\begin{equation}
\label{firstenfin}
F(h)= \bigcup_{\sigma \in T} \, [\,  \mu(L_{\sigma} 
\backslash \{  \sigma \} ) \, ,\,  \mu(L_{\sigma}  )\, ]
\quad  {\rm and} \quad 
F(h')= \bigcup_{\sigma \in T} \, [\,  \mu'(L_{\sigma} \backslash  
\{ \sigma \} ) \, ,\,  \mu'(L_{\sigma}  )\, ] \; .
\end{equation}
Since $h$ and $h'$ satisfy (Min), $F(h)$ and $F(h')$ are sets of full Lebesgue measure. Thus, they are 
dense in resp. $[0, \zeta (h)] $ and  $[0, \zeta (h')] $. Observe that it is possible to find a
non-negative application 
$\widetilde{\varphi}$ on $F (h)$ that is non-decreasing and such that for any
$\sigma \in T$

\begin{itemize}

\item{(a)}   $ \quad \widetilde{\varphi} (\mu (L_{\sigma} \backslash \{ \sigma \})) = 
\mu' (L_{\sigma} \backslash \{ \sigma \}) \; ;$

\item{(b)}   $\quad $ if $\mu(\{  \sigma \}) >0$, then $\widetilde{\varphi}$ is 
left-continuous on $(\mu (L_{\sigma} \backslash \{ \sigma \}), \mu
(L_{\sigma})]$ and $\widetilde{\varphi} (\mu (L_{\sigma} )) \leq
\mu' (L_{\sigma}) \; .$
\end{itemize}

\begin{remark}
\label{nonunici}
If $\mu (\{ \sigma \})\mu' (\{ \sigma \}) >0$, then observe that we can find 
infinitely many $\widetilde{\varphi}$ satisfying (a) and (b). \cq 
\end{remark}

\vspace{3mm}

We define $\varphi $ by 
$$ \varphi (t)= 
\sup \{ \widetilde{\varphi} (s) \; , \; s\in F (h) \quad {\rm
  and} \quad s\leq t \} \; , \quad t\in [0, \mu (T)]  \; . $$
Observe that $\varphi$ and $\widetilde{\varphi}$ coincide on $F (h)$. Consequently 
$\varphi (0)=\widetilde{\varphi}(0)= 0 $.

   Let $t\in (0, \mu (T) ]$. We claim that there exists an increasing sequence 
$s_n \in [0, \mu (T)] $, $n\geq 1$, converging to $t$ and such that 
\begin{equation}
\label{suisuite}
\lim_{n\rightarrow \infty} \varphi (s_n) = \varphi (t) \; .
\end{equation}
Since $\varphi$ is non-decreasing, the previous claim easily implies that $\varphi$
is left-continuous at $t$. Let us prove (\ref{suisuite}): the result
is clear if $t\notin F(h)$; the only non-trivial case to consider then, is 
when $t=\mu (L_{\sigma} \backslash \{ \sigma \}) $, $\sigma \in T\backslash \{
\rho\}$. It is always possible to find an increasing sequence $\sigma_n \in \lgeo \rho ,
\sigma \lgeo $, $n\geq 1$, that converges to $\sigma $ and such that $\mu (\{ \sigma_n
\}) =0$ , $n\geq 1$. Thus, Lemma \ref{increa} implies 
\begin{equation}
\label{encore}
\bigcup_{n\geq 1} L_{\sigma_n} \backslash \{ \sigma_n\} = \bigcup_{n\geq 1} L_{\sigma_n}  =L_{\sigma}
\backslash \{ \sigma\} \; .
\end{equation}
Set for any $n\geq 1$, $s_n= \mu (L_{\sigma_n} \backslash \{ \sigma_n\}
)$. Clearly, $\varphi (s_n)=\widetilde{\varphi} (s_n)=\mu '(L_{\sigma_n}
\backslash \{ \sigma_n\})$ and (\ref{encore}) implies 
$$  \lim_{n\rightarrow \infty} \varphi (s_n) = 
\lim_{n\rightarrow \infty } 
\mu ' (L_{\sigma_n} \backslash \{ \sigma_n\}) =\mu ' (L_{\sigma } \backslash
\{ \sigma \}) = \varphi (t) $$
which implies (\ref{suisuite}).

\vspace{4mm}

Thus, we have constructed a non-decreasing, left-continuous mapping $\varphi :[0, \zeta (h)]\rightarrow [0, \infty )$ 
that coincides with $\widetilde{\varphi}$ on $F(h)$ and such that $\varphi(0)=0$. Moreover, (a) and (b) imply that 
$\phi =\phi' \circ \varphi $ on $F(h)$ by (\ref{firstenfin}). It easily 
implies $h=h'\circ \varphi $ on $[0, \zeta (h)]$ since $h$ and 
$h'\circ \varphi $ are left-continuous and since $F(h)$ is dense in 
$[0, \zeta (h)]$.

\vspace{4mm}

Let us prove the uniqueness result and the other points of Theorem
\ref{timechange}. To that end, we need the following
proposition.
\begin{proposition}
\label{firstcontr}
Let $(T, d, \rho)$ be a rooted compact real tree and let $\leq $ be a linear order satisfying (Or1) and 
(Or2). Let $\mu$ and $\mu'$ be two finite Borel measures on $T$ that both satisfy (Mes). Denote by $h$ and $h'$ the height 
functions associated with resp. $(T, d, \rho , \leq, \mu)$ and $(T, d, \rho , \leq, \mu' )$ by 
Theorem \ref{main} ($h$ and $h'$ then satisfy (Min)). Assume that $\varphi:[0, \zeta (h)] \rightarrow [0, \infty
)$ is non-decreasing, left-continuous and such that 
$$\varphi (0)=0 \quad  {\rm and} \quad h=h'\circ \varphi \; . $$ 
Then for any $\sigma \in T$, we get  
$$ \varphi (\mu (L_{\sigma} \backslash \{ \sigma \})) = 
\mu' (L_{\sigma} \backslash \{ \sigma \}) \quad {\rm and} \quad 
 \varphi (\mu (L_{\sigma} )) \leq
\mu' (L_{\sigma}) \; .$$
\end{proposition}

\vspace{4mm}

\noindent
{\bf Proof of Proposition \ref{firstcontr}:} To simplify notation, we set 
$f=\jmath_{h}^{-1} \circ \jmath_{h'}$; then $f$ maps the rooted ordered
compact real tree $(T_{h'}, d_{h'}, \rho_{h'}, \leq_{h'} )$ onto 
$(T_{h}, d_{h}, \rho_{h}, \leq_{h} )$. We first want to prove
\begin{equation} 
\label{identiti}
 \jmath_{h} \circ p_h = \jmath_{h'} \circ p_{h'} \circ \varphi \; .
\end{equation}
First, let us fix $s_1,s_2 \in [0, \zeta (h)]$ and let us set 
$\sigma_1=p_h (s_1)$ and $\sigma_2=p_h (s_2)$. For any $\sigma \in T_h $,
denote by $t_{\sigma} \in [0, \zeta (h)]$ a time such 
that $p_h( t_{\sigma})=\sigma $; assume that $t_{\sigma_1}=s_1$ and 
$t_{\sigma_2}=s_2$. We then define $G$ from $T_h$ to $T_h$ by 
\begin{equation} 
\label{gedef}
G(\sigma ) = f(p_{h'} (\varphi (t_{\sigma})) ) \; , \quad \sigma \in T_h \; .
\end{equation}
We first get for any $\sigma , \sigma' \in T_h$
$$ d_h (G(\sigma ), G(\sigma '))= d_{h'} ( p_{h'} (\varphi (t_{\sigma})) , p_{h'} (\varphi (t_{\sigma'})) ) 
=  d_{h'} ( \varphi (t_{\sigma}) , \varphi (t_{\sigma'}) )$$
for  $f$ is an isometry. Then, note that 
$$  d_{h'} ( \varphi (t_{\sigma}) , \varphi (t_{\sigma'}) ) \geq d_h (t_{\sigma} , t_{\sigma'} ) =d_h(\sigma , \sigma') $$
since $h=h' \circ \varphi $. Thus, for any $\sigma, \sigma' \in T_h$
$$  d_h (G(\sigma ), G(\sigma ')) \geq d_h(\sigma , \sigma') \; . $$
Since the metric space $(T_h ,d_h)$ is compact, standard arguments imply that $G$ is actually a bijective isometry 
(see Theorem 1.6.15 (2) in \cite{BuBu}). It first implies that $p_{h'} \circ \varphi$ is surjective: 
\begin{equation} 
\label{surjj}
 p_{h'} \circ \varphi \, \left( [0, \zeta (h)] \right) = T_{h'}  \; .
\end{equation}
It also implies that 
$d_h (G(\sigma_1 ), G(\sigma_2)) = d_h(\sigma_1 , \sigma_2) $. So, 
we have proved 
\begin{equation} 
\label{egagali}
 d_h(s_1 , s_2)= d_{h'}(\varphi (s_1) , \varphi(s_2) ) \; , \quad s_1,s_2 \in
 [0, \zeta (h)] \; .
\end{equation}

\vspace{4mm}

Let us prove that $G$ preserves $\leq_h $. To that end, we need to prove the following lemma. 

\begin{lemma}
\label{interconst}
Let $h,h'\in \cH$ and $\varphi : [0, \zeta (h)] \rightarrow [0, \infty) $ be as in Proposition \ref{firstcontr}. Let 
$t\in [0, \zeta (h))$ be such that $\varphi (t) < \varphi (t+ )$. Then 
$$ h'(u)=h' (\varphi (t))=h(t) \; , \quad u \in [\varphi (t) ,  \varphi (t+)] \; .$$
\end{lemma}

\noindent
{\bf Proof of Lemma \ref{interconst}:} We introduce 
$$s_0= \inf \{ s \in  [\varphi (t) ,  \zeta (h')] \; : \;  h'(s) \neq h' (\varphi (t)) \}  \; ,$$
with the convention that $ \inf \emptyset =\infty$. Suppose that $s_0 < \varphi (t+)$. Then,  
$F (h') \cap (s_0 ,  \varphi (t+))$ is non-empty for 
$F (h') $ is dense in $[0, \zeta (h')]$ by (Min). Consequently, we can find 
$s \in F (h') \cap (s_0 ,  \varphi (t+))$ such that $h'(s) \neq h' (\varphi (t))$. 
There exists $u\in [0, \zeta (h)]$ such that $p_{h'} (\varphi (u ) )= p_{h'} (s)$ 
since $p_{h'} \circ \varphi$ is surjective. Set $\sigma= p_{h'} (s)$. If 
$\varphi (u) <s $, then $ \varphi (u) \leq   \varphi (t) $ and since $s\in F (h')$, 
it implies 
$$[\varphi (u) , s] \subset [\mu' (L_{\sigma }  \backslash \{ \sigma  \}) , \mu' (L_{\sigma }) ]$$
by (\ref{firstenfin}). Consequently, we get $ \sigma= p_{h'} (r)$, $r\in
[\varphi (u) , s]$. But $s_0 \in [\varphi (u) , s]$, which rises a
contradiction. 
Thus, $\varphi (u) >s $. It implies that $u>t$ and $\varphi (u)>\varphi (t+)$. Since $d_{h'} (\varphi (u) ,s)=0$, we get 
\begin{equation} 
\label{majhei}
 h' (\varphi (u)) =h'(s) =\inf_{r\in [s,\varphi (u) ]} h'( r )\leq h' ( \varphi(t+)+)=h(t+) \; . 
\end{equation}
Now deduce from (\ref{egagali}) that for any $\epsilon >0$
$$ d_h(t, t +\epsilon )= d_{h'} (\varphi (t), \varphi(t+\epsilon ) ) \; . $$
Observe that 
$$ \lim_{\epsilon \rightarrow 0} d_{h}(t, t +\epsilon )= h(t)-h(t+) $$
and that 
\begin{eqnarray*}
\lim_{\epsilon \rightarrow 0} d_{h'}(\varphi (t), \varphi(t+\epsilon ) ) & = &
h'(\varphi(t))+ h' ( \varphi(t+)+) -2\,\left(  h' ( \varphi(t+)+)   \wedge   
 \inf_{r\in [\varphi (t), \varphi (t+) ]} h'( r ) \right) \\
& = & h(t)+ h ( t+)
-2\, \left( h ( t+)  \wedge    \inf_{r\in [\varphi (t), \varphi (t+) ]}h'( r ) \right) 
\; . 
\end{eqnarray*}
Consequently, 
\begin{equation} 
\label{minhei}
h(t+) \leq  \inf_{r\in [\varphi (t), \varphi (t+) ]}h'( r ) \; . 
\end{equation}
Then (\ref{majhei}) and (\ref{minhei}) both imply that $h' (s)=h(t+)$ for 
$\inf_{r\in [s, \varphi (u) ]} h'( r ) \leq h' (s)$. Since we have supposed that 
$h'(s) \neq h' (\varphi (t))$, we get 
$$ h'(s)=h(t+) < h(t)=h'(\varphi (t)) \; . $$ 
Now set 
$$ u_0= \sup \{ u \in  [0 ,  \varphi (t)) \; : \;  h'(u) \leq h'(s)=h (t+) \}  \; .$$
Clearly, 
$$ h'(u_0)= \inf_{r\in [u_0,s ]} h'( r )=h'(s)= h (t+) \; .$$
Thus, $d_{h'}(u_0, s)=0$. Since $s\in F(h')$, we get 
$[u_0, s] \subset [\mu' (L_{\sigma }  \backslash \{ \sigma  \}) , \mu' (L_{\sigma }) ]$. 
Thus, for any $r\in [u_0, s]$, we get $p_{h'} (r)=\sigma $. But $s_0 \in [u_0, s]$ and 
$$h'(s_0)=h'(\varphi (t))>d(\rho, \sigma )=h'(s) \; .$$
Thus, $p_{h'}(s_0)\neq \sigma $, which rises a contradiction. It then proves that  $s_0 \geq  \varphi (t+)$, 
which implies the lemma. \cqfd  

\vspace{4mm}

Let us complete the proof of the proposition. 
Let $t\in (0, \zeta (h)]$ such that for any $s\in [0,t)$, $p_h (s) \neq p_h(t)$. 
Then, the previous lemma implies that 
\begin{equation} 
\label{entranan}
p_{h'}(u)\neq p_{h'}(\varphi (t)) \; , \quad u \in [0,\varphi (t)) \; . 
\end{equation}
(Suppose indeed that there exists $u \in [0,\varphi (t)) $ such that
$p_{h'}(u)=  p_{h'}(\varphi (t)) $ which is equivalent to $d_{h'} (u,\varphi
(t))$; there 
exists $s\in [0, t)$ such
that $\varphi (s) \leq u \leq \varphi (s+) $; the previous lemma implies that 
$d_{h'}(\varphi (s),u)=0$; thus $d_{h'}(\varphi (s),\varphi (t))=0$;
(\ref{egagali}) implies $d_{h}( s,t)=0$, which rises a contradiction.)

\vspace{4mm}

Consequently, we get for any $s_1,s_2 \in [0, \zeta (h)]$, 
\begin{equation} 
\label{ordcons}
p_h (s_1) \leq_h \, p_h (s_2) \quad \Longrightarrow \quad  p_{h'} (\varphi( s_1) ) \leq_{h'} \, p_{h'} (\varphi( s_2) ) \; .
\end{equation}
Thus, by definition of $G$
\begin{equation}
\label{croiincrea}
\sigma_1 \leq_h \sigma_2 \; \Longrightarrow  \; G(\sigma_1 )\leq_h G(\sigma_2) \; .
\end{equation}

\vspace{3mm}

Let us prove that (\ref{croiincrea}) implies (\ref{identiti}): 
Let us suppose that there exits $\sigma \in T_h $ such that $G(\sigma )\neq \sigma $. Then, eiher 
$G(\sigma) <_h \sigma $, either $\sigma  <_h G(\sigma)$. Define recursively 
the sequence $G^n( \sigma )$ , 
$n\geq 1$, by $G^{n+1}( \sigma )=G(G^n( \sigma ))$ and $G^1( \sigma )=G( \sigma )$. In both cases 
the sequence $G^n( \sigma )$ , 
$n\geq 1$, is $(\leq_h )$-monotone and thus convergent, by Proposition \ref{monot}. Now, observe that 
$d_h(G(\sigma) , \sigma)=d_h(G^{n+1}(\sigma) , G^n (\sigma))$, which contradicts $G(\sigma )\neq \sigma $.

\vspace{3mm}

So, we have proved that $G(\sigma)=f(p_{h'}(\varphi (t_{\sigma}) ))=\sigma $. By (\ref{egagali}), this identity 
does not depend on the choice of $t_{\sigma}$ in $p_h^{-1} (\{ \sigma\})$. Consequently 
$$ f(p_{h'}(\varphi (t ) ))= p_h(t) \; , \quad t\in [0, \zeta (h)] \; , $$ 
which implies (\ref{identiti}). 

\vspace{4mm} 

Let us now complete the proof of Proposition \ref{firstcontr}: recall from (\ref{debutprems}) that for any $\sigma \in T$
$$ \mu (L_{\sigma} \backslash \{ \sigma \} )= \inf \{ t\geq 0 \; : \; \jmath_h (p_h (t))=\sigma  \} \; .$$
Then, (\ref{entranan}) and (\ref{identiti}) imply 
\begin{equation} 
\label{boutun}
\varphi ( \mu (L_{\sigma} \backslash \{ \sigma \} ) )= \inf \{ t\geq 0 \; : \; \jmath_{h'} (p_{h'} (t))=\sigma  \} = 
\mu' (L_{\sigma} \backslash \{ \sigma \} ). 
\end{equation}
This proves the first point of Proposition \ref{firstcontr}.

  Let us prove the
last one.  Fix $t\in [\mu (L_{\sigma} \backslash \{ \sigma \} ), \mu
(L_{\sigma})] $. (\ref{identiti}) implies that 
$$\sigma= \jmath_{h'}( p_{h'}(\varphi (t ) ))=\jmath_{h}( p_{h}(t )) \; , \quad 
t\in [\mu (L_{\sigma} \backslash \{ \sigma \} ), \mu (L_{\sigma})] \; . $$
Suppose that 
$\varphi(t) \notin F (h')$. Thus, $\varphi(t) > \mu' (L_{\sigma} )$ by (\ref{firstenfin}). Since 
$ p_{h'}(\mu' (L_{\sigma}))= p_{h'}(\varphi(t)) =\sigma $, there exists $s\in (\mu' (L_{\sigma}) , \varphi(t) ) $ such 
that $ p_{h'}(s) \neq \sigma $. Set $ p_{h'}(s) =\sigma ' $. It is easy 
to check that 
$$  p_{h'}^{-1} (\{ \sigma ' \})  \subset (\mu' (L_{\sigma}) , \varphi(t) ) \; . $$
Since $ p_{h'}\circ \varphi $ is surjective, we can find $u<t$ such that 
\begin{equation} 
\label{cocontr}
\varphi(u) \in ( \mu' (L_{\sigma} ),\varphi(t) ) \quad {\rm and}  \quad  
\sigma' =  \jmath_{h'}( p_{h'}(\varphi (u ) )) \; .
\end{equation}
It implies  
$$ \varphi  (\mu (L_{\sigma} \backslash \{ \sigma \} )) = \mu' (L_{\sigma} \backslash \{ \sigma \} ) \leq 
\mu' (L_{\sigma}) < \varphi (u) \leq \varphi (t) \; . $$
Thus, $u\in [\mu (L_{\sigma} \backslash \{ \sigma \} ) , t]$. (\ref{preee}) then implies that 
$$ \sigma =\jmath_{h}( p_{h}(u )) = \jmath_{h'}( p_{h'}(\varphi (u ) )) \; , $$
which contradicts (\ref{cocontr}). Thus, it proves that $\varphi(t) \in F (h')$ and (\ref{firstenfin}) implies that 
$\varphi(t) \leq \mu'(L_{\sigma})$, which completes the proof of the proposition. \cqfd

\vspace{4mm}

Let us complete the proof of Theorem \ref{timechange}: Observe that $(iv)$ is Lemma \ref{interconst}. 
Let us prove $(ii)$: if $\mu'$ has no atom, then 
Proposition \ref{firstcontr} implies 
$$ \varphi (F(h)) =F(h') = \{ \mu' (L_{\sigma} ) \; : \; \sigma \in T \} \; , $$
which implies that $\varphi$ is continuous since $F(h)$ and $F(h')$ are dense. 

\vspace{4mm}

   Let us prove $(iii)$: assume that $\mu$ has no atom.  
let $t_1 <t_2$ be in $(0, \mu(T))$; there exist 
$\sigma_1 <\sigma_2$ in $T$ such that 
$$ t_1 < \mu ( L_{\sigma_1}) <\mu ( L_{\sigma_2}) <t_2 \; ,$$
for  $F(h)= \{ \mu (L_{\sigma }); \sigma \in T \}$ is dense 
in $[0, \mu(T)]$; it implies by (Inc) that 
$$ \varphi (t_1) \leq \varphi (\mu ( L_{\sigma_1}))=
\mu' ( L_{\sigma_1}\backslash \{ \sigma_1\}) < \mu' ( L_{\sigma_2}\backslash \{ \sigma_2\}) = 
\varphi (\mu ( L_{\sigma_2})) \leq  \varphi (t_2)$$
and thus $\varphi (t_1)< \varphi (t_2)$, which completes the proof of  $(iii)$.

\vspace{4mm}

Finally, let us prove $(i)$. Assume that $\mu$ and $\mu'$ do not share 
any atom. Thus, Proposition \ref{firstcontr} and 
(\ref{firstenfin}) imply that $\varphi$ is uniquely
determined on $\{0\}\cup F (h)$. Consequently, $\varphi$ is uniquely
determined on $[0, \zeta(h)]$ for $F (h)$ is dense in $[0, \zeta(h)]$ by
(Min) and for $\varphi$ is left-continuous. This proves one implication of $(i)$;
the converse of $(i)$ is a consequence of Remark \ref{nonunici}. \cqfd

\section{Properties of height functions.}
\label{propriri}

   In this section we give some simple properties concerning the regularity of height functions 
in terms of properties of the corresponding trees. We also make the connection with 
an earlier probabilistic approach by Aldous.

Let $(T, d)$ be a real tree. Recall the definition of the length measure $\ell_T$ on $T$ from the Introduction section.
Observe that  $\ell_T$ only relies on the metric 
structure. That the tree has finite length should be ``read'' from any  height function coding the tree. 
More precisely, let $h\in \cH_M$. We set for any 
$0\leq a\leq b\leq M$ 
$$ v(h,[a,b])= \sup \, \sum_{1\leq i\leq n}  | h(t_i)-h(t_{i-1})| \; ,$$
where the supremum is taken over all subdivisions $t_0=a < t_1 < \ldots < t_n=b$. The 
(possibly infinite) quantity $ v(h,[a,b])$ is then the total variation of $h$ over $[a,b]$. Let $r\in [0,M]$
and let $t_0=0 < t_1 < \ldots < t_n=r$. We denote by ${\rm Span}_h (t_1, \ldots , t_n)$ the subtree of 
$T_h$ spanned by the vertices $p_h(t_1), \ldots , p_h(t_n)$ and the root $\rho_h$: 
$${\rm Span}_h (t_1, \ldots , t_n) = \bigcup_{1\leq i\leq n} \lgeo \rho_h , p_h(t_i) \rgeo  \; . $$
Fix $\epsilon >0$. For any $1\leq i\leq n-1$ , we can find $s_i (\epsilon ) \in [t_i, t_{i+1}] $ such that 
\begin{equation}
\label{subdi}
 h(s_i (\epsilon )) \leq \frac{\epsilon}{n} + \inf_{s\in [t_i, t_{i+1}]} h(s) \; . 
\end{equation}
Now, think of the rooted ordered subtree  ${\rm Span}_h (t_1, \ldots , t_n)$ as a planar tree, namely a tree embedded in the 
clockwise oriented half-plane; imagine a particle that continuously moves on it at unit speed, that starts at the root 
$\rho_h$ and that backtraks as less as possible. The total amount of time needed by the particle to cover the tree and 
to go back to the root is twice the total length of ${\rm Span}_h (t_1, \ldots , t_n)$. More precisely 
the function recording the distance 
of the particle from the root is the piecewise linear continuous function with slope $+1$ or $-1$ that goes through the 
values 
$$ 0 \, , \,  h(t_1) \, , \, \inf_{s\in [t_1, t_{2}]} h(s) \, , h(t_2) \, , \,  
\ldots , \, \inf_{s\in [t_{n-1}, t_{n}]} h(s) \, , \, 
h(t_n)  \, , \,  0 \; . $$
If we look at the particle until it visits for the last (and perhaps also the first) time $p_h(t_n)$, then all the point of 
${\rm Span}_h (t_1, \ldots , t_n)$ have been visited twice or more except the points of 
$\rgeo \rho_h , p_h(t_n) \lgeo \backslash {\rm Br}(T)$,
the leaves of 
${\rm Span}_h (t_1, \ldots , t_n)$ and possibly the root. Deduce from the previous observations that 
\begin{eqnarray*}
 \sum_{1\leq i\leq n}  | h(t_i)-h(t_{i-1})| & \leq  & d_h(\rho_h , p_h(t_1)) + d_h (p_h(t_1) , p_h(t_2)) + \ldots + 
d_h (p_h(t_{n-1}) , p_h(t_n)) \\
&  \leq  & 2 \ell_T ({\rm Span}_h (t_1, \ldots , t_n)) - h(r)
\end{eqnarray*}
since $d_h(\rho_h , p_h( t_n))=h(r)$. Now deduce from (\ref{subdi}) that 
\begin{eqnarray*}
 h(t_1)+|h(t_1)-h(s_1 (\epsilon))| +  |h(t_2)-h(s_1 (\epsilon))| + 
  \ldots  + |h(t_n)- h(s_{n-1} (\epsilon))|   \quad \geq \\ 
 d_h(\rho_h , p_h(t_1)) + d_h (p_h(t_1) , p_h(t_2)) + \ldots + 
d_h (p_h(t_{n-1}) , p_h(t_n)) \; -\; \epsilon \; . 
\end{eqnarray*}
Consequently 
$$v(h,[0,r]) = 2 \ell_T (p_h ([0,r]) ) - h(r) $$
for 
$$  \sup \;  \ell_T ({\rm Span}_h (t_1, \ldots , t_n)) =  \ell_T (p_h ([0,r]) ) \; ,$$
where the supremum is taken over all the subdivisions $t_0=0 < t_1 < \ldots < t_n=r$. 
This implies the following proposition. 
\begin{proposition}
\label{vartot}
Let $(T, d)$ be a compact real tree. $\ell_T (T)$ is a finite quantity iff there exists 
a height function $h\in \cH$ with bounded variation such that $(T,d)$ and  $(T_h ,d_h)$ are isometric. 
\end{proposition}
\begin{remark}
\label{linhaus}
It is easy to check that if the length $\ell_T (T)$ is finite then Hausdorff
and packing dimensions agree and are equal to 
$1$. \cq 
\end{remark}

\vspace{4mm}

   We now discuss continuity properties of height processes. Let us first prove the 
following lemma. 
\begin{lemma}
\label{continuous} Let $h \in \cH$. We can always find a continuous $c\in \cH$
such that $(T_h,d_h ,\rho_h , \leq_h)$ and $(T_c,d_c, \rho_c , \leq_c)$ are isometric. 
\end{lemma}
\noindent
{\bf Proof}: Let us first mention that $c$ is in general not unique and that it may never satisfy (Min) 
(see Comment \ref{comaa}). Here we provide one possible function $c$ by interpolating the jumps of $h$ in an order-preserving way. 
Denote by $t_n$ , $n\geq 1$, a sequence of $[0, \zeta (h)]$
containing all the jump-times of $h$; set 
$$ \psi (t)=t+\sum_{n\geq 1} 2^{-n}\un_{[0, t]} (t_n) \quad 
{\rm and} \quad \Lambda (s)= \inf \{t\in [0, M] \; : \; \psi (t)>s  \}. $$
Clearly $\psi$ is increasing and right-continuous on $[0, M]$. Thus,
$\Lambda$ is well defined and continuous on $[0, M+1)$. Moreover,  
$\lim_{s\rightarrow M+1} \Lambda (s) =M$. Set $a= \psi (\Lambda (s)-)$ and 
$b=\psi (\Lambda (s)) $. If $a < b $, then $a\leq s\leq b$ and for any $u\in
(a, b)$, we get $\Lambda (u)=\Lambda(s)$. In that case define 
$$\theta (u)= \left( h(\Lambda (s))- h(\Lambda (s)+) \right) \, 
\frac{u-a}{b-a}  \; . $$
If otherwise $a=b$, then set $\theta (s)=0$. Then, define 
$c(s)=h(\Lambda (s))-\theta (s)$, $s\in [0, M+1]$. Check that $c$ is continuous 
and that $(T_h,d_h ,\rho_h , \leq_h)$ and $(T_c , d_c ,\rho_c , \leq_c)$ are isometric. \cqfd

\vspace{4mm}

   Let $(T, d, \rho , \leq , \mu)$ be a structured tree such that $\leq $ satisfies (Or1) and (Or2) and such that $\mu $ 
satisfies (Mes). Recall that $\phi : [0, \mu(T) ] \rightarrow T$ stands for the exploration mapping associated with 
$\mu$ defined in Definition \ref{exploration}. Fix $t\in [0, \mu(T) ]$.  
It is easy to check that $\phi (t)\neq \phi (t+)$ iff no subtrees are
grafted on the ``right side'' of the branch $\rgeo \phi(t+) , \phi (t)\rgeo $. Namely, $\phi (t)\neq \phi (t+)$ iff  
$$ \left\{ \, \sigma \in T \; : \; \phi (t)< \sigma \quad {\rm and } \quad \sigma \wedge \phi (t) 
\in  \rgeo \phi (t+) , \phi (t) \lgeo \, \right\} \; = \; \emptyset . $$
Thus, the height process $h\in \cH$ associated with the structured tree $(T, d, \rho , \leq , \mu)$ by Theorem \ref{main} is 
continuous iff for any $\sigma_1 \in T$ and for any $\sigma_2 \in \rgeo \rho , \sigma_1 \lgeo $, 
\begin{equation}
\label{critun}
\left\{  \, \sigma \in T \; : \;  \sigma_1 <\sigma  \quad {\rm and }
\quad \sigma \wedge \sigma_1 \in \rgeo \sigma_2 , \sigma_1 
\lgeo \, \right\}
 \; \neq  \; \emptyset .
\end{equation}
It implies that the leaves of $T$ are dense: 
\begin{equation}
\label{critde}
 \overline{{\rm Lf} (T)}=T \; .
\end{equation}
(Indeed let $\gamma \in \rgeo \sigma_2 , \sigma_1 
\lgeo $ and fix $\epsilon >0$; (\ref{critun}) implies that ${\rm Br }(T)$ is dense in 
$\rgeo \sigma_2 , \sigma_1  \lgeo $; since $(T,d)$ is compact, there are only finitely many 
connected components of $T\backslash \rgeo \sigma_2 , \sigma_1 \lgeo $ with a diameter larger than 
$\epsilon $; Thus the set of points in $\rgeo \sigma_2 , \sigma_1  \lgeo $ 
on which are grafted the connected components of $T\backslash \rgeo \sigma_2 , \sigma_1 \lgeo $ with diameter 
$\leq \epsilon $ is dense in $\rgeo \sigma_2 , \sigma_1  \lgeo $; consequently we can find a leaf 
$\sigma \in {\rm Lf} (T)$ in such a component such that $d(\sigma, \gamma) \leq 2 \epsilon $; it implies 
that the leaves are dense in the skeleton of $T$, which proves (\ref{critde}).)

Conversely, we prove the following proposition.
\begin{proposition}
\label{shuffleconti}
Let $(T, d, \rho )$ be a compact rooted real tree such that $ \overline{{\rm
    Lf} (T)}=T $. Then a.s. for any finite Borel measure $\mu $ whose topological support is
$T$, the height process $h_{{\rm Sh}}$
associated with the structured tree $(T, d, \rho , \leq_{{\rm Sh}}
,\mu)$  by Theorem \ref{main} is continuous. 
\end{proposition}
\noindent
{\bf Proof}: Clearly,  $\overline{{\rm Lf} (T)}=T $ implies 
$$\forall \sigma_1 \in T \; , \; \forall \sigma_2 \in \rgeo \rho ,
\sigma_1 \lgeo \; , \; \exists \sigma 
\in T \backslash \lgeo \sigma_2,
\sigma_1 \rgeo
\; : \;  \sigma  \wedge \sigma_1 \in \rgeo \sigma_2 , \sigma_1 \lgeo  \; . $$
Arguments similar to those used in the proof of 
Proposition \ref{Shufflemes} imply that  a.s. the ordered tree $(T, d, \rho ,
\leq_{{\rm Sh}})$ satisfies  
(\ref{critun}). The details are left to the reader. \cqfd

\vspace{4mm}

\begin{remark}
\label{lamesu}
Although (\ref{critde}) does not depend on any measure on $T$, note 
that if there exists a measure $\mu $ on $(T, d, \rho)$ such that 
\begin{equation}
\label{crittrois}
{\rm supp} \, \mu =T \quad {\rm and} \quad \mu ( {\rm Sk}(T))=0 \; ,  
\end{equation}
then $(T, d, \rho)$ satisfies (\ref{critde}). Observe that if in addition 
$\mu $ is non-atomic, then any height function 
$h\in \cH$ coding $(T, d, \rho , \leq , \mu)$ fulfils 
(Min) and is therefore unique. \cq 
\end{remark}

\vspace{4mm}

Conversely we have the following proposition.

\begin{proposition}
\label{existconti}
Let $(T, d, \rho )$ be a compact rooted real tree such that 
$ \overline{{\rm Lf} (T)}=T $. Then there exists a probability measure $\mu$ on $T$ 
that satisfies (CT1), (CT2) and (CT3). 
\end{proposition}
\noindent
{\bf Proof}: We construct such a probability measure thanks to a specific 
splitting of $T$ that we first explain: Let 
$Q_{\varnothing}=(q^{(\varnothing )}_n ; n\geq 1)$ be a dense sequence of distinct
leaves of $T$. Let $C^o_k$ , $k\geq 1$ be the connected components of
$T\backslash \lgeo \rho , q^{(\varnothing )}_1 \rgeo $ listed in such a way that
for any $1\leq k<l$,  
$$ \min \; \{ n\geq 1 \; : \; q^{(\varnothing )}_n \in C^o_k \} \; < \; 
\min \; \{ n\geq 1 \; : \; q^{(\varnothing )}_n \in C^o_l \} \; .$$
Fix $k\geq 1$. Denote by $C_k$ the closure of $C^o_k$ and denote by 
$\sigma_k$ the vertex of $\lgeo \rho ,
q^{(\varnothing )}_1 \rgeo $ such that $C_k= C^o_k \cup \{\sigma_k \}$. We also
define $Q_{k}=(q^{(k)}_i ; i\geq 1)$ by $ q^{(k)}_i= q^{(\varnothing )}_{n(i)}$,
where $n(i)$, $i\geq 1$, is the increasing sequence of indices $n\geq 1$ such
that $q^{(\varnothing )}_n \in C^o_k $. Then, we have defined 
$$ {\rm Split} \left( \, ((T,d,\rho ); Q_{\varnothing}) \, \right) =  \left( \,
  ((C_k ,d,
  \sigma_k); Q_{k} ) \; ; \; k\geq 1\right) \; .$$
We recursively define $\left( (C_u ,d, \sigma_u); Q_{u} \right) $ for any word
$u\in \U$ in the following way: 
$$ {\rm Split}  \left( \, ((C_u ,d, \sigma_u); Q_{u} ) \, \right) =  \left( \,
  ((C_{(u,k)} ,d,
  \sigma_{(u,k)}); Q_{(u,k)}) \; ; \; k\geq 1\right) \; ,$$
where $(u,k)$ stands for the concatenation of the word $u$ with the single
letter word $k$. Observe that for any $u=(v,w)$ with $v,w\in \U$ we get 
\begin{equation}
\label{nested}
C_u \subset C_v
\end{equation}
and 
\begin{equation}
\label{racinul}
 \sigma_u \in C^o_v \quad {\rm if} \quad w\neq \varnothing \; ,
\end{equation}
where $ C^o_v $ stands for the interior of the compact set $ C_v $. Note
that for any $n\geq 1$
\begin{equation}
\label{excu1}
\bigcup_{\substack{u\in \U \\ |u|=n}}C^o_u =T\backslash
\bigcup_{\substack{u\in \U \\ |u|\leq n}} \lgeo \rho , q^{(u)}_1 \rgeo 
\end{equation}
and by (\ref{racinul}) we also get 
\begin{equation}
\label{excu2}
\bigcup_{\substack{u\in \U \\ |u|=n}}C_u \subset \bigcup_{\substack{v\in \U \\
    |v|=n-1}} C^o_v \; .
\end{equation}
Thus 
\begin{equation}
\label{excu3}
\bigcap_{n\geq 1} \bigcup_{\substack{u\in \U \\ |u|=n}}C_u = T\backslash
\bigcup_{u\in \U } \lgeo \rho , q^{(u)}_1 \rgeo \; .
\end{equation}
Now observe that 
\begin{equation}
\label{excu4}
\{ q^{(u)}_1 \; ; \; u\in \U\} =\{ q^{(\varnothing )}_n  \; ; \; n\geq 1 \} \; .
\end{equation}
Then we get by (\ref{intskel}) 
\begin{equation}
\label{excu5}
\bigcap_{n\geq 1}\bigcup_{\substack{u\in \U \\ |u|=n}}C_u = {\rm
  Lf}(T)\backslash \{ q^{(\varnothing )}_n  \; ; \; n\geq 1 \} \; .
\end{equation}
Denote by $\U_{\infty}$ the set of the positive integers valued
sequences. Let $v_{\infty}=(v_{\infty}(n); n\geq 1)$ be in $\U_{\infty}$. Set 
$u_n=(v_{\infty}(1), \ldots , v_{\infty}(n))$ and define the non-empty compact
$C_{v_{\infty }}$ by 
$$ C_{v_{\infty }} = \bigcap_{n\geq 1} C_{u_n} \subset 
{\rm Lf}(T)\backslash \{ q^{(\varnothing )}_n  \;
; \; n\geq 1 \}  .$$ 
Suppose that $C_{v_{\infty }}$ contains two distinct leaves $\sigma$ and
$\sigma'$. There exist $n,n'\geq 1$ such that 
$$ d(\sigma, \sigma \wedge q^{(\varnothing )}_n ) < \frac{1}{3} d(\sigma, \sigma
\wedge \sigma' ) \quad {\rm and } \quad d(\sigma', \sigma' \wedge q^{(\varnothing )}_{n'} ) < \frac{1}{3} d(\sigma', \sigma
\wedge \sigma' ) \; .$$ 
It implies that there exist two distinct words $u,u'\in \U$ such that 
$$ \sigma \in C^o_u \; , \; \sigma' \in C^o_{u'} \quad {\rm and } \quad C^o_u
\cap C^o_{u'}=\emptyset \; ,  $$ 
which rises a contradiction. Consequently $C_{v_{\infty }}$ reduces to a
single point denoted by $\xi(v_{\infty })$. Moreover, deduce from (\ref{excu5}) 
that $\xi$ define a bijective map from $\U_{\infty}$ onto ${\rm
Lf}(T)\backslash \{ q^{(\varnothing )}_n  \; ; \; n\geq 1 \}$. In addition
observe that for any $u=(k_1, \ldots , k_n) \in \U$ with 
$k_1, \ldots , k_n \geq 1$, we get 
$$  \xi^{-1} \left( C_u \cap ({\rm
Lf}(T)\backslash \{ q^{(\varnothing )}_n  \; ; \; n\geq 1 \}) \, \right) =
\{v_{\infty } \in \U_{\infty} \; : \; v_{\infty }(i)=k_i \; , \; 1\leq i\leq n
\} \; .$$ 
It implies that $\xi $ is measurable when $\U_{\infty}$ is equipped with the 
sigma-field generated by the applications $v_{\infty }\rightarrow
v_{\infty }(n)$ , $n\geq 1$, and when 
${\rm Lf}(T)\backslash \{ q^{(\varnothing )}_n ; n\geq 1 \}$ 
is equipped with the trace of the
Borel sigma-field.

  Let ${\bf p}=(p_i;i\geq 1)$ be a probability distribution on the positive 
integers such that $p_i >0$ , $i\geq 1$. Let $V=(\kappa_n ; n\geq 1)$ be a
sequence of i.i.d random variables distributed in accordance with ${\bf p}$. 
Denote by $\mu$ the distribution of $\xi(V)$. Clearly $\mu$ satisfies
(CT3). Let $\sigma=\xi(v_{\infty})$. Observe that 
$$ \mu (\{ \sigma\})= \bP (\xi(V)= \xi(v_{\infty}))= \lim_{n\rightarrow \infty} 
\bP (v_{\infty }(i)=\kappa_i \; ; \; 1\leq i\leq n)=0 \; .$$
Thus, $\mu$ satisfies (CT2).

For any $n\geq 1$ we set 
$u_n=(v_{\infty}(1) ,  \ldots ,  v_{\infty}(n))$, then $\{ \sigma\}= \bigcap
C_{u_n} $, by definition of $\xi (v_{\infty})=\sigma $. It implies that 
the diameter of $C_{u_n} $ goes to zero. 
Thus, for any $\epsilon >0$, there exists $n(\epsilon) $ such
that $ C_{u_{n(\epsilon)}} $ is contained in the open ball 
$B(\sigma , \epsilon)$ centered at $ \sigma$ with radius
$\epsilon$. Consequently, 
$$ \mu (B(\sigma , \epsilon))\geq 
\mu (C_{u_{n(\epsilon)}} )=p_{v_{\infty }(1)}
p_{v_{\infty} (2)} \ldots p_{v_{\infty }(n(\epsilon))} >0\; , $$
which implies (CT1). This completes the proof of the proposition. \cqfd

\vspace{4mm}

Let us consider a continuum tree $(T,d,\rho, \mu)$. We now make 
the connection with an earlier work of Aldous (namely Theorem 15 in 
\cite{Al2}) that
provides a randomized construction of the height 
function of continuum trees. This construction detailed in the proof of 
Theorem 15 in \cite{Al2} can be rephrased as follows: 

\begin{itemize}

\item Let $\Sigma_n$, $n\geq 1$, be an i.i.d. sequence of points in $T$ with 
distribution  $\mu$. Since $(T, d, \rho , \mu)$ is a continuum tree,
then a.s. the $\Sigma_n$'s are distinct leaves and they form a dense subset of $T$.

\item We equip the continuum tree $(T,d,\rho,\mu)$ with the random uniform shuffling 
$\leq_{{\rm Sh}}$ that is assumed to be independent of the sequence $\Sigma_n$ , 
$n\geq 1$. 

\item We set $\Sigma_0 = \rho$. For any $n\geq 0$, we define a random 
number $U_n$ in $[0, 1]$ as follows: 

-- $\; $ We set $U_0=0$; we also assume that $U_1$ is independent of the sequence $\Sigma_n$ , $n\geq 1$, and 
that $U_1$ is uniformly distributed in $[0, 1]$.

-- $\; $ Suppose that $U_1, \ldots , U_n$ have been defined;
  there are two cases.  Either there 
exists a pair $k_1, k_2 \in \{ 0, \ldots , n\}$ such that 
$\Sigma_{n+1}$ is the unique point 
$\sigma \in  \{\Sigma_0 , \dots , \Sigma_{n+1} \}$ such that 
$\Sigma_{k_1 }  <_{{\rm Sh}} \sigma <_{{\rm Sh}} \Sigma_{k_2 }$; in
  that case, pick 
$U_{n+1}$ uniformly at random in the closed interval whose ends are $U_{k_1}$
  and $ U_{k_2}$. Either 
$$\forall  k \in \{ 0, \ldots , n\} \; , \quad \Sigma_{k }  <_{{\rm
      Sh}}\Sigma_{n+1} \; ; $$
in that case, pick 
$U_{n+1}$ uniformly at random in the interval $[\max_{0 \leq k\leq n} U_{k} \, , \,  1]$.
\end{itemize}

\noindent
Now set for any $t\in [0, 1]$, 
$$ f(t)=\limsup_{\epsilon \rightarrow 0} \left\{  
d(\rho, \Sigma_n ) \; , n\geq 0 \; : \;  \; U_n \in [t-\epsilon , t+\epsilon ] \right\} \; .$$
Then Theorem 15 \cite{Al2} implies that a.s. $f$ is a continuous 
function such that $(T_f, d_f)$ and $(T, d)$ are isometric and it is clear 
from Theorem \ref{main} that $f$ is 
the (unique, by Proposition \ref{shuffleconti} and Remark \ref{lamesu}) 
height function associated with the structured tree 
$(T, d, \rho , \leq_{{\rm Sh}}, \mu )$.

    Consequently, all height functions constructed thanks to Theorem 15 in 
\cite{Al2} coincide with the construction given by Theorem 
\ref{main}: in particular, it is the case of the normalized Brownian 
excursion that encodes the Continuum Random Tree; it is also 
the case of the height functions of the 
Inhomogeneous Continuum Random Trees given in \cite{AlMiPi2} and of the height 
functions of the genealogical tree of stable fragmentations in \cite{HaaMi}.

    L\'evy trees introduced by Le Gall and Le Jan in \cite{LGLJ1}
   generalize the Brownian tree. They are constructed via 
the coding  by the so called {\it Height Process} that is a local-time functional 
of a spectrally positive L\'evy process. L\'evy trees can be seen as 
family trees of {\it continuous states branching processes} that have 
been introduced by 
Jirina and Lamperti (see \cite{Ji, La1} and also \cite{Bi2}). The L\'evy trees 
are the scaling limits of the discrete Galton-Watson 
trees (see \cite{Du2, DuLG, DuLG2} for a detailed account on 
that topics). When the underlying branching process a.s. dies out in finite 
time, then the Height Process is continuous with compact support and the L\'evy tree coded by 
this process is a continuum tree ${\bf (T , d, \rho , {\rm m} )}$. Moreover, 
given ${\bf (T , d, \rho , {\rm m} )}$ the order induced by the Height Process 
corresponds to a uniform random shuffling. Consequently, if we fix the structured tree 
${\bf (T , d, \rho , \leq  , {\rm m} )}$ coded by a sample path of the Height 
Process, then the height function given by Theorem \ref{main} coincides with 
the Height Process itself.

  In all these examples of random trees, order does not really 
matter. Let us end the paper with an example of random tree where the role played by
the order is crucial. Let
$X=(X(t), t\geq 0)$ be a L\'evy process without negative jumps and started at 
$X(0)=x >0$. 
We assume that $X$ does not drift to $+\infty$ so that the
stopping time $M$ given by 
$$ M= \inf \{t\geq 0 \; : \; X(t)=0 \}$$
is a.s. finite. Let us set 
$$ {\bf h} (t)=X(M-t) \; , \quad t\in [0, M] \; . $$ 
Then, $ {\bf h} \in \cH$. We denote by $({\bf T , d, \rho , \leq , \mu })$ the
random structured tree coded by ${\bf h}$. When $X$ is a compound Poisson
process with unit drift, then $X$ can be interpreted as the load of a Last-In-First-Out 
M/G/1 queueing system and the underlying tree is given by the following
rule: we say that Client (a) is the child of Client (b) if Client (b) was
currently served when Client (a) arrived in the line (see \cite{LGLJ1, Li} for more details). 
The underlying tree can also be seen as the life-time tree
of a Crump-Mode-Jagers branching process (see \cite{Jag} or \cite{DuLa} for a
connections with L\'evy processes).

   Here we consider the case of a L\'evy process $X$ for which points
are regular and instantaneaous, namely a.s. $$ \forall \epsilon >0
   \; , \quad   \inf_{0 \leq s \leq \epsilon } X(s) <X(0) < \sup_{0 \leq s \leq \epsilon }
X(s) \; . $$
It is equivalent for the L\'evy process to have infinite variation
paths (we refer to the book of Bertoin \cite{Be} Chapter VII 
Corollary 5 for details). By an easy
time-reversal argument, we can show that for any $t_0 >0$,  a.s. we
get 
$$ \forall \epsilon >0 \; , \quad   \inf_{t_0 \leq s \leq t_0 +\epsilon } X(s) < X(t_0) \quad {\rm and } \quad  
\inf_{t_0 -\epsilon  \leq s \leq t_0 } X(s) < X(t_0) \; .$$
This implies that a.s. $\mu$ is a non-atomic measure and that $\mu ({\rm Sk}({\bf T}))=0$. 
Thus, $({\bf T , d, \rho , \mu })$ is a continuum random tree.

Now, fix $({\bf T , d, \rho , \leq , \mu })$ and denote by $\leq_{\rm sh}$ a random
uniform shuffling of $({\bf T , d, \rho , \mu })$. Then, Proposition \ref{shuffleconti} implies that the new height
function ${\bf h}_{\rm Sh}$ associated by Theorem
\ref{main} with $({\bf T , d, \rho , \leq_{\rm sh} , \mu })$ is
continuous. Thus, it is a continuous rearrangement 
of the L\'evy process $X=(X(t); 0\leq  t\leq M)$ coding the same measured compact rooted real tree. Excepted in the
Brownian case, the distribution of ${\bf h}_{\rm Sh}$ does not seem to be simple to characterize.


\begin{thebibliography}{10}

\bibitem{AlMiPi1}
{\sc Aldous, D., Miermont, G., and Pitman, J.}
\newblock {B}rownian bridge asymptotics for random p-mappings.
\newblock {\em Electronic J. Probab. 9\/} (2004), 37--56.

\bibitem{AlMiPi2}
{\sc Aldous, D., Miermont, G., and Pitman, J.}
\newblock The exploration process of inhomogeneous continuum random trees and
  an extension of {J}eulin's local time identity.
\newblock {\em Probab. Th. Related Fields 129\/} (2004), 182--218.

\bibitem{AlMiPi3}
{\sc Aldous, D., Miermont, G., and Pitman, J.}
\newblock Weak convergence of random p-mappings and the exploration process of
  inhomogeneous continuum random trees.
\newblock {\em Probab. Th. Related Fields 133\/} (2005), 1--17.

\bibitem{Al1}
{\sc Aldous, D.~J.}
\newblock The continuum random tree {I}.
\newblock {\em Ann. Probab. 19\/} (1991), 1--28.

\bibitem{Al2}
{\sc Aldous, D.~J.}
\newblock The continuum random tree {III}.
\newblock {\em Ann. Probab. 21\/} (1993), 248--289.

\bibitem{Be}
{\sc Bertoin, J.}
\newblock {\em L\'evy Processes}.
\newblock Cambridge Univ. Press, 1996.

\bibitem{Bi2}
{\sc Bingham, N.~H.}
\newblock Continuous branching processes and spectral positivity.
\newblock {\em Stochastic Process. Appl. 4\/} (1976), 217--242.

\bibitem{Bun}
{\sc Buneman, P.}
\newblock A note on the metric properties of trees.
\newblock {\em J. Combinatorial Theory Ser. B 17\/} (1974), 48--50.

\bibitem{BuBu}
{\sc Burago, D.~Burago, Y., and Ivanov, S.}
\newblock {\em A Course in Metric Geometry}, vol.~33.
\newblock AMS, Boston, 2001.

\bibitem{Chis}
{\sc Chiswell, I.}
\newblock {\em Introduction to {$\Lambda$}-trees}.
\newblock World Scientific Publishing Co., Inc, River Edge, 2001.

\bibitem{Croy}
{\sc Croydon, D.}
\newblock Measure and heat kernel estimates for the continuum random tree.
\newblock {\em preprint\/} (2005).

\bibitem{Dress84}
{\sc Dress, A.}
\newblock Trees, tight extensions of metric spaces, and the cohomological
  dimension of certain groups: A note on combinatorial properties of metric
  spaces.
\newblock {\em Adv. Math. 53\/} (1984), 321--402.

\bibitem{DMT96}
{\sc Dress, A., Moulton, V., and Terhalle, W.}
\newblock {T}-theory: an overview.
\newblock {\em European J. Combin. 17\/} (1996), 161--175.

\bibitem{DT96}
{\sc Dress, A., and Terhalle, W.}
\newblock The real tree.
\newblock {\em Adv. Math. 120\/} (1996), 283--301.

\bibitem{Du2}
{\sc Duquesne, T.}
\newblock A limit theorem for the contour process of conditioned
  {G}alton-{W}atson trees.
\newblock {\em Ann. Probab. 31}, 2 (2003), 996--1027.

\bibitem{DuLa}
{\sc Duquesne, T., and Lambert, A.}
\newblock Work in progress.
\newblock {\em -\/} (2005).

\bibitem{DuLG}
{\sc Duquesne, T., and Le~Gall, J.-F.}
\newblock {\em Random Trees, {L}\'evy Processes and Spatial Branching
  Processes}.
\newblock Ast\'erisque no 281, 2002.

\bibitem{DuLG2}
{\sc Duquesne, T., and Le~Gall, J.-F.}
\newblock Probabilistic and fractal aspects of {L}\'evy trees.
\newblock {\em To appear in Probab. Theorey and Rel. Fields\/} (2004).

\bibitem{DuLG3}
{\sc Duquesne, T., and Le~Gall, J.-F.}
\newblock The {H}ausdorff measure of stable trees.
\newblock {\em preprint\/} (2005).

\bibitem{DuWi1}
{\sc Duquesne, T., and Winkel, M.}
\newblock Growth of {L}{\'e}vy trees.
\newblock {\em preprint\/} (2005).

\bibitem{Ev00}
{\sc Evans, S.}
\newblock Snakes and spiders: {B}rownian motion on real trees.
\newblock {\em Probab. Theory Related Fields 117}, 3 (2000), 361--386.

\bibitem{EvPitWin}
{\sc Evans, S., Pitman, J., and Winter, A.}
\newblock Rayleigh processes, real trees, and root growth with re-grafting.
\newblock {\em To appear in Probab. Th. Rel. Fields\/} (2005).

\bibitem{EvWin}
{\sc Evans, S., and Winter, A.}
\newblock Subtree prune and re-graft: a reversible real tree valued {M}arkov
  process.
\newblock {\em preprint\/} (2005).

\bibitem{Fels}
{\sc Felsenstein, J.}
\newblock {\em Inferring Phylogenies}.
\newblock Sinauer Associates, Sunderland, Massachusett, 2003.

\bibitem{HaaMi}
{\sc Haas, B., and Miermont, G.}
\newblock The genealogy of self-similar fragmentations with negative index as a
  continuum random tree.
\newblock {\em Electr. J. Probab. 9\/} (2004), 57--97.

\bibitem{LyoHam}
{\sc Hambly, B., and Lyons, T.}
\newblock Uniqueness for the signature of a path of bounded variation and
  continuous analogues for the free group.
\newblock {\em Preprint\/} (2004).

\bibitem{Jag}
{\sc Jagers, P.}
\newblock General branching processes as {M}arkov fields.
\newblock {\em Stoch. Proc. Appl. 32\/} (1989), 213--224.

\bibitem{Ji}
{\sc Jirina, M.}
\newblock Stochastic branching processes with continous state-space.
\newblock {\em Czech. Math. J. 8\/} (1958), 292--313.

\bibitem{Krebs}
{\sc Krebs, W.}
\newblock Brownian motion on the continuum tree.
\newblock {\em Probab. Theory Rel. Fields 101}, 3 (1995), 421--433.

\bibitem{La1}
{\sc Lamperti, J.}
\newblock The limit of a sequence of branching processes.
\newblock {\em Z. Wahrsch. Verw. Gebiete 7\/} (1967), 271--288.

\bibitem{LGman}
{\sc Le~Gall, J.-F.}
\newblock 2005.
\newblock Manuscript notes.

\bibitem{LG1}
{\sc Le~Gall, J.-F.}
\newblock Brownian excursions, trees and measure-valued branching processes.
\newblock {\em Ann. Probab. 19\/} (1991), 1399--1439.

\bibitem{LG93}
{\sc Le~Gall, J.-F.}
\newblock A class of path-valued {M}arkov processes and its applications to
  superprocesses.
\newblock {\em Prob. Th. Rel. Fields 95\/} (1993), 25--46.

\bibitem{LG2}
{\sc Le~Gall, J.-F.}
\newblock The uniform random tree in a {B}rownian excursion.
\newblock {\em Probab. Theory and Related Fields 96\/} (1993), 369--383.

\bibitem{LGLJ1}
{\sc Le~Gall, J.-F., and Le~Jan, Y.}
\newblock Branching processes in {L}\'evy processes: the exploration process.
\newblock {\em Ann. Probab. 26-1\/} (1998), 213--252.

\bibitem{LGLJ2}
{\sc Le~Gall, J.-F., and Le~Jan, Y.}
\newblock Branching processes in {L}\'evy processes: {L}aplace functionals of
  snakes and superprocesses.
\newblock {\em Ann. Probab. 26\/} (1999), 1407--1432.

\bibitem{Li}
{\sc Limic, V.}
\newblock A {L}{I}{F}{O} queue in heavy traffic.
\newblock {\em Ann. Appl. Probab. 11\/} (2001), 301--331.

\bibitem{MayOv}
{\sc Mayer, J., and Oversteegen, L.}
\newblock A topological characterization of {$\R $}-trees.
\newblock {\em Trans. Amer. Math. Soc. 320\/} (1990), 395--415.

\bibitem{Mier03}
{\sc Miermont, G.}
\newblock Self-similar fragmentations derived from the stable tree {I}:
  splitting at heights.
\newblock {\em Probab. Theory Relat. Fields 127}, 3 (2003), 423--454.

\bibitem{Mier05}
{\sc Miermont, G.}
\newblock Self-similar fragmentations derived from the stable tree {II}:
  splitting at nodes.
\newblock {\em Probab. Theory Relat. Fields 131}, 3 (2005), 341--375.

\bibitem{Ne}
{\sc Neveu, J.}
\newblock Arbres et processus de {G}alton-{W}atson.
\newblock {\em Ann. Inst. H. Poincar\'e 26\/} (1986), 199--207.

\bibitem{SemSte}
{\sc Semple, C., and Steel, M.}
\newblock {\em Phylogenetics}, vol.~24 of {\em Oxford Lecture Series in
  Mathematics and its Applications}.
\newblock Oxford University Press, Oxford, 2003.

\bibitem{SimPer}
{\sc Sim{\~o}es~Pereira, J. M.~S.}
\newblock A note on the tree realizability of a distance matrix.
\newblock {\em J. Combinatorial Theory 6\/} (1969), 303--310.

\bibitem{Weil}
{\sc Weill, M.}
\newblock Regenerative real trees.
\newblock {\em preprint\/} (2005).

\bibitem{Zar}
{\sc Zareckii, K.~A.}
\newblock Constructing a tree on the basis of a set of distances between the
  hanging vertices.
\newblock {\em Uspehi Mat. Nauk 20}, 6 (1965), 90--92.

\end{thebibliography}
\end{document}